\RequirePackage{fix-cm}
\documentclass[smallextended]{svjour3}       % onecolumn (second format)
\smartqed  % flush right qed marks, e.g. at end of proof
\usepackage{graphicx,amsmath,amsfonts}
\usepackage{amsmath}
\usepackage{amssymb}
\usepackage{mathtools}
\usepackage{float}
\usepackage{bbm}
\usepackage{longtable}
\usepackage{graphicx}
\usepackage{hyperref}
\usepackage{pifont} 
\usepackage{stmaryrd}

\usepackage{subfig}
\usepackage{wrapfig}
\usepackage{mathrsfs}

\usepackage{xcolor}
\usepackage[normalem]{ulem}

\DeclarePairedDelimiter{\ceil}{\lceil}{\rceil}

\def\Bila{\mathbf{B}}

\def\S{\mathbf{S}}

\def\SS{\mathcal{S}}

\def\P{\mathcal{P}}

\def\D{\mathcal{D}}

\def\K{\mathcal{K}}

\def\R{\mathbb{R}}

\def\A{\mathcal{A}}

\def\I{\mathcal{I}}

\def\J{\mathcal{J}}

\def\F{\mathcal{F}}

\def\bu{\textbf{u}}

\newcommand {\nor} [1]{\parallel #1 \parallel}

\newcommand{\eequ}{\end{equation}}

\newcommand{\bequ}{\begin{equation}}

\newcommand{\eequd}{\end{eqnarray*}}

\newcommand{\bequd}{\begin{eqnarray*}}

\definecolor{Darkgreen}{rgb}{0,0.5,0}

\newcommand{\re}[1]{\textcolor{black}{#1}}

\newcommand{\dt}[1]{\textcolor{black}{#1}}

\begin{document}

\title{Re-polarisation of macrophages within a multi-scale moving boundary tumour invasion model}

\titlerunning{Re-polarisation of macrophages within cancer invasion}        % if too long for running head

\author{Szabolcs Suveges	\and
        Raluca Eftimie		\and
        Dumitru Trucu
}


\institute{Szabolcs Suveges \at
              Division of Mathematics\\
              University of Dundee\\
              Dundee, DD1 4HN, UK\\
              \email{ssuveges@dundee.ac.uk}           %  \\
           \and
           Raluca Eftimie \at
           Laboratoire Mathématiques de Besançon, UMR - CNRS 6623,\\
           Université de Bourgogne Franche-Comté,\\
           Besançon, 25000, France\\
            \email{raluca.eftimie@univ-fcomte.fr}    
            \and
            Dumitru Trucu \at
            Division of Mathematics\\
           University of Dundee\\
           Dundee, DD1 4HN, UK\\
            \email{trucu@maths.dundee.ac.uk}    
}

\date{Received: date / Accepted: date}
\maketitle

\begin{abstract}
\re{Cancer invasion of the surrounding tissue is a multiscale process that involves not only tumour cells but also other immune cells in the environment, such as the tumour-associated macrophages (TAMs). The heterogeneity of these immune cells, with the two extremes being the pro-inflammatory and anti-tumour M1 cells, and the anti-inflammatory and pro-tumour M2 cells, has a significant impact on cancer invasion as these cell interact in different ways with the tumour cells and with the ExtraCellular Matrix (ECM). Experimental studies have shown that cancer cells co-migrate with TAMs, but the impact of these different TAM sub-populations (which can change their phenotype and re-polarise depending on the microenvironment) on this co-migration is not fully understood.}
In this study, we extend a previous multi-scale moving boundary mathematical model, by introducing the M1-like macrophages alongside with their exerted multi-scale effects on the tumour invasion process. With the help of this model we investigate numerically the impact of re-polarising the M2 TAMs into the anti-tumoral M1 phenotype and how such a strategy affects the overall tumour progression. \re{In particular, we investigate numerically whether the M2$\to$M1 re-polarisation could depend on time and/or space, and what would be the macroscopic effects of this spatial- and temporal-dependent re-polarisation on tumour invasion. }

\end{abstract}

\keywords{Cancer invasion, Macrophages, Macrophage re-polarisation, Multi-scale modelling, Cell adhesions, WENO schemes, Convolution}

%%%%%%%%%%%%%%%%%%%%%%%%%

\section{Introduction}

\re{The last few decades have seen a shift in the focus of cancer research: from a research that was focused on individual tumour cells to a research that is now focused on tumour microenvironment (TME) and the interactions between different types of cells inside the TME \cite{Henke2020}. The TME is formed of tumour's vasculature, connective tissue, infiltrating immune cells and the extracellular matrix (ECM). In recent years the ECM has received considerable attention, due to its role in cancer evolution and response to therapies \cite{Henke2020}}.
%As opposed to normal cells, it is well known that cancer cells are capable of tissue invasion and metastasis\cite{Hanahan2000,Hanahan2011} as well as they can establish non-local adhesive bonds with other neighbouring cells and with the extracellular matrix (ECM). 

The ECM is a complex network of macromolecules (such as fibrous proteins, water and minerals), which is an essential part of any healthy tissue \cite{Filipe2018_unexploredECM_cancer}. To maintain its functionality, the ECM is subject to continuous remodelling (via synthesis and degradation). However, this carefully orchestrated process is disturbed by the progressing tumour, resulting in drastic changes in the structure of the ECM. Moreover, in the tumour micro-environment, the excessive degradation of the ECM is related to the over-secretion of several enzymes capable of degrading the ECM, for instance, the matrix metalloproteinases (MMPs). Such enzymes are not explicitly tied to cancer cells, but instead, several other stromal cells secrete them as well \cite{MadsenBugge2015_SourceMDEcancer}.

One of these stromal cell populations is the macrophages, which can form up to $50 \%$ of the tumour mass \cite{Kelly_1988,Vinogradov2014_Macroph50PercentCancer}. Macrophages are a heterogeneous set of differentiated immune cells that are polarised in accordance with the surrounding stimuli in order to exhibit different properties and functions. 
A specific class of them that is present in the tumour micro-environment is called the tumour-associated macrophages (TAMs) that are in general correlated with poor prognosis in various type of cancers \cite{Caro2015,Cassetta2016,Hallam2011,Li2019,Qiu2018,Wang2016,Wu2016,Yu2019,Yuan2017,Zhao2017}. These activated macrophages are often categorised as the pro-inflammatory M1-like macrophages and the anti-inflammatory M2-like macrophages. Hence, while the former one expresses anti-tumoral functions, the latter one has pro-tumoral properties. Due to their fundamentally different effect on the tumour, the ratio between the M1 and M2 TAMs can significantly affect the prognosis \cite{Cui2013,Pantano2013,Zhang2012,Yin2017,Yuan2017}. Moreover, due to the plasticity of the macrophages, their phenotype can be changed in response to the environment \cite{Becker_2013,Rocha_2014} and so targeting the ratio between the two extreme phenotypes could be an important tool for fighting against tumour development. To this end, \re{in this study} we explore the re-polarisation of M2 TAMs into the M1 phenotype in tumours with different characteristics, paying special attention to these two populations near the tumour interface. Furthermore, since tumours are only detectable above a certain size, we investigate not only the spatial dependency of macrophages re-polarisation, but also the temporal dependency of this re-polarisation, and how they affect the rate of tumour growth.

Nutrients (e.g. oxygen and glucose) are essential for any cells to live and function properly, including cancer cells and macrophages. In order to supply these vital nutrients to every tissue of the body, nutrients are extravasated from the blood flow and diffuse through the ECM. However, the progressing tumour destroys significant parts of the vasculature, excessively influencing the blood flow and therefore, the level of supplied nutrients. Thus, as the tumour progresses, most of these nutrients are supplied through the vasculature, located in the peritumoral regions \cite{Zhang2019}. This creates areas inside the tumour where the level of nutrients become low and so such areas first become hypoxic and then necrotic. Consequently, cells situated in necrotic regions are being deprived of nutrients, leading to cell death and ultimately giving rise to a necrotic core. \re{Hypoxic conditions can also modify the polarisation of macrophages and influence the malignant behaviour of some cancer cells \cite{Ke2019_HypoxiaMacrophPolarise}}.

Over the last decade, there have been substantial advances in mathematical modelling to understand the dynamics of the cancer cell populations as well as their interactions with their surroundings \cite{Anderson_et_al_2000,Anderson2009,Anderson_2005,Chaplain_Lolas_2005,Chaplain2006,Deakin_Chaplain_2013,Deisboeck_et_al_2011,Domschke_et_al_2014,Kiran_2009,KnutsdottirPalssonKeshet2014,Macklin_et_al_2009,MahlbacherLowengrubFriesboes2018,Shuttleworth_2019,Suveges_2020,szymanska_08,Dumitru_et_al_2013,Xu_2016}. Although these models mainly \re{focused on} the interplay between the cancer cells and the neighbouring ECM, the importance of the macrophages during tumour development cannot be overlooked and as a consequence, some mathematical models started to investigate the role of macrophages on the overall tumour progression \cite{Dallon1999,KnutsdottirPalssonKeshet2014,MahlbacherLowengrubFriesboes2018,McDougall2006,Owen2004,Owen1997,Owen1998,Suveges_2020,Webb2007}. Initially, these models \cite{Owen1997,Owen1998,Webb2007} focused only on the anti-tumour role of the macrophages, and only later models started to explore also their pro-tumour role \cite{KnutsdottirPalssonKeshet2014,MahlbacherLowengrubFriesboes2018}. Furthermore, most of these initial mathematical models investigated tumour progression only at one spatio-temporal scale \cite{Anderson_et_al_2000,Chaplain_Lolas_2005,Chaplain2006,szymanska_08}, but the relevance of various biological processes occurring on different scales cannot be ignored. Consequently, recent models started to capture these interactions \cite{Deisboeck_et_al_2011,Domschke_et_al_2014,Peng2016,Shuttleworth_2019,Suveges_2020,Dumitru_et_al_2013} in order to investigate the different aspects of tumour development in a multi-scale fashion. However, due the novelty of these approaches, they have not been extended to capture also the complex roles of macrophages in the tumour development.

In this work, we further extend an existing multi-scale moving boundary modelling platform \cite{Shuttleworth_2019,Suveges_2020,Dumitru_et_al_2013}. We build upon a mathematical model \cite{Suveges_2020} that accounts not only for the proteolytic processes occurring at the leading edge of the tumour but also for the fibre and non-fibre components of ECM, as well as for the presence of M2-like macrophages. \re{In the new model introduced in this study} we capture the contribution of M1 TAMs to the overall tumour dynamics by considering their presence both on the macro-scale and on the micro-scales. Specifically, we extend the macro-scale dynamics by introducing an equation describing the spatio-temporal evolution of the M1 TAMs population, as well as by considering several of its effects on the rest of the tumour dynamics. Besides, we take into account their contribution to the proteolytic micro-scale processes as well as their impact on the rearrangement of the ECM fibre micro-constituents. Hence, M1 TAMs not only influence the directionality of the moving tumour boundary but also the spatial bias of the overall ECM structure to withstand incoming cell fluxes.
\re{Moreover, in this study} we take into consideration the effects of the depleted nutrients on different cell functions. This enables us to introduce regulations on cell proliferation, on cell death and on the macrophage polarisation based on the available nutrients which are supplied through the peritumoral blood vessels. Ultimately, with the extended model that we propose in this work, we aim to investigate whether an effective re-polarisation strategy is spatial/temporal dependent or not within a fibrous tissue environment.

The structure of the paper is as follows. In Section~\ref{section:Macro_Scale_Model} we present the extended macro-scale dynamics. Then, in Section~\ref{section:Micro_Scale_Model} we outline the two micro-scale processes. In Section~\ref{section:Numerical_Approach}, we describe the convolution-driven numerical schemes that we use to solve the multi-scale moving boundary model. We present our numerical simulations in Section~\ref{section:Numerical_Results}, by focusing especially on the spatial and temporal dependency of the re-polarisation. At last, we summarise and discuss the results in Section~\ref{section:Conclusions}.

%--------------------------------------------------------------
%	  						 Model
%--------------------------------------------------------------
\section{Multi-Scale Modelling of the Tumour Dynamics}\label{section:Model}
%--------------------------------------------------------------
%	  				Intro for the model
%--------------------------------------------------------------
\dt{Building on} the multi-scale moving boundary framework initially introduced in \cite{Dumitru_et_al_2013} and later expanded in \cite{Shuttleworth_2019,Suveges_2020}, \dt{in this work}, we \dt{explore} not only the M2-like macrophages \cite{Suveges_2020} but also the M1 phenotype, \dt{by assessing the impact that these bring within the interlinked tissue-scale (macro-scale) and cell-scale (micro-scale) tumour dynamics.} Moreover, the biological context of the \dt{cacer} macro-dynamics is \dt{further} broadened by \dt{considering the presence of the} nutrients, \dt{such as} include glucose and oxygen, \dt{which are key constituents of the tumour microenvirnment and play an impotan role within overall tumour progression}. 
%Therefore, in this work, we further extend the multi-scale moving boundary framework by considering the dynamics of M1 TAMs population as well as the dynamics of nutrients and their effects on the overall tumour progression.

\subsection{Macro-Scale Dynamics}\label{section:Macro_Scale_Model}
%--------------------------------------------------------------
%				Intro for the Macro-Scale
%--------------------------------------------------------------
\dt{As this work extends the modelling} framework \dt{introduced in} \cite{Shuttleworth_2019,Suveges_2020,Dumitru_et_al_2013}, we start this section by introducing some of its \dt{key} features. Thus, on the macro-scale we focus on the expanding tumour region $\Omega(t)$ that progresses within \dt{a} maximal tissue cube $Y \dt{\subset} \R^{d}$, for $d = 2,3$ and over the time interval $[0, T]$ (\emph{i.e.}, $\Omega(t) \subset Y, \forall t\in [0,T]$). \dt{In this context}, at any macro-scale spatio-temporal point $(x,t) \in \Omega(t) \times [0,T]$, we consider \dt{a mixed cell population consisting of distributions of}: $(a)$ cancer cell\dt{s}  $c(x,t)$; $(b)$ \dt{M1-like macrophages, $M_{1}(x,t)$, briefly addressed here as  M1 TAM};  and $(c)$ \dt{M2-like macrophages, $M_{2}(x,t)$, which are briefly referred to as M2 TAM}. \dt{This mixture of cancer cells and macrophages exercise their naturally multiscale dynamics within an extracellular matrix, which, as in \cite{Shuttleworth_2019,Shuttleworth2020b,Shuttleworth2020a}, is regarded as consisting of two major phases, namely a fibrous and a non-fibrous one.}   Specifically, on the one hand, we have the fibre ECM phase, accounting for all major fibrous proteins (such as collagen and fibronectin), \dt{whose micro-scale structure enables a }spatial bias \dt{for withstanding incoming spatial cell fluxes, inducing this way an intrinsic ECM fibres spatial orientation \cite{Shuttleworth_2019,Shuttleworth2020b,Shuttleworth2020a,Suveges_2020}. Therefore, the spatio-temporal distribution of the oriented ECM fibres} at the macro-scale point $(x,t)$ \dt{is described by a vector field} $\theta_{f}(x,t)$, \dt{where $\theta_{f}(\cdot,t):\R^{d}\to\R^{d}$, with its Euclidean norm $F(x,t):=\nor{\!\!\theta_{f}(x,t)\!\!}_{_{2}}$ representing the amount of fibres at $(x,t)$}. \dt{Then, on the other hand,} besides these fibrous proteins, the ECM also contains many other components such as non-fibrous proteins (for instance amyloid fibrils), enzymes, polysaccharides and extracellular $Ca^{2+}$ ions. Hence, in the second ECM phase, we bundle together these constituents and refer to it as the non-fibre ECM phase, \dt{and its distribution at each $(x,t)\in \Omega(t) \times [0,T]$ is denoted} by $l(x,t)$. \dt{Therefore, }for compactness, we denote the global five-dimensional tumour vector by $\bu$ that is given by
\begin{equation}
	\bu := (c(x,t), M_{1}(x,t), M_{2}(x,t), l(x,t), F(x,t))^{\intercal},
	\label{5D_Tumour_Vector}
\end{equation}
as well as denoting the total space occupied at $(x,t)$ by $\rho(\bu)$ and define it as
\begin{equation}
	\rho(\bu) := c(x,t) + F(x,t) + l(x,t) + M_{1}(x,t) + M_{2}(x,t),
	\label{Total_Space_Occupied}
\end{equation}
for all $(x,t) \in \Omega(t) \times [0,T]$. \dt{Furthermore}, the last component of the macro-scale dynamics is the nutrients density $\sigma(x,t)$, whose level \dt{within the tumour microenvironment} is depleted by the \dt{invading cancer.}

\dt{Finally}, the macrophages \dt{are considered }to infiltrate the tumour through the \emph{outer-boundary} \cite{Suveges_2020}, which is denoted by $\partial \Omega_{o}(t) \subset \partial\Omega(t)$ {(see Fig. \ref{fig:Re_Polarisation_Domain}), and is} defined in \ref{outerBoundarySet}.

\subsubsection{Nutrient}\label{Nutrients}
\re{The growth of tumours leads to the inability of normal blood vessels to provide enough oxygen and nutrients for the tumour cells \cite{Siemann2015_ModulationTumVasculatOxygen}. Moreover, during tumour growth, some of the pre-existing blood vessels disintegrate, are obstructed or are compressed, also leading to a reduction in the available level of oxygen, nutrients, growth factors, etc \cite{Vaupel1989_OxygenNutrientsTum}. Consequently, the tumour initiates its own vasculature network developed in the proximity of tumour boundary, which helps bring in various nutrients that help tumour growth.}

\re{To incorporate this aspect into our mathematical model, we consider} the overall influx of nutrients through the outer tumour boundary $\partial \Omega_{o}(t)$ by using Dirichlet boundary conditions. Although different cell types uptake the supplied nutrients at different rates, for simplicity here we assume that all present cell populations (cancer cells, M1 and M2 TAMs) uptake nutrients at the same constant rate $d_{\sigma} > 0$. The spatial transport of the nutrients is modelled by a diffusion \re{term} with constant coefficient $D_{\sigma} > 0$. Since this diffusion occurs more rapidly than cell diffusion (i.e., cell random walk), we use a quasi steady state reaction-diffusion equation (similar to the one for instance in \cite{Macklin2009}) for the nutrients $\sigma(x,t)$:
\begin{subequations}
	\begin{align}
		0 = & D_{\sigma}\Delta \sigma - d_{\sigma}(c + M_{1} + M_{2}) \sigma, \label{Nutrient_Equation_eq}\\[0.2cm]
		\sigma(x,t) = & \sigma_{nor}, \qquad \forall x \in \partial \Omega_{o}(t), \forall t \in [0,T].\label{Nutrient_Equation_boundary}
	\end{align}
	\label{Nutrient_Equation}
\end{subequations}
Here, $\sigma_{nor}$ is the normal level of nutrients along the tumour interface, which is considered here to be a constant.

Since cells require nutrients to function properly, here we introduce four smooth bounded \emph{effect-functions} that we use to model the effects of the nutrients on the different cell functions. To construct these effect-functions, let us \re{define} \re{three critical nutrient levels:
\begin{itemize}
\item First, we define the necrotic threshold $\sigma_{n} > 0$: if $\sigma \leq \sigma_{n}$ this leads to necrotic tumour cell death in that area \cite{Lee2017_TumNecrosis-TumProgress}. 
\item Then, we define the nutrient level sufficient for cells to function properly: $\sigma_{p} > 0$.
\item  Finally, we define the normal level of nutrients: $\sigma_{nor}$. We also use this threshold value in the Dirichlet boundary condition in \eqref{Nutrient_Equation}. 
 \end{itemize}
 }
 \noindent Hence, we have the following relationship between these three values, namely: 
 $\sigma_{nor} > \sigma_{p} > \sigma_{n}.$

\begin{figure*}
\centering
  \includegraphics[width=0.9\textwidth]{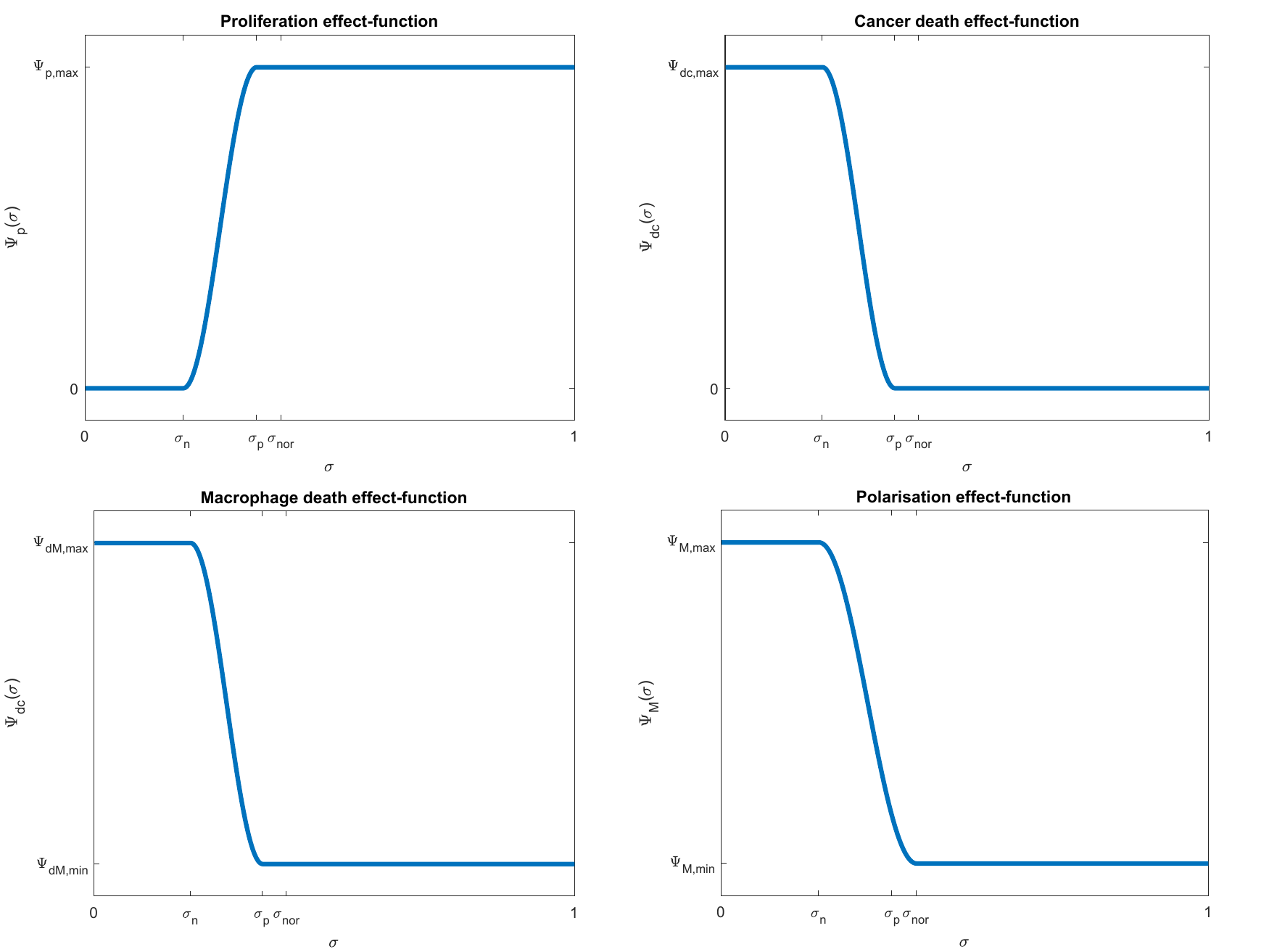}
\caption{Plots of the four nutrients effect-functions $\Psi_{p}(\sigma)$, $\Psi_{dc}(\sigma)$, $\Psi_{dM}(\sigma)$ and $\Psi_{M}(\sigma)$ defined in \eqref{Nutrient_Functions_Proliferation}, \eqref{Nutrient_Functions_Death},  and \eqref{Nutrient_Functions_Polarisation}, respectively.}
\label{fig:Nutrient_Functions}
\end{figure*}

Starting with the effect of nutrients on cell proliferation, we first assume that this effect is identical for cancer cells and both macrophage populations. Hence, we consider a maximal proliferation enhancement rate $\Psi_{p,max} > 0$ which corresponds to nutrient levels $\sigma \geq \sigma_{p}$. Also, we assume no proliferation in the necrotic regions. Thus, we define the proliferation effect-function as follows:
\begin{equation}
	\Psi_{p}(\sigma) :=
	\begin{cases}
		0 & \text{if } \sigma \leq \sigma_{n}, \\
		\Psi_{p,max} & \text{if } \sigma \geq \sigma_{p}, \\
		\Phi (\sigma, \Psi_{p,max}, 0, \sigma_{p} - \sigma_{n}) & \text{otherwise}, \\
	\end{cases}
	\label{Nutrient_Functions_Proliferation}
\end{equation}
where $\Phi (\sigma, \cdot, \cdot, \cdot)$ is a generic cosine function that describes a smooth transition between the two extrema, \emph{i.e.,} for any level $\sigma_{n} < \sigma < \sigma_{p}$. Thus, $\Phi (\sigma, \cdot, \cdot, \cdot)$ is defined by
\begin{equation}
	\Phi (\sigma, \Phi_{max}, \Phi_{min}, \Phi_{L}) := \dfrac{\Phi_{max} - \Phi_{min}}{2} \Bigg [ \cos \bigg( \dfrac{\pi (\sigma - \sigma_{n} - \Phi_{L})}{\sigma_{p} - \sigma_{n}} \bigg) + 1 \Bigg] + \Phi_{min},
	\label{Nutrient_Functions_Transition}
\end{equation}
where $\Phi_{min}$ is the minimum, $\Phi_{max}$ is the maximum and $\Phi_{L}$ controls the phase shift of the cosine function. We illustrate the proliferation effect-function, defined in \eqref{Nutrient_Functions_Proliferation} in Fig.~\ref{fig:Nutrient_Functions}.

As opposed to cell proliferation, we distinguish between the death of cancer cells and death of macrophages, since cancer cells resist death \cite{Hanahan2000,Hanahan2011} while macrophages do not. Therefore, we consider a maximal enhancement death rate $\Psi_{d,max} > 0$ in necrotic regions for both populations, and while we \re{assume} no death for cancer cells, we consider a minimal level of macrophage enhancement death rate $\Psi_{dM,min} > 0$ in regions where the nutrient level is $\sigma \geq \sigma_{p}$. Hence, cancer cell death and macrophage death effect-functions are defined as (using the transition function \eqref{Nutrient_Functions_Transition}):
\begin{subequations}
\begin{align}
	\Psi_{dc}(\sigma) := &
	\begin{cases}
		\Psi_{d,max} & \text{if } \sigma \leq \sigma_{n}, \\
		0 & \text{if } \sigma \geq \sigma_{p}, \\
		\Phi (\sigma, \Psi_{dc,max}, 0, 0) & \text{otherwise}, \\
	\end{cases} \label{Nutrient_Functions_Death_Cancer}
	\\
	\Psi_{dM}(\sigma) := &
	\begin{cases}
		\Psi_{d,max} & \text{if } \sigma \leq \sigma_{n}, \\
		\Psi_{dM,min} & \text{if } \sigma \geq \sigma_{p}, \\
		\Phi (\sigma, \Psi_{dM,max}, \Psi_{dM,min}, 0) & \text{otherwise}. \\
	\end{cases} \label{Nutrient_Functions_Death_Macrophages}
\end{align}
\label{Nutrient_Functions_Death}
\end{subequations}
These death effect-functions are illustrated in Fig.~\ref{fig:Nutrient_Functions}.

Finally, it was shown \re{experimentally} \cite{Tan_2016} that \re{the} nutrient level affects the macrophages polarisation rate, and so here we consider also a polarisation effect-function. Thus, we consider a maximal enhancement rate $\Psi_{M,max} > 0$ in necrotic regions and a minimal enhancement rate $\Psi_{M,min} > 0$ in regions with normal nutrient levels. Using the smooth transition function \eqref{Nutrient_Functions_Transition}, we define the polarisation effect-function $\Psi_{M}(\sigma)$ by
\begin{equation}
	\Psi_{M}(\sigma) :=
	\begin{cases}
		\Psi_{M,max} & \text{if } \sigma \leq \sigma_{n}, \\
		\Psi_{M,min} & \text{if } \sigma \geq \sigma_{nor}, \\
		\Phi (\sigma, \Psi_{M,max}, \Psi_{M,min}, 0) & \text{otherwise}. \\
	\end{cases}
	\label{Nutrient_Functions_Polarisation}
\end{equation}
This function is illustrated in Fig.~\ref{fig:Nutrient_Functions}.
%------------ Figure for the nutrient functions
%------------

%--------------------------------------------------------------
%		 			   M1 and M2 TAMs
%--------------------------------------------------------------
\subsubsection{Dynamics of both M1 and M2 TAMs}\label{Section_M1M2}
\dt{Focusing now on} the macrophage population, specifically to the two extreme phenotypes, M1 and M2 TAMs, in this paper we focus on the \dt{dynamics of the} two TAMs populations exclusively inside the tumour domain $\Omega(t)$, the \dt{evolution of the} macrophages \dt{distribution} in the surrounding tumour stroma \dt{being} beyond the scope of this current work. Hence, since macrophages are recruited to tumour sites as an immune response through the peritumoral vasculature, here this influx is represented by a source of M1 TAMs  \dt{that is localised} along the outer tumour boundary $\partial \Omega_{o}(t)$, \dt{which is enabled} by the immediate activation of macrophages into M1 TAMs as they enter the tumour \cite{Suveges_2020}. For simplicity, we assume that both profile and maximal magnitude $M_{0} > 0$ of this source are identical along the tumour interface, and so this influx term is given by
\begin{equation}
	M_{I} := M_{0}(\chi_{\partial \Omega_{o}(t)} \ast \psi_{\rho})(x).
	\label{M_influx}
\end{equation}
Here, $\psi_{\rho}$ is the standard mollifier defined in Appendix \ref{Standard_Mollifier} with appropriately chosen range $\rho > 0$, $``\ast"$ is the convolution operator \cite{Damelin_2011} and $\chi_{\partial \Omega_{o}(t)}$ is the characteristic function of the outer boundary $\partial \Omega_{o}(t)$.

As in \cite{Suveges_2020}, we recognise that the stiffness of the ECM plays a role in the proliferation of the macrophages \cite{Hayenga_2015} and that according to  several biological studies \cite{Cassetta_2019,Chitu_2011,Jenkins_2011}, cancer cells trigger the proliferation of macrophages by producing survival and proliferation factors. Hence, we denote by $\mu_{MF} > 0$ the \dt{macrophages} proliferation enhancement rate \dt{due to} the fibres \dt{ while} the proliferation effect-function $\Psi_{p}(\sigma)$ defined in \eqref{Nutrient_Functions_Proliferation} account\dt{s} for the effect of nutrients on the macrophage proliferation. For simplicity, both the fibre enhancement rate and the effect-function are assumed to \dt{remain unchanged} for both phenotype. Thus, we formulate the proliferation laws for the M1 and M2 TAMs as
\begin{subequations}
	\begin{align}
		P_{M_{1}}(\bu) := & \mu_{M}\Psi_{p}(\sigma)(1 + \mu_{MF}F)M_{1}c(1 - \rho(\bu))^{+}, \label{Macrophage_Proliferation_M1}\\[0.2cm]
		P_{M_{2}}(\bu) := & \mu_{M}\Psi_{p}(\sigma)(1 + \mu_{MF}F)M_{2}c(1 - \rho(\bu))^{+}, \label{Macrophage_Proliferation_M2}
	\end{align}
	\label{Macrophage_Proliferation}
\end{subequations}
respectively. In \eqref{Macrophage_Proliferation}, $\mu_{M}$ is the positive baseline proliferation rate \dt{ while the term} $(1 - \rho(\bu))^{+} := \max(0, 1 - \rho(\bu))$ ensures that there is no overcrowding.

Furthermore, \dt{for both M1 and M2 TAMs,} we consider a natural death rate $d_{M} > 0$ that is regulated by the available nutrients \dt{through the} death effect-function $\Psi_{dM}(\sigma)$ introduced in \eqref{Nutrient_Functions_Death_Macrophages}\dt{,} and so the death terms \dt{for each of the two phenotypes} are defined by
\begin{equation}
	Q_{M_{1}}(\bu) := d_{M} \Psi_{dM}(\sigma) M_{1}, \qquad
	Q_{M_{2}}(\bu) := d_{M} \Psi_{dM}(\sigma) M_{2},
	\label{Macrophage_Death}
\end{equation}
for the M1 and M2 TAMs populations, respectively.

\begin{figure*}
\centering
  \includegraphics[width=0.5\textwidth]{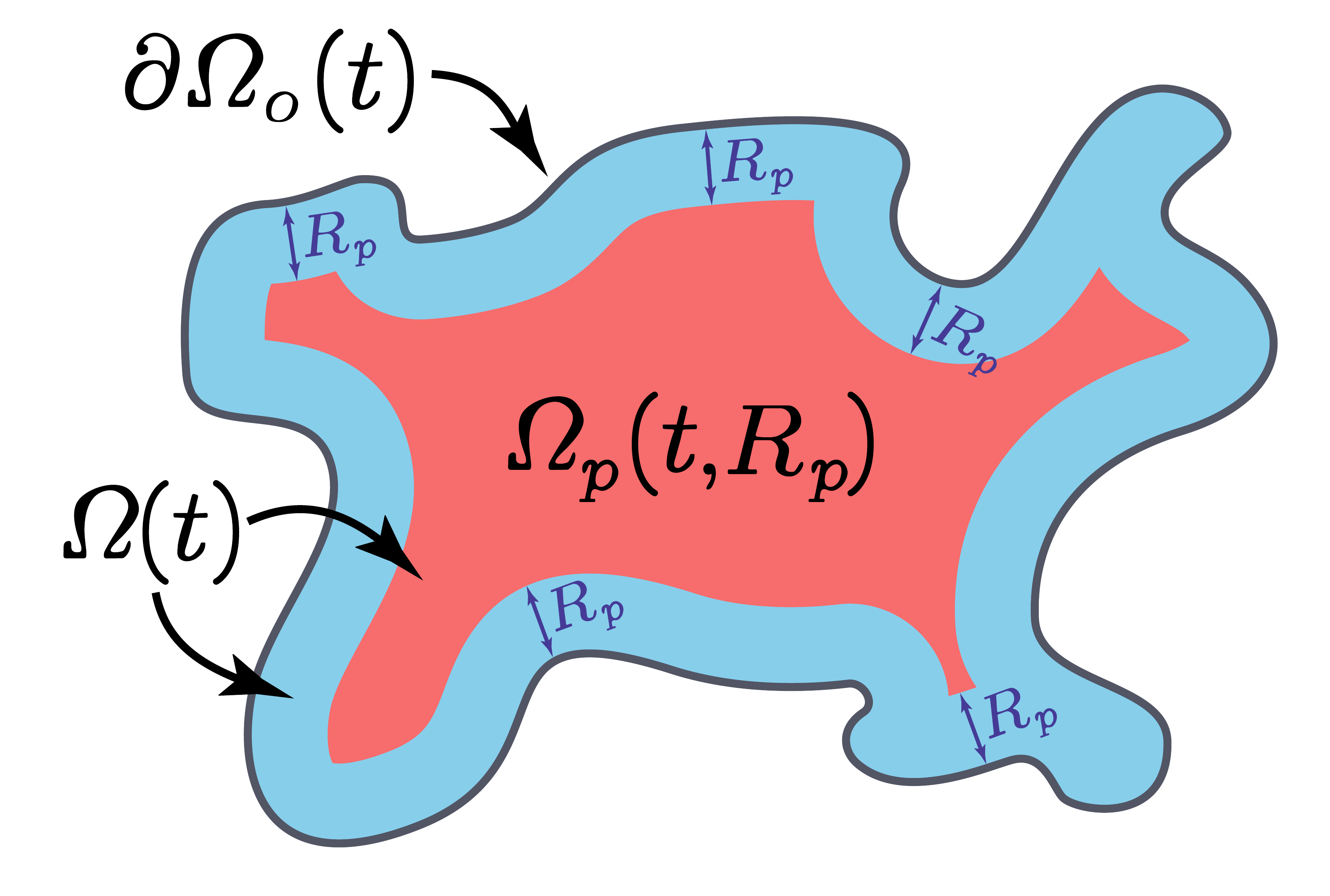}
\caption{Schematic of the re-polarisation domain $\Omega_{p}(t,R_{p})$ that is highlighted with red.}
\label{fig:Re_Polarisation_Domain}
\end{figure*}

Due to the versatility of the macrophages, their phenotype can be switched from one to another \cite{Becker_2013,Rocha_2014}. \dt{In the present work}, we focus on two factors that drive the polarisation of M1 TAMs into M2 TAMs, \dt{which are detailed as follows}. On the one hand, cytokines secreted by the cancer cells \dt{were} shown \cite{Pastuszak_Lewandoska_2018,Sica_2007} to trigger the polarisation process. \dt{O}n the other hand, the nutrient level was also shown \cite{Tan_2016} to \dt{a}ffect this process. \dt{As a consequence}, we describe the polarisation of M1 TAMs to M2 TAMs by
\begin{equation}
	T_{12}(\bu) := p_{12} \Psi_{M}(\sigma) c M_{1},
	\label{Macrophage_Polarisation}
\end{equation}
where $p_{12} > 0$ is \dt{a} constant proliferation rate\dt{,} and $\Psi_{M}$ is the polarisation effect-function defined in \eqref{Nutrient_Functions_Polarisation}. Further, in vitro, it has been demonstrated \cite{Davis_2013} that the M2-like macrophages can be re-polarised back into the M1 phenotype which may be a viable strategy against tumour development. To \dt{that} end, we explore \dt{here mathematically} the possibilities of \dt{the re-polarisation} strategy \dt{through a} re-polarisation term \dt{of the form}
\begin{equation}
	T_{21}(\bu) := 
	\begin{cases}
		\hfil 0 & \text{if } t < t_{p}, \\
		p_{21} M_{2} (\chi_{_{\dt{\Omega_{p}(t, R_{p})}}} \ast \psi_{\rho})(x) & \text{if } t \geq t_{p}.
	\end{cases}
	\label{Macrophage_Re_Polarisation}
\end{equation}
Here, $p_{21} > 0$ is the constant re-polarisation rate, $t_{p} > 0$ is the activation time and $\chi_{_{\dt{\Omega_{p}(t, R_{p})}}}$ is the characteristic function of the re-polarisation domain $\Omega_{p}(t, R_{p})$ that is defined in Appendix \ref{Appendix_Re_Polarisation_Domain} and illustrated in Fig.~\ref{fig:Re_Polarisation_Domain}. This re-polarisation term \eqref{Macrophage_Re_Polarisation} allows us to \dt{examine} whether \dt{or not} we \dt{would} need to \dt{account for} spati\dt{o}-temporal dependencies through the domain $\Omega_{p}(t, R_{p})$ \dt{and} activation time\dt{s} $t_{p} > 0$ \dt{in order to obtain} an effective re-polarisation strategy.
%------------ Figure for the re-polarrisation domain
%------------

The motility of both macrophages \dt{phenotipes is driven both by random and directed} movement. \dt{Based on recent biological evidence \cite{Hayenga_2015}, }increased stiffness of the substrate leads to an increase in macrophages’ speed, \dt{aspect explored in our modelling through a diffusion enhancement} that corresponds to with the level of ECM fibres. To \dt{that} end, we consider a stiffness-dependent macrophage diffusion coefficient $D^{M}(\bu)$ of the form
\begin{equation}
	D^M(\bu) := D_{M}(1 + D_{MF}F),
	\label{Macrophage_Diffusion_Function}
\end{equation}
where \dt{$D_{M} > 0$ is the baseline macrophage diffusion rate, and}  $D_{MF} > 0$ is the diffusion enhancement rate due to the presence of fibres. On the other hand, besides random movement, macrophages also exercise directed migration \dt{due to a both adhesive interactions with the surrounding cells and the ECM as well as an underlying cross-talk between themselves and the cancer cells}. \dt{A} similar \emph{``non-local flux term"} to the one introduced in \cite{Suveges_2020} \dt{is used here to} explore the \dt{complex} interactions of the cells distributed at $x \in \Omega(t_{0})$ with other cells within a sensing region $\Bila(0,R)$, \dt{and this accounts for: (1)} cell-cell TAMs self-adhesion \cite{Cui_2018}; \dt{(2)} nutrients level mediated movement\cite{Murdoch2004}; \dt{and (3)} the contribution of the cancer cells to the directional movement of the macrophages \cite{Chen_2011,Condeelis2006,Dutta_2018,Xuan_2014}. Specifically, \dt{the} contribution of the cancer cells to the directional movement of the macrophages account not only for the biological evidence that cancer cells can bind themselves to TAMs \cite{Chen_2011} but also for the fact that cancer cells can attract TAMs \cite{Condeelis2006,Dutta_2018,Xuan_2014} by secreting various chemokines. Although we neither model \dt{explicitly} the involved chemokine activities within this cross-talk nor the chemo\dt{-}attractant activities involved with the nutrients, here we appropriately account for both of them \dt{through the following} non-local flux term:
\begin{equation}
	\begin{split}
		\A_{M}(x,\! t,\! \textbf{u},\! \S_{MM}\!)\! :=\!\! \frac{1}{R}\!\! \!\! \int\limits_{\Bila(0,R)}\!\!\!  \!\!\! \!\K(y) n(y) \! \Big [ & S_{M\sigma} \big(1 \!-\! \sigma(x\!+\!y,t) \big) \!+\! \S_{Mc} c(x\!+\!y, t) \\
		& \!+\! \S_{MM}\! \big (\! M_{1}(x\!+\!y, t) \!+\! M_{2}(x\!+\!y, t)\!\big ) \! \Big  ] \! \big [ 1 \!-\! \rho(\textbf{u}) \big ]^{+},
	\end{split}
	\label{General_TAMs_Adhesion}
\end{equation}
where $R$ represents the radius of the sensing region $\Bila(0,R)$. Further, $\S_{Mc} > 0$ is the combined strength of the macrophage-cancer adhesion, and $\S_{M\sigma} > 0$ denotes strength of the macrophage-nutrient relationship, with both $\S_{Mc}$ and $\S_{M\sigma}$ \dt{being} assumed to \dt{maintain their individual values unchanged when considering the cases of} M1 and M2 TAMs populations. Furthermore, $\S_{MM}$ denotes the self-adhesion strength that differs for M1 and M2 TAMs \cite{Cui_2018}, \emph{i.e.},
\begin{itemize}
\item for M1 TAMs $\S_{MM} = \S_{M_{1}M} > 0$, and 
\item for M2 TAMs $\S_{MM} = \S_{M_{2}M} > 0$,
\end{itemize} 
with $\S_{M_{1}M} \neq \S_{M_{2}M}$. \dt{Finally}, to account in \eqref{General_TAMs_Adhesion} for the gradual weakening of these different adhesions as we move away from the centre $x$ within $\Bila(x,R)$, we use a radially symmetric kernel $\K(\cdot)$ that is given by
\begin{equation*}
	\K(y) = \psi \Big( \frac{y}{R} \Big), \qquad \forall y\in \Bila(0,R),
\end{equation*}
where $\psi(\cdot)$ is the standard mollifier defined in Appendix \ref{Standard_Mollifier}. Moreover, in \eqref{General_TAMs_Adhesion}, $[1 - \rho(\bu) \big]^{+}$ ensures that overcrowded tumour regions do not contribute to macrophage migration and $n(\cdot)$ is the unit radial vector given by
\begin{equation}
	n(y) :=
	\begin{cases}
		\dfrac{y}{\nor{y}_{2}} & \text{if } y \in \Bila(0, R) \setminus \{0\}, \\
		\hfil 0 & \text{if } y = 0.
	\end{cases}
	\label{Unit_Radial_Vactor}
\end{equation}
\dt{Thus, aggregating now all these cell movement aspects explored in} \eqref{M_influx} - \eqref{General_TAMs_Adhesion}\dt{,} the dynamics of the two distinct \dt{macrophages} phenotypes are mathematically formulated as
\begin{subequations}
	\begin{align}
		\frac{\partial M_{1}}{\partial t} \!= & \nabla \! \cdot \! [D^{M}\!(\bu)\! \nabla M_{1} \!-\! M_{1} \A_{M}(x,\!t,\!\bu,\S_{M_{1}M})] \!+\! P_{M_{1}}\!(\bu) \!-\! Q_{M_{1}}\!(\bu) \nonumber\\ 
		& -\! T_{12}\!(\bu) \!+\! T_{21}\!(\bu)\!+\! M_{I},\label{Both_TAMs_Equation_M1} \\[0.2cm]
		\frac{\partial M_{2}}{\partial t} \!= & \nabla \! \cdot \! [D^{M}\!(\bu)\! \nabla M_{2} \!-\! M_{2} \A_{M}(x,\!t,\!\bu,\S_{M_{2}M})] \!+\! P_{M_{2}}\!(\bu) \!-\! Q_{M_{2}}\!(\bu),\nonumber\\
		&+\! T_{12}\!(\bu) \!-\! T_{21}\!(\bu), \label{Both_TAMs_Equation_M1}
	\end{align}
	\label{Both_TAMs_Equation}
\end{subequations}
where $\S_{M_{1}M} > 0$ and $\S_{M_{2}M} > 0$ are the self-adhesion strengths of M1 and M2 TAMs, respectively.

%--------------------------------------------------------------
%		 				Cancer cells
%--------------------------------------------------------------
\subsubsection{Dynamics of the cancer cell population}
The third cell population that we consider \dt{at macro-scale} is the cancer cell population.  \dt{Crucially important for cancer development and invasion, the cancer cell proliferation is a complex process that is regulated by several processes involving nutrients and macrophages. From the modelling perspective, while we consider the proliferation process as being of logistic type \cite{Laird1964,Laird1965,Tjorve2017}, we explore the influence of nutrients and macrophages as follows. On the one hand,} similar to both TAMs populations, we consider the proliferation effect-function $\Psi_{p}(\sigma)$ defined in \eqref{Nutrient_Functions_Proliferation} to \dt{explore the influence of the available nutrients on} the rate of cancer cell proliferation. On the other hand, \dt{biological evidence shows that while} M2 TAMs promote \dt{cancer cell proliferation} \cite{Hu2015}, M1 TAMs inhibit\dt{s this} \cite{Liu_2018}. Thus, expand\dt{ing here} the proliferation law introduced in \cite{Suveges_2020} by accounting for the negative effect of M1 TAMs\dt{,} we obtain  leading to the following proliferation law:
\begin{equation}
	P_{c}(\bu) := \mu_{c}\Psi_{p}(\sigma) [1 - \mu_{cM_{1}}M_{1} + \mu_{cM_{2}}M_{2}] c [1-\rho(\bu)]^{+},
	\label{Cancer_Proliferation}
\end{equation}
where \dt{$\mu_{c} > 0$ is a baseline proliferation rate  that is being regulated by the available nutrients, being enhanced  by the M2 TAMs at a rate $\mu_{cM_{2}}>0$ and at the same time weakened by the presence of the M1 TAMs at a rate $\mu_{cM_{1}} > 0$. Again, here the term} $(1 - \rho(\bu))^{+}$ ensures that there is no overcrowding.

\dt{Besides proliferation, i}t is well known that cancer cells resist death \cite{Hanahan2000,Hanahan2011}. However, due to the peritumoral vasculature as well as the excessive degradation of the ECM, the efficiency of the nutrients delivery significantly reduces \dt{inside the tumour}, leading to necrosis \cite{Zhang2019}. In addition, numerous studies have shown \cite{Blankenstein1991,Lamagna2006,MacMicking1997,McBride1986,Nathan1987} that classically activated M1-like macrophages can produce significant amounts of pro-inflammatory cytokines\dt{,} and thereby have the ability to kill cancer cells. To \dt{that} end, \dt{we assume here} a baseline death rate $d_{c} > 0$ that is regulated not only by the cancer cell death effect-function $\Psi_{dc}(\sigma)$ \dt{introduced in} \eqref{Nutrient_Functions_Death_Cancer}, but also by the M1 TAMs \dt{at a rate} $d_{cM_{1}} > 0$. \dt{This results in the following mathematical representation of} the cancer cell death \dt{process, namely} 
\begin{equation}
	Q_{c}(\bu) := d_{c} [\Psi_{dc}(\sigma) + d_{cM_{1}} M_{1}] c.
	\label{Cancer_Death}
\end{equation}

\begin{figure*}
\centering
  \includegraphics[width=0.9\textwidth]{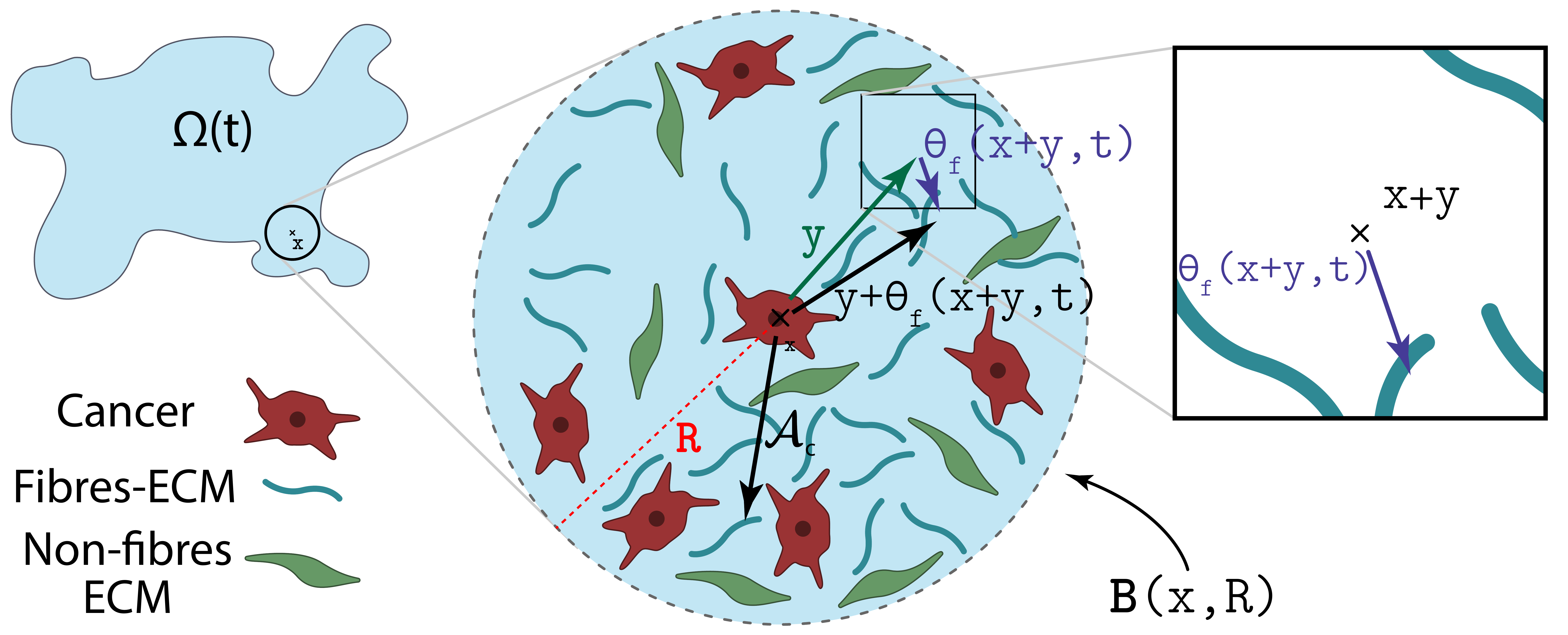}
\caption{Illustration of how the orientation of the fibres $\theta_{f}(\cdot, \cdot)$ biases the cell-fibre adhesion.}
\label{fig:Fibre_Adheison}
\end{figure*}

Similar to the macrophages, for the cancer cell population we also \dt{account for the diffusion enhancement that the} spatial \dt{distribution of}  ECM fibres \dt{enables} \cite{Allena2016,Cheung2009,Ebata2020,Isenberg2009,Lo2000,Raab2012,ReinhartKing2005,Saez2007,Trichet2012}. \dt{Furthermore, the random movement of the cell population is also affected by the presence of} both macrophage populations. While in general, the M2 TAMs were shown to promote cancer cell motility,\cite{Afik_2016} recent biological evidence \cite{Lee2020,Liu_2018} indicates that the M1 phenotype has a negative effect on the cancer cell motility. \dt{Therefore, the diffusion coefficient for the random movement of the cancer cells can be formulated mathematically as} 
\begin{equation}
	D^{c}(\bu) := D_{c}(1 + D_{cM_{2}}M_{2} + D_{cF}F - D_{cM_{1}}M_{1}).
	\label{Cancer_Diffusion_Function}
\end{equation}
\dt{where $D_{c} > 0$ is a baseline diffusion rate, $D_{cF} > 0$ is the ECM fibres enhancement coefficient, $D_{cM_{1}} > 0$  represents the weakening effect due to the presence of M1 TAM, and  $D_{cM_{2}} > 0$ accounts for the positive motility effect due to the presence of M2 TAM}.\\\\
In addition to random motility, we also take into account the directed movement of the cancer cells induced by various adhesion mediated processes \cite{Chen_2011,Condeelis2006,Huda2018,Petrie2009,Weiger2013,Wu3949}. Specifically, we extend \dt{here} the non-local flux term \dt{introduced in }\cite{Suveges_2020} (where the interaction of the cancer cells with other cancer cells, with M2 TAMs as well as with the two phase ECM within a sensing region $\Bila(0,R)$ were explored) by \dt{accounting for the} cancer cell M1 TAMs adhesion relationship. On the one hand, as we mentioned in Section~\ref{Section_M1M2}, cancer cells can bind themselves to TAMs and on the other hand, biological evidence \cite{Condeelis2006} suggest that cancer cells are also attracted by TAMs. \dt{As these two macrophages$-$cancer cells interactions (i.e., M1 TAM$-$cancer cells and M2 TAM$-$cancer) have both a non-local character}, \dt{for simplicity we assume that these are similarly distributed over the sensing region $\Bila(x,R)$ and have the same strength.} \dt{Therefore, the directed movement of the cancer cell population is governed by the following} non-local flux term:
\begin{equation}
	\begin{split}
		\A_{c}(x, t, \textbf{u}, \theta_{f}) := \frac{1}{R} \int\limits_{\Bila(0,R)}\!\!\! \!\!\! \K(y) \! \Big [ & n(y) \big ( \S_{cc} c(x\!+\!y, t) \!+\! \S_{cl} l(x\!+\!y, t)  \\
		& \!+\! \S_{cM} (M_{1}(x\!+\!y,t) \!+\! M_{2}(x\!+\!y,t)) \big) \\
		& \!+\! \widehat{n}(y,\! \theta_{f}(x+y,t)) \S_{cF} F(x\!+\!y, t) \Big ] \! \big [ 1 \!-\! \rho(\textbf{u}) \big ]^{+},
	\end{split}
	\label{Cancer_Adhasion}
\end{equation}
where $R$, $n(\cdot)$ and $\K(\cdot)$ are the same as in \eqref{General_TAMs_Adhesion}. Further, in \eqref{Cancer_Adhasion} $\widehat{n}(\cdot, \cdot)$ is the unit radial vector biased by the orientation of the fibres, \emph{i.e.,}
\begin{equation}
	\widehat{n}(y, \theta_{f}\dt{(x+y)}) :=
	\begin{cases}
		\dfrac{y + \theta_{f}(x+y,t)}{\nor{y + \theta_{f}(x+y,t)}_{2}} & \text{if } y \in \Bila(0, R) \setminus \{0\}, \\
		\hfil 0 & \text{if } y = 0,
	\end{cases}
	\label{Unit_Radial_Vector_Fibres}
\end{equation}
\dt{with} one of these biased vector $y + \theta_{f}(x+y,t)$ \dt{being} shown in Fig.~\ref{fig:Fibre_Adheison}.
%------------ Figure for the adhesion and adhesion
%------------
Moreover, in \eqref{Cancer_Adhasion} $\S_{cM} > 0$ represents the strength of the adhesion relationship between the cancer cells and M1 and M2 TAMs, $S_{cF} > 0$ is the strength of the cell-fibre ECM adhesion \cite{Wolf2009,Wolf_Friedl_2011} and $S_{cl} > 0$ corresponds to strength of adhesion between the cancer cells and the non-fibre ECM phase (that includes for instance amyloid fibrils\dt{,} which can support cell-adhesion processes \cite{Ghosh2017,Gras2009,Gras2008,Jacob2016}). \dt{Furthrmore, based on the} biological evidence \cite{Gu2014,Hofer2000} \dt{which} suggests that the emergence of strong and stable cancer self adhesion bonds are positively correlated with the high level of extracellular $Ca^{+2}$ ions (which is one of the constituents of the non-fibre ECM phase), \dt{proceeding as in \cite{Shuttleworth_2019,Shuttleworth2020b,Shuttleworth2020a,Suveges_2020} } the cancer cells self-adhesion coefficient $\S_{cc}$ is taken here as
\begin{equation*}
	\S_{cc}(x,t) :=\S_{min} + (\S_{max} - \S_{min}) \exp \bigg[ 1-\dfrac{1}{1-(1-l(x, t))^{2}} \bigg], 
\end{equation*}
where $\S_{max} > 0$ and $\S_{min} > 0$ correspond to maximum and minimum levels of $Ca^{+2}$ ions. Therefore, $\S_{cc}$ smoothly increases from a minimal to a maximum value in order to fully explore the varying strengths of cell-cell adhesion.

Thus, using \eqref{Cancer_Proliferation} to \eqref{Cancer_Adhasion} the spatio-temporal dynamics of the cancer population $c(x,t)$ is expressed as
\begin{equation}
	\begin{split}
		\frac{\partial c}{\partial t} = & \nabla \cdot [D^{c}(\bu) \nabla c - c\A_{c}(x,t,\bu,\theta_{f})] + P_{c}(\bu) - Q_{c}(\bu).
	\end{split}
	\label{Cancer_Equation}
\end{equation}

%--------------------------------------------------------------
%		 				Both ECM
%--------------------------------------------------------------
\subsubsection{Two-phase ECM \dt{macro-scale} dynamics}
Besides the cancer cells, both macrophage phenotypes contribute to the degradation of the ECM by secreting proteolytic enzymes \cite{Dollery2006,Goswami_2017,Madsen_2017,Newby2008,Rath_2019} (e.g., various classes of matrix metalloproteinases). To \dt{that} end, we extend the dynamics of the fibre, and non-fibre ECM components used in \cite{Suveges_2020} by incorporating the effects of the M1 phenotype. Thus, the dynamics of the non-fibre $l(x,t)$ as well as the fibre ECM $F(x,t)$ are formalized as
\begin{subequations}
	\begin{align}
		\dfrac{\partial l}{\partial t} = & - l(\beta_{lc}c + \beta_{lM_{1}}M_{1} + \beta_{lM_{2}}M_{2}) + (\gamma_{0} + \gamma_{M_{2}}M_{2})(1 - \rho(\bu))^{+}, \label{NonFibre_ECM_Equation}\\
		\dfrac{\partial F}{\partial t} = & - F(\beta_{Fc}c + \beta_{FM_{1}}M_{1} + \beta_{FM_{2}}M_{2}), \label{Fibre_ECM_Equation}
	\end{align}
	\label{Both_ECM_Equation}
\end{subequations}
where $\beta_{lc}$, $\beta_{lM_{1}}$, $\beta_{lM_{2}}$ are the positive degradation rates of the non-fibre ECM phase due to the cancer cells, M1 and M2 TAMs, respectively. Similarly, $\beta_{Fc}$, $\beta_{FM_{1}}$, $\beta_{FM_{2}}$ are all positive and describe the degradation rates of the fibre component of the ECM due to the cancer cells, M1 and M2 TAMs, respectively. Finally, in \eqref{Both_ECM_Equation} $\gamma_{0} > 0$ represents the constant rate of remodelling and $\gamma_{M_{2}} > 0$ is the remodelling enhancement rate induced by the M2 TAM population \cite{Afik_2016,Goswami_2017,Springer_2016}.

%--------------------------------------------------------------
%		 				Summary
%--------------------------------------------------------------
\subsubsection{The full macro-scale dynamics}
In summary, using \eqref{Nutrient_Equation}, \eqref{Both_TAMs_Equation}, \eqref{Cancer_Equation} and \eqref{Both_ECM_Equation} the macro-scale dynamics is given by the following coupled PDEs
\begin{subequations}
	\begin{align}
		\frac{\partial c}{\partial t} \!= & \nabla \!\! \cdot \! [D^{c}(\bu) \nabla c \!-\! c\A_{c}(x,\!t,\!\bu,\theta_{f})] \!+\! P_{c}(\bu) \!-\! Q_{c}(\bu), \label{Full_Macro_System_a}\\[0.2cm]
		\frac{\partial M_{1}}{\partial t} \!= & \nabla \!\! \cdot \! [D^{M}\!(\bu)\! \nabla \!M_{1} \!-\! M_{1} \A_{M}(x,\!t,\!\bu,\S_{M_{1}M}\!)] +\!\! P_{M_{1}}\!(\bu) \!-\! Q_{M_{1}}\!(\bu) \nonumber\\
		&-T_{12}\!(\bu) \!\!+\! \!T_{21}\!(\bu)\! \!+\! \!M_{\!I}, \label{Full_Macro_System_b}\\[0.2cm]
%		& \!+\! M_{I}, \\
		\frac{\partial M_{2}}{\partial t} \!= & \nabla \! \!\cdot \! [D^{M}\!(\bu)\! \nabla M_{2} \!-\! M_{2} \A_{M}(x,\!t,\!\bu,\S_{M_{2}M}\!)] +\! P_{M_{2}}\!(\bu) \!-\! Q_{M_{2}}\!(\bu)\nonumber\\
		&+ T_{12}\!(\bu) \!-\! T_{21}\!(\bu), \label{Full_Macro_System_c}\\[0.2cm]
		\dfrac{\partial l}{\partial t} \!= & - l(\beta_{lc}c \!+\! \beta_{lM_{1}}M_{1} \!+\! \beta_{lM_{2}}M_{2}) \!+\! (\gamma_{0} \!+\! \gamma_{M_{2}}M_{2})(1 \!-\! \rho(\bu)), \label{Full_Macro_System_d}\\[0.2cm]
		\dfrac{\partial F}{\partial t} \!= & - F(\beta_{Fc}c \!+\! \beta_{FM_{1}}M_{1} \!+\! \beta_{FM_{2}}M_{2}), \label{Full_Macro_System_e}\\[0.2cm]
		0 = & D_{\sigma}\Delta \sigma - d_{\sigma}(c \!+\! M_{1} \!+\! M_{2}), \label{Full_Macro_System_f}
	\end{align}
	\label{Full_Macro_System}
\end{subequations}
\dt{in the presence of appropriate initial conditions (such as those specified in \eqref{initalConditionMacroSystem17Dec2020})} with zero-flux boundary conditions for $c$, $M_{1}$, $M_{2}$, $l$ and $F$, as well as Dirichlet boundary condition \eqref{Nutrient_Equation} for the nutrients $\sigma$.

\subsection{Processes on the Micro-Scales and Links Between the Scales}\label{section:Micro_Scale_Model}
%--------------------------------------------------------------
%				Intro for the Micro-Scales
%--------------------------------------------------------------
\dt{As the} process of cancer invasion is truly a multi-scale phenomena \cite{Weinberg2006}, the macro-scale \dt{dynamics is} tightly linked together with several micro-scale processes. 
%\dt{and building on the multiscale modelling framework developed in Refs. \cite{Shuttleworth_2019,Shuttleworth2020b,Shuttleworth2020a,Suveges_2020,Dumitru_et_al_2013}, we explore this here}. 
Among, the micro-scale \dt{processes of important for cancer invasion, \dt{of main interest for us in this work are }the micro-scale rearrangement of} ECM fibre micro-constituents as well as the cell-scale proteolytic processes that take place at the leading edge of the tumour. In the following, we outline the details of these two micro-dynamics as well as the naturally occurring double feedback loop \dt{that} links them to the tumour macro-scale dynamics \eqref{Full_Macro_System}.

%--------------------------------------------------------------
%		 			Fibre micro-scale
%--------------------------------------------------------------
\subsubsection{Fibres on the micro-scale \dt{and their bottom-up and top-down links to macro-dynamics}}\label{Fibre_Micro_Scale}
Following \cite{Shuttleworth_2019}, the macroscopic \dt{oriented} ECM fibres are represented not only through their amount $F(x,t)$, but also via their spatial bias that characterize their ability to withstand incoming forces. Indeed, both of these characteristics of the ECM fibres are induced by the distribution of the micro-fibres $f(z,t)$ on a cell-scale micro-domain $\delta Y(x) := x + \delta Y$ \dt{of} appropriate micro-\dt{scale} size $\delta > 0$, \dt{ and are captured via} a vector field representation \cite{Shuttleworth_2019} $\theta_{f}(x,t)$ of the ECM fibres which is defined by
\begin{equation}
	\theta_{f}(x,t) := \dfrac{1}{\lambda (\delta Y(x))} \int\limits_{\delta Y(x)} f(z,t) dz \cdot \dfrac{\theta_{x,\delta Y(x)}(x,t)}{\nor{\theta_{x,\delta Y(x)}(x,t)}_{2}},
	\label{Orientation_Calculation}
\end{equation}
where $\lambda(\cdot)$ is the Lebesgue measure in $\R^{d}$ and $\theta_{f,\delta Y(x)}(\cdot,\cdot)$ is the revolving barycentral orientation with respect to the measure $f(z,t)\lambda(\cdot)$ \dt{given by\cite{Shuttleworth_2019}} 
%Thus, $\theta_{f,\delta Y(x)}(\cdot,\cdot)$ is uniquely defined by the mass distribution of the micro-fibres $f(z,t)$ on $\delta Y(x)$ and so it is given by the Bochner-mean-value\cite{yosida1980} of the \emph{barycentral vector-valued function}
\begin{equation*}
	\theta_{f,\delta Y(x)}(x,t) := \dfrac{\int\limits_{\delta Y(x)} f(z,t) (z-x) dz}{\int\limits_{\delta Y(x)} f(z,t) dz}.
\end{equation*}
Finally, the second characteristic, namely the ECM fibres amount, is given by the Euclidean norm of the vector field $\theta_{f}(x,t)$, \emph{i.e.,}
\begin{equation*}
	F(x,t) := \nor{\theta_{f}(x,t)}_{2},
\end{equation*}
and so it describes the mean-value of the micro-fibres distributed on $\delta Y(x)$. Therefore, the micro-scale naturally emerges a link with the macro-scale since the representation of the ECM fibres on the macro-scale is induced by the micro-scale fibre distribution, and so we refer to this as the \emph{fibres bottom-up} link.

\begin{figure*}
\centering
  \includegraphics[width=0.65\textwidth]{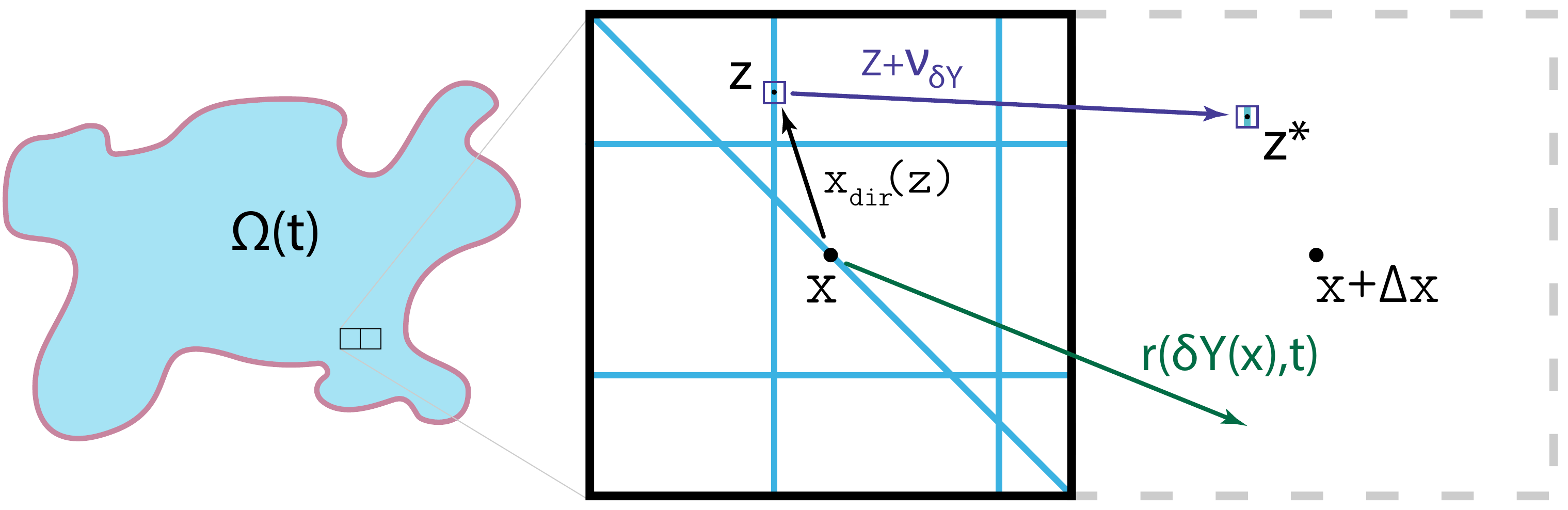}
\caption{Typical examples of the relevant vectors $x_{dir}(z) \dt{:}= z - x$, $r(\delta Y(x),t)$ and $\nu_{\delta Y(x)}(z,t)$, allowing the redistribution of each micro-point $(z,t)$.}
\label{fig:Fibre_Rearrangement_Vectors}
\end{figure*}

\begin{figure*}
\centering
  \includegraphics[width=0.9\textwidth]{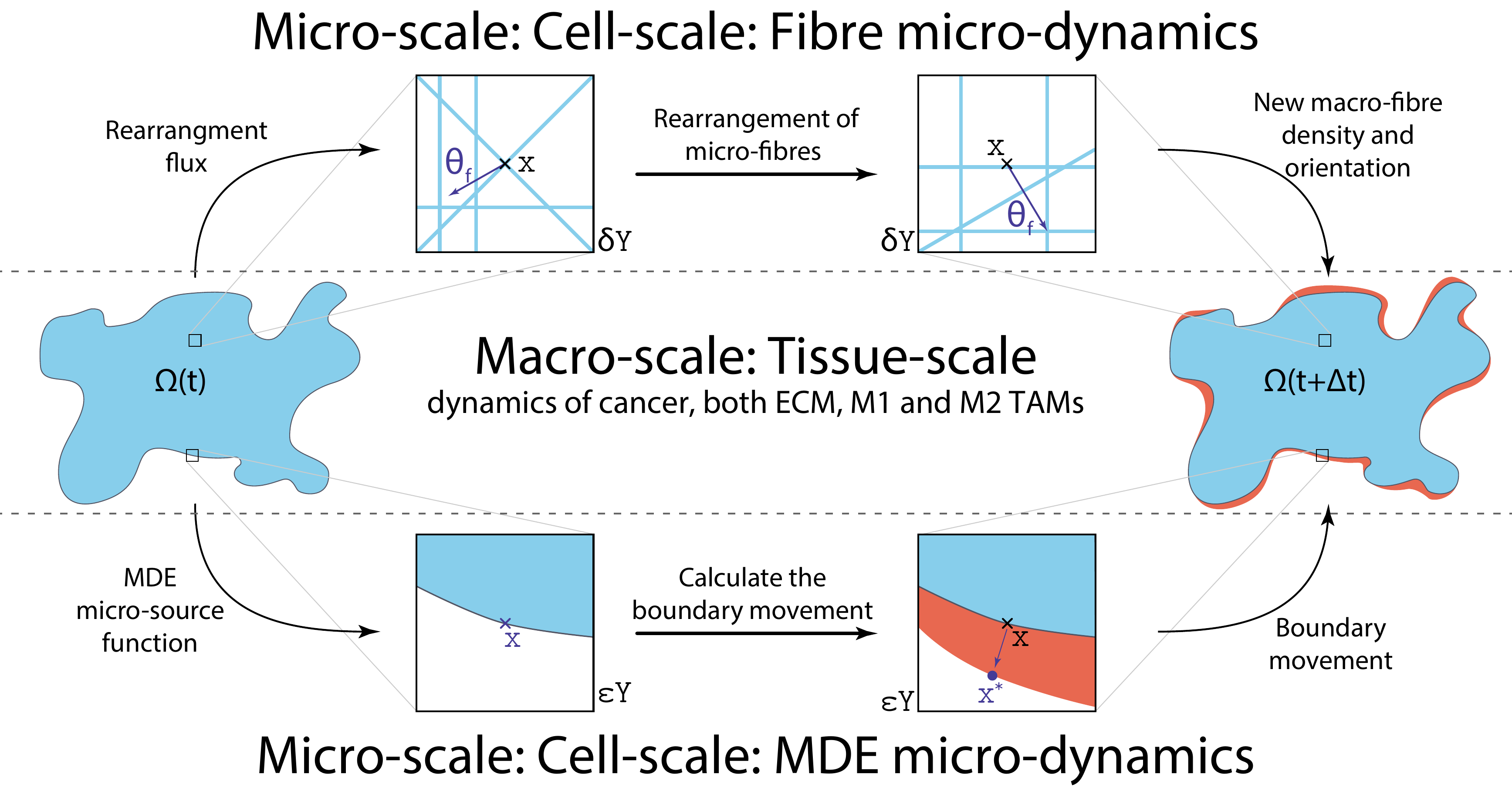}
\caption{Schematics of the four links between the macro and both micro-scales as well as how they are linked together.}
\label{fig:Scale_Links}
\end{figure*}

On the other hand, the macro-scale spatial fluxes, generated by the tumour macro-dynamics \eqref{Full_Macro_System}, \dt{trigger a} rearrangement of the ECM fibres micro constituents on each micro-domain $\delta Y(x)$. \dt{Indeed, the collective migration of  the cancer cells, M1 TAMS, and M2 TAMs lead naturally to the emergence of the associated spatial fluxes $\F_{c}$, $\F_{M_{1}}$, and $\F_{M_{2}}$ given by} 
%Hence, the micro-fibre rearrangement process is initiated \dt{collectively by the spatial fluxes associated} by the cancer cell $\F_{c}$, M1 TAMs $\F_{M_{1}}$ and M2 TAMs $\F_{M_{2}}$ fluxes that are given by
\begin{equation*}
	\begin{split}
		\F_{c}(x,t) := & D^{c}(M_{1}, M_{2}, F) \nabla c - c\A_{c}(x,t,\bu,\theta_{f}), \\
		\F_{M_{1}}(x,t) := & D^{M}(F) \nabla M_{1} - M_{1} \A_{M}(x,t,\bu,\S_{M_{1}M}), \\
		\F_{M_{2}}(x,t) := & D^{M}(F) \nabla M_{2} - M_{2} \A_{M}(x,t,\bu,\S_{M_{2}M}),
	\end{split}
\end{equation*}
respectively. \dt{The combined action of these fluxes upon the ECM fibres distributed at $(x,t)\in \Omega_{t}\times [0,T]$ is felt uniformly by its constituent micro-fibres $f(z,t)$ distributed on the associated micro-domain} $\delta Y(x)$, consequently \dt{inducing} a \emph{micro-fibres rearrangement vector} \dt{similar to the ones proposed in \cite{Shuttleworth_2019,Suveges_2020} of the form}
\begin{equation}
	\begin{split}
		r(\delta Y(x),t) := & \omega_{c}(x,t)\F_{c}(x,t) + \omega_{M_{1}}(x,t)\F_{M_{1}}(x,t) \\
		& + \omega_{M_{2}}(x,t)\F_{M_{2}}(x,t) + \omega_{F}(x,t)\theta_{f}(x,t),
	\end{split}
	\label{Rearrangement_Vector}
\end{equation}
which \dt{triggers a} spatial redistribution of the micro-fibres in $\delta Y(x)$. In \eqref{Rearrangement_Vector}, the non-linear weights $\omega_{c}$, $\omega_{M_{1}}$, $\omega_{M_{2}}$ and $\omega_{F}$ are the associated mass factions of cancer cells, M1 TAMs, M2 TAMs\dt{,} and ECM fibres \dt{at $(x,t)$,} and so these are 
\begin{equation*}
	\begin{split}
		\omega_{c}(x,t) := & \dfrac{c(x,t)}{c(x,t) + M_{1}(x,t) + M_{2}(x,t) + F(x,t)}, \\
		\omega_{M_{1}}(x,t) := & \dfrac{M_{1}(x,t)}{c(x,t) + M_{1}(x,t) + M_{2}(x,t) + F(x,t)}, \\
		\omega_{M_{2}}(x,t) := & \dfrac{M_{2}(x,t)}{c(x,t) + M_{1}(x,t) + M_{2}(x,t) + F(x,t)}, \\
		\omega_{F}(x,t) := & \dfrac{F(x,t)}{c(x,t) + M_{1}(x,t) + M_{2}(x,t) + F(x,t)},
	\end{split}
\end{equation*}
respectively. 
Ultimately, the rearrangement vector \eqref{Rearrangement_Vector} induces a \emph{relocation vector} $\nu_{\delta Y(x)}(z,t)$ \dt{(detailed in Appendix \ref{Fibre_Reallocation_Process}),} and as a result we can appropriately calculate the new positions of any micro-scale node $z \in \delta Y(x)$ that are given by
\begin{equation*}
	z^{*} := z + \nu_{\delta Y(x)}(z,t).
\end{equation*}
In Fig.~\ref{fig:Fibre_Rearrangement_Vectors}, we illustrate the micro-fibre rearrangement process by considering a typical example of these vectors while for further details we refer the reader to Appendix \ref{Fibre_Reallocation_Process} and \cite{Shuttleworth_2019,Suveges_2020}. \dt{Finally, this rearrangement of the ECM fibres at micro-scale triggered by the emergent macro-scale spatial fluxes ($\F_{c}$, $\F_{M_{1}}$, and $\F_{M_{2}}$) establishes the \emph{fibre top-bottom} link.}
%------------ Figure for the vectors that we use to calculate the fibre rearrangement
%------------
In Fig.~\ref{fig:Scale_Links}, we illustrate the fibres bottom-up and top-down links, connecting the macro-scale and the ECM fibre micro-scale.
%------------ Figure for the links between the scales
%------------

%--------------------------------------------------------------
%		 			 MDE micro-scale
%--------------------------------------------------------------
\subsubsection{MDE boundary \dt{micro-dynamics and its connections to macro-dynamics}}\label{MDE_Micro_Scale}
The second micro-process that we consider is the proteolytic molecular processes which are driven by the secretion of \emph{matrix-degrading enzymes} (MDEs) (such as matrix-metalloproteinases) and takes place along the leading edge of the tumour. \dt{Indeed, b}esides cancer cells, both M1 and M2 TAMs \dt{within the proliferating outer rim of the tumour} secrete MDEs \cite{Dollery2006,Goswami_2017,Madsen_2017,Newby2008,Rath_2019}, \dt{and this way the tissue-scale dynamics induces a source for cell-scale proteolytic activity. Once secreted, the MDE molecular population exercises a spatial transport across the tumour interface within a cell scale neighbourhood of the tumour boundary, causing} degradation of the peritumoural ECM, \dt{and ultimately} leading to changes in the tumour morphology \cite{Weinberg2006}. \dt{Thus, following here the approach introduced in Ref~\cite{Dumitru_et_al_2013}, we explore the emerging spatio-temporal molecular MDEs micro-dynamics on an appropriate cell-scale neighbourhood of the tumour interface $\partial \Omega(t)$ enabled by the union of a covering bundle of cubic micro-domains $\{ \epsilon Y \}_{\epsilon Y \in \P(t)}$, with each $\epsilon Y$ being of micro-scale size $\epsilon>0$. This allows us to decompose the MDE micro-dynamics occurring on $\bigcup\limits_{\epsilon Y \in \P(t)} \epsilon Y$ into a union of MDE micro-processes occurring on each $\epsilon Y$} \cite{Dumitru_et_al_2013}.\\

\dt{Therefore, considering the tumour evolution over a time perspective $[t_{0},t_{0}+\Delta t]$, for an arbitrary instance $t_{0}\in[0,T]$, and of appropriate micro-scale range $\Delta t>0$, on any of the micro-domains $\epsilon Y\in \P(t_{0})$ we denote by $m(y, \tau)$ the spatio-temporal distribution of MDEs at micro-scale point $(y,\tau)\in \epsilon Y\times [0,\Delta t]$.} 
\dt{In this context,} at any spatio-temporal  $(y,\tau) \in (\epsilon Y \cap \Omega(t_{0})) \times [0, \Delta t]$, a source of MDEs arises as a collective contribution of the cancer cell and both macrophage populations \dt{from the outer proliferating rim of the tumour that are situated} within a distance $\gamma_{h}\dt{>0}$ from \dt{$y\in\epsilon Y$}. Hence, \dt{denoting this micro-scale MDE source by $h(y, \tau)$, this} can be formalised mathematically via the non-local expression
\begin{equation}
	h(y,\tau) = 
	\begin{cases}
		\dfrac{\int\limits_{\Bila(y, \gamma_{h}) \cap \Omega (t_{0})} h_{\Sigma}\dt{(x,\tau)} \; dx}{\lambda (\Bila(y, \gamma_{h}) \cap \Omega (t_{0}))} & y \in \epsilon Y \cap \Omega(t_{0}), \\[10pt]
		\hfil 0 & y \notin \epsilon Y \setminus (\Omega(t_{0}) + \{ z \in Y \; | \; \| z \|_{2} < \rho \}),
	\end{cases}
	\label{MDE_Source}
\end{equation}
where $\Bila(y, \gamma_{h}) := \{ z \in Y | \nor{y-z}_{_{\infty}} \leq \gamma_{h} \}$ denotes the $\nor{\cdot}_{_{\infty}}$ ball with appropriately chosen radius $\gamma_{h} > 0$ and $0 < \rho < \gamma_{h}$ is a small mollification range which smooths \dt{out} the source function $h(\cdot, \cdot)$. Further, in \eqref{MDE_Source} $h_{\Sigma}$ is given by
\begin{equation*}
	h_{\Sigma}\dt{(x,\tau)}:=\alpha_{c} c(x, t_{0}+\tau) + \alpha_{M_{1}} M_{1}(x, t_{0}+\tau) + \alpha_{M_{2}} M_{2}(x, t_{0}+\tau),
\end{equation*}
where $\alpha_{c} > 0$, $\alpha_{M_{1}} > 0$ and $\alpha_{M_{2}} > 0$ are constant \dt{secretion rates of the MDEs by the} cancer cells, M1 and M2 TAMs respectively. As the MDE micro-source is naturally induced by the macro-scale, this establishes a \emph{MDE top-down} link between the tumour macro-dynamics and MDE-micro-dynamics occurring at the tumour interface. Finally, \dt{under the presence of} the MDE source $h(\cdot, \cdot)$, the MDE micro-dynamics \dt{is given} by
\begin{equation}
\begin{array}{l}
	\dfrac{\partial m}{\partial \tau} = D_{m} \Delta m + h(y, \tau), \quad\\[0.3cm]
	m(y,0) = 0, \quad \quad\\[0.3cm]
	\dfrac{\partial m}{\partial n} = 0,
\end{array}
	\label{MDE_Equation}
\end{equation}
where $D_{m} > 0$ is the constant diffusion coefficient and $n$ is the outward normal vector. \dt{Finally, as the patterns of significant degradation of the peritumoural ECM correspond the regions where significant levels MDEs are transported during the micro-dynamics, }the micro-dynamics \eqref{MDE_Equation} \dt{enables us to capture to changes in tumour morphology by} determining the direction \dt{of movement} and magnitude of the displacement of invading cancer in the surrounding tissue \cite{Dumitru_et_al_2013}. Thus, the MDE micro-process \dt{induces changes in the shape of the} tumour boundary $\partial\Omega(t_{0})$, affecting \dt{directly} the macro-dynamics, \dt{which is continued on  a modified tumour macro-domain $\Omega(t_{0}+\Delta t)$, see Appendix \ref{Appendix_MDE_Micro_Process} for details. Hence, a} \emph{MDE bottom-up} link \dt{between the  proteolytic boundary micro-dynamics and macro-scale cancer dynamics is this way established.} In Fig.~\ref{fig:Scale_Links}, we illustrate both MDE top-down and bottom-up links that connect the macro-scale and the MDE micro-scale.
%--------------------------------------------------------------
%	  						Numerics
%--------------------------------------------------------------
\section{Numerical Approach}\label{section:Numerical_Approach}
%--------------------------------------------------------------
%	  				Intro for the numerics
%--------------------------------------------------------------
In this section, we \dt{describe} the numerical \dt{approach developed} to solve the tumour macro-scale dynamics \eqref{Full_Macro_System}. First, to solve the quasi-steady nutrients $\sigma$ equation \eqref{Full_Macro_System_f}, we use the usual successive over-relaxation method with relaxation parameter $\omega = 0.5$ and tolerance of $10^{-5}$. For the rest of the dynamics \eqref{Full_Macro_System_a}-\eqref{Full_Macro_System_e}, we use the method of lines approach to discretise the system \eqref{Full_Macro_System} first in space and then for the resulting ODEs we use a non-local predictor-corrector scheme \cite{Shuttleworth_2019}. In this context, we carry out accurate approximation of the two distinct spatial operators, namely the diffusion and adhesion operators by using fast convolution-driven approaches. Specifically, while for the diffusion we construct a convolution-based second-order central difference scheme \cite{Suveges_2020}, for the adhesion operators, we formulate a convolution-driven \emph{fifth-order weighted essentially non-oscillatory} (WENO5) finite difference scheme. 

For completeness, we also detail the numerical approach of approximating the adhesion integrals \eqref{General_TAMs_Adhesion} and \eqref{Cancer_Adhasion} in \ref{Appendix_Adhesion} as well as the details of MDE micro-scale calculations is presented in \ref{Appendix_MDE_Calculation}, while for further details we refer the reader to \cite{Shuttleworth_2019,Suveges_2020,Dumitru_et_al_2013}.

To facilitate the description of the numerical approaches, let us first introduce some basic notations. For simplicity, we discretise the maximal tissue cube $Y \in \R^2$ uniformly in each direction, \emph{i.e.,} let $L$ be the length of $Y$ and $h_{L} = \Delta x = \Delta y$ to be the uniform spatial step-size. Therefore, the discretisation of $Y$ can be represented through the macro-spatial locations $\{ (x_{i}, y_{j}) \}_{i,j=1..N}$, with $N = L / h_{L} + 1$. Then the discretised global five-dimensional tumour vector at any time-step $n$ is denoted by $\bu^{n} = [c^{n}, M_{1}^{n}, M_{2}^{n}, l^{n}, F^{n}]^{\intercal}$, \dt{while the discretised} nutrient field is denoted by $\sigma^{n}$. Further, we denote the diffusion coefficient functions and adhesion integrals for cancer cell and TAMs by $(D^{c})^{n}$, $(D^{M})^{n}$, $(\A_{c})^{n}$ and $(\A_{M})^{n}$, respectively. Finally, \dt{ with these notations, we are able to} rewrite the macro-scale dynamics \eqref{Full_Macro_System} in \dt{the following convenient} form
\begin{equation}
	\begin{split}
		\frac{\partial \bu}{\partial t} = & \nabla \cdot [\D(\bu) \nabla \bu] - \nabla \cdot \F(\bu) + \SS(\bu), \\
		0 = & D_{\sigma}\Delta \sigma - d_{\sigma}(c + M_{1} + M_{2}),
	\end{split}
	\label{Full_Macro_System_Compact}
\end{equation}
where
\begin{equation*}
\begin{split}
	\D(\bu) := &
	\begin{bmatrix} 
    	D^{c}(\bu) \\
    	D^{M}(\bu) \\
    	D^{M}(\bu) \\
    	0 \\
    	0 \\
    \end{bmatrix},
    \qquad
   	\F(\bu) :=
    \begin{bmatrix} 
    	c \A_{c}(x,t,\bu,\theta_{f}) \\
    	M_{1} \A_{M}(x,t,\bu,\S_{M_{1}M}) \\
    	M_{2} \A_{M}(x,t,\bu,\S_{M_{2}M}) \\
    	0 \\
    	0 \\
    \end{bmatrix},
	\\
    \SS(\bu) := &
    \begin{bmatrix} 
    	P_{c}(\bu) - Q_{c}(\bu) \\
    	P_{M_{1}}(\bu) - Q_{M_{1}}(\bu) - T_{12}(\bu) + T_{21}(\bu) + M_{I} \\
    	P_{M_{2}}(\bu) - Q_{M_{2}}(\bu) + T_{12}(\bu) - T_{21}(\bu) \\
    	- l(\beta_{lc}c + \beta_{lM_{1}}M_{1} + \beta_{lM_{2}}M_{2}) + (\gamma_{0} + \gamma_{M_{2}}M_{2})(1 - \rho(\bu)) \\
    	- F(\beta_{Fc}c + \beta_{FM_{1}}M_{1} + \beta_{FM_{2}}M_{2}) \\
    \end{bmatrix}.
\end{split}
\end{equation*}

%--------------------------------------------------------------
%					Diffusion operators
%--------------------------------------------------------------
\subsection{Diffusion Operators: $\nabla \cdot [\D(\bu) \nabla \bu]$}\label{Diffusion_Operators_Numeics}
Starting with the discretisation of the diffusion operators $\nabla \cdot [\D(\bu) \nabla \bu]$ in \eqref{Full_Macro_System_Compact}, during the computations, we first detect whether a spatial node $(i,j)$ is inside or outside the expanding tumour domain $\Omega(t_{0})$ via an indicator function $\I (\cdot, \cdot) : X \times X \rightarrow \{0,1\}$, with $X = \{ 1,...,M \}$ that is defined by 
\begin{equation}
	\I(i, j) := 
	\begin{cases}
	1 & \text{if } (x_{i}, x_{j}) \in \Omega (t_{0}), \\
	0 & \text{otherwise}.
	\end{cases}
	\label{Indicator_All}
\end{equation}
Similarly, we detect any boundary nodes using a boundary indicator function defined by
\begin{equation}
	\I_{B}(i, j) := 
	\begin{cases}
	1 & \text{if } (x_{i}, x_{j}) \in \partial \Omega (t_{0}), \\
	0 & \text{otherwise},
	\end{cases}
	\label{Indicator_Boundary}
\end{equation}
where $\partial \Omega (t_{0})$ is boundary of $\Omega (t_{0})$, and for convenience we also define the interior indicator function by
\begin{equation}
	\I_{In}(i, j) := 
	\begin{cases}
	1 & \text{if } (x_{i}, x_{j}) \in \Omega (t_{0}) \setminus \partial \Omega (t_{0}). \\
	0 & \text{otherwise}.
	\end{cases}
	\label{Indicator_Interior}
\end{equation}
These two functions, given in \eqref{Indicator_Boundary} and \eqref{Indicator_Interior}, enable us to split the domain into two parts, namely to boundary and strictly inside parts. The motivation behind this is to use two different computational technique on these parts which eventually reduces the computational cost. Hence, while for any interior node we can use the same discrete universal numerical operator (convolution), for a boundary point we need to apply unique operators due to the zero-flux boundary condition and the continuously changing tumour domain (that may result in an irregular tumour shape). Therefore, to achieve the reduction in computational cost, we accordingly split the spatial operators $\nabla \cdot [\D(\bu) \nabla \bu]$ into two components as well, and so at any spatial node $(x_{i}, y_{j}) \in \Omega(t_{0}) \subset Y$, the diffusion operators can be represented as
\begin{equation}
	\nabla \cdot [\D(\bu) \nabla \bu](i,j) := 
	\begin{cases}
	\nabla \cdot [\D(\bu) \nabla \bu]^{In}(i,j) & \text{if } \I_{In}(i, j) = 1, \\
	\nabla \cdot [\D(\bu) \nabla \bu]^{B}(i,j) & \text{if } \I_{B}(i, j) = 1, \\
	\hfil 0 & \text{otherwise}.
	\end{cases}
	\label{Splitted_Diffusion}
\end{equation}
In this context, the usual two dimensional second-order central difference scheme for a non-constant diffusion operator is given by
\begin{equation}
	\begin{split}
		\big( \nabla& \cdot [\D(\bu)\nabla \bu] \big)^{n}_{i,j} =\\
		& \dfrac{1}{\Delta x^2} \bigg(  \dfrac{\D(\bu^{n}_{i+1,j}) + \D(\bu^{n}_{i,j})}{2} \big( \bu^{n}_{i+1,j} - \bu^{n}_{i,j} \big) 
		 - \dfrac{\D(\bu^{n}_{i,j}) + \D(\bu^{n}_{i-1,j})}{2} \big( \bu^{n}_{i,j} - \bu^{n}_{i-1,j} \big) \\
		&\quad\,\,\, + \dfrac{\D(\bu^{n}_{i,j+1}) + \D(\bu^{n}_{i,j})}{2} \big( \bu^{n}_{i,j+1} - \bu^{n}_{i,j} \big) 
		 - \dfrac{\D(\bu^{n}_{i,j}) + \D(\bu^{n}_{i,j-1})}{2} \big( \bu^{n}_{i,j} - \bu^{n}_{i,j-1} \big) \!\!\bigg),
	\end{split}
	\label{Diffusion_Numerics_Normal}
\end{equation}
where $\bu^{n}_{i,j} = [c^{n}_{i,j}, M^{n}_{1_{i,j}}, M^{n}_{2_{i,j}}, l^{n}_{i,j}, F^{n}_{i,j}]^{\intercal}$. Further, we observe that \eqref{Diffusion_Numerics_Normal} can be equivalently expressed by sum of discrete convolutions and so for all interior node $(x_{i},y_{j})$ the scheme is given by
\begin{equation}
	\begin{split}
		\big( \nabla \cdot [\D(\bu) \nabla \bu]^{In} \big)^{n} = \dfrac{1}{\Delta x^2} \sum_{k=1}^{2} \Bigg( & \Big( \widetilde{\K}^{2k-1}_{A} \ast \D(\bu^{n}) \Big) \circ \Big( \widetilde{\K}^{k}_{F} \ast \bu^{n} \Big) \\
		& - \Big( \widetilde{\K}^{2k}_{A} \ast \D(\bu^{n}) \Big) \circ \Big( \widetilde{\K}^{k}_{B} \ast \bu^{n} \Big) \Bigg),
	\end{split}
	\label{Diffusion_Numerics_Convolution}
\end{equation}
where $\ast$ is the discrete convolution and $\circ$ is the Hadamard product, that is defined in \ref{Appendix_Hadamard} (and for further details we refer the reader to Ref~\cite{Golub_2013}). Also, in \eqref{Diffusion_Numerics_Convolution} each $\widetilde{\K}^{k}_{A}$, with $k=1\dots4$ describe the average between the point $(i,j)$ and one of its immediate neighbour and so they are defined by
\begin{equation}
	\widetilde{\K}^{1}_{A} = [0, 0.5, 0.5], \quad
	\widetilde{\K}^{2}_{A} = [0.5, 0.5, 0], \quad
	\widetilde{\K}^{3}_{A} = \big( \widetilde{\K}^{1}_{A} \big)^{\intercal}, \quad
	\widetilde{\K}^{4}_{A} = \big( \widetilde{\K}^{2}_{A} \big)^{\intercal}.
	\label{Kernels_Agerage}
\end{equation}
Moreover, in \eqref{Diffusion_Numerics_Convolution} $\widetilde{\K}^{k}_{F}$ and $\widetilde{\K}^{k}_{B}$ with $k = 1,2$ are induced by the forward and backward differences, respectively. Hence, they are defined in both direction $i$ (if $k = 1$) as well as in direction $j$ (if $k = 2$) by
\begin{equation}
	\widetilde{\K}^{1}_{F} = [0, -1, 1], \quad
	\widetilde{\K}^{1}_{B} = [-1, 1, 0], \quad
	\widetilde{\K}^{2}_{F} = \big( \widetilde{\K}^{1}_{F} \big)^{\intercal}, \quad
	\widetilde{\K}^{2}_{B} = \big( \widetilde{\K}^{1}_{B} \big)^{\intercal}.
	\label{Kernels_Forward_Backward}
\end{equation}
%\begin{remark}
%It is well known that convolution flips one of the input function by definition. Hence, the induced vectors defined in \eqref{Kernels_Agerage} and \eqref{Kernels_Forward_Backward} are already flipped accordingly to make equations \eqref{Diffusion_Numerics_Normal} and \eqref{Diffusion_Numerics_Convolution} equivalent.
%\end{remark}

However, for boundary nodes, we cannot use the form \eqref{Diffusion_Numerics_Convolution} due to the imposed zero-flux boundary condition and the continuously changing tumour domain. This is because, the calculation of the diffusion operators at the boundary nodes involves the approximation of the values at any node that does not belong to the tumour domain, \emph{i.e.,} for any node $(x_{i},y_{j}) \notin \Omega(t_{0})$. In some cases, due to the irregular domain, such values may not be uniquely defined because multiple nodes can require the value of the same ghost/outside node, however with different values. Consequently, for any boundary node, instead of the convolutional form \eqref{Diffusion_Numerics_Convolution}, we rather use the usual form \eqref{Diffusion_Numerics_Normal} and symmetrically reflect the values of the interior nodes to the ghost/outside nodes in a node by node fashion.

%--------------------------------------------------------------
%				  Adhesion operators
%--------------------------------------------------------------
\subsection{Adhesion Operators: $\nabla \cdot \F(\bu)$}
The other spatial operators that contribute to the motility are the adhesion operators $\nabla \cdot \F(\bu)$ \dt{in} \eqref{Full_Macro_System_Compact}. The procedure to approximate these differential operators is based on the standard WENO5 scheme proposed in \cite{Jiang1996,Liu1994}. However, it was shown \cite{Zhang2006} that these standard WENO5 scheme suffer from slight post-shock oscillations. \dt{As a consequence, we adopt here the modified WENO5 scheme proposed in \cite{Zhang2006}, which uses modified smoothness indicators and is referred to as the ZSWENO scheme}.

Similarly to the diffusion operator case, here we split the domain into two parts as well. However, ZSWENO schemes induce larger stencils compared to the second-order central difference scheme and so, we need to split the domain $\Omega(t_{0})$ into two parts differently. We refer these two parts as the \emph{inside} and \emph{layer} parts. First, the former one is detected by using an inside indicator function $\I_{I}(i, j)$ that we define by
\begin{equation}
	\I_{I}(i, j) := 
	\begin{cases}
	1 & \text{if } (\I \ast \K_{I})_{i,j} = 1, \\
	0 & \text{otherwise},
	\end{cases}
	\label{Indicator_Inside}
\end{equation}
where $\I$ is defined in \eqref{Indicator_All} and $\K_{I}$ is given by
\begin{equation*}
	\K_{in} = \frac{1}{13}
	\begin{bmatrix}
		0 & 0 & 0 & 1 & 0 & 0 & 0 \\
    	0 & 0 & 0 & 1 & 0 & 0 & 0 \\
    	0 & 0 & 0 & 1 & 0 & 0 & 0 \\
    	1 & 1 & 1 & 1 & 1 & 1 & 1 \\
    	0 & 0 & 0 & 1 & 0 & 0 & 0 \\
    	0 & 0 & 0 & 1 & 0 & 0 & 0 \\
    	0 & 0 & 0 & 1 & 0 & 0 & 0 \\
	\end{bmatrix},
\end{equation*}
which is induced by the fifth-order ZSWENO stencils. Then the layer part of the domain is detected by a layer indicator function formulated as
\begin{equation}
	\I_{L}(i, j) := 
	\begin{cases}
	1 & \text{if } \I(i,j) - \I_{I}(i, j) = 1, \\
	0 & \text{otherwise}.
	\end{cases}
	\label{Indicator_Layer}
\end{equation}
Similarly to the diffusion, using these indicator functions \eqref{Indicator_Inside} and \eqref{Indicator_Layer}, we split the adhesion operator into two parts and so at any spatial node $(x_{i}, y_{j}) \in \Omega(t_{0})$ this operator is represented as
\begin{equation}
	\nabla \cdot \F(\bu)(i,j) := 
	\begin{cases}
	\nabla \cdot \F(\bu)^{I}(i,j) & \text{if } \I_{I}(i, j) = 1, \\
	\nabla \cdot \F(\bu)^{L}(i,j) & \text{if } \I_{L}(i, j) = 1, \\
	\hfil 0 & \text{otherwise}.
	\end{cases}
	\label{Splitted_Adhesion}
\end{equation}
\dt{Therefore, for} the inside operator $\nabla \cdot \F(\bu)^{I}$, we use the usual conservative form
\begin{equation}
	\nabla \cdot \F(\bu)^{I} = \dfrac{1}{\Delta x} \Big( \hat{F}_{i+\frac{1}{2},j} - \hat{F}_{i-\frac{1}{2},j} + \hat{G}_{i,j+\frac{1}{2}} - \hat{G}_{i,j-\frac{1}{2}} \Big),
	\label{WENO_Conservative_Form}
\end{equation}
where $\hat{F}_{i+\frac{1}{2},j}$, $\hat{F}_{i-\frac{1}{2},j}$, $\hat{G}_{i,j+\frac{1}{2}}$ and $\hat{G}_{i,j-\frac{1}{2}}$ are the numerical fluxes at $(x_{i+\frac{1}{2}}, y_{j})$, $(x_{i-\frac{1}{2}}, y_{j})$, $(x_{i}, y_{j+\frac{1}{2}})$ and $(x_{i}, y_{j-\frac{1}{2}})$, respectively. Also in \eqref{WENO_Conservative_Form}, for compact notation, $\hat{F}$ and $\hat{G}$ denotes the $x$ and $y$ components of the vector field $\F(\cdot)$, respectively. 

For brevity, we will only focus on defining $\hat{F}_{i+\frac{1}{2},j}$ and $\hat{F}_{i-\frac{1}{2},j}$, and note that \dt{corresponding  calculations for $\hat{G}_{i, j+\frac{1}{2}}$ and $\hat{G}_{i, j-\frac{1}{2}}$ follows identical steps.} In this context, we split the $x$ component of $\F(\cdot)$ into two parts
\begin{equation*}
	\hat{F}(\bu) = \hat{F}^{+}(\bu) + \hat{F}^{-}(\bu),
\end{equation*}
where $d\hat{F}^{+}(\bu) / d \bu > 0$ and $d\hat{F}^{-}(\bu) / d \bu \leq 0$. Then, we define these parts by using the popular Rusanov-type flux splitting method \cite{Rusanov1962} that is given by
\begin{equation}
	\hat{F}^{\pm}(\bu) = \dfrac{1}{2} \Big( \hat{F}(\bu) \pm \alpha \bu \Big),
	\label{Rusanov_Flux_Splitting}
\end{equation}
where we approximate $\alpha := \max | d\hat{F}(\bu) / d \bu |$, \emph{i.e.,} the spectral radius of the Jacobian generated by $\hat{F}(\cdot)$ in a Jacobian-free manner, detailed in Section~\ref{Approximation_of_Propagation_Speed}. To this end, let us denote by $\hat{F}^{+}_{i \pm \frac{1}{2},j}$ and $\hat{F}^{-}_{i \pm \frac{1}{2},j}$ the numerical fluxes obtained by splitting $\hat{F}(\cdot)$ into the positive and negative parts, respectively. Then these numerical fluxes used in \eqref{WENO_Conservative_Form} are given by the sum of their associated parts
\begin{equation*}
	\hat{F}_{i \pm \frac{1}{2},j} = \hat{F}^{+}_{i \pm \frac{1}{2},j} + \hat{F}^{-}_{i \pm \frac{1}{2},j},
\end{equation*}
where following the standard ZSWENO procedure, $\hat{F}^{+}_{i \pm \frac{1}{2},j}$ and $\hat{F}^{-}_{i \pm \frac{1}{2},j}$ are given by the weighted combination of the three third-order \emph{essentially non-oscillatory} (ENO) approximations \cite{Jiang1996,Kim2005,Liu1994,Zhang2006}. Hence, they are given by 
\begin{equation}
	\hat{F}^{\pm}_{i+\frac{1}{2},j} = \sum_{k=0}^{2} \omega^{\pm}_{k,i+\frac{1}{2},j} \hat{F}^{\pm}_{k,+}, \qquad
	\hat{F}^{\pm}_{i-\frac{1}{2},j} = \sum_{k=0}^{2} \omega^{\pm}_{k,i-\frac{1}{2},j} \hat{F}^{\pm}_{k,-},
	\label{WENO_Fluxes}
\end{equation}
where $\omega^{\pm}_{k,i\pm\frac{1}{2},j}$ are the non-linear weights that we will define later, $\hat{F}^{\pm}_{k,+}$ are the ENO approximations given by
\begin{equation}
	\begin{split}
		\hat{F}^{+}_{0,+} = & \dfrac{1}{3} \hat{F}^{+}(\bu_{i-2,j}) - \dfrac{7}{6} \hat{F}^{+}(\bu_{i-1,j}) + \dfrac{11}{6} \hat{F}^{+}(\bu_{i,j}), \\
		\hat{F}^{+}_{1,+} = & - \dfrac{1}{6} \hat{F}^{+}(\bu_{i-1,j}) + \dfrac{5}{6} \hat{F}^{+}(\bu_{i,j}) + \dfrac{1}{3} \hat{F}^{+}(\bu_{i+1,j}), \\
		\hat{F}^{+}_{2,+} = & \dfrac{1}{3} \hat{F}^{+}(\bu_{i,j}) + \dfrac{5}{6} \hat{F}^{+}(\bu_{i+1,j}) - \dfrac{1}{6} \hat{F}^{+}(\bu_{i+2,j}), \\
		\hat{F}^{-}_{0,+} = & \dfrac{1}{3} \hat{F}^{-}(\bu_{i+3,j}) - \dfrac{7}{6} \hat{F}^{-}(\bu_{i+2,j}) + \dfrac{11}{6} \hat{F}^{-}(\bu_{i+1,j}), \\
		\hat{F}^{-}_{1,+} = & - \dfrac{1}{6} \hat{F}^{-}(\bu_{i+2,j}) + \dfrac{5}{6} \hat{F}^{-}(\bu_{i+1,j}) + \dfrac{1}{3} \hat{F}^{-}(\bu_{i,j}), \\
		\hat{F}^{-}_{2,+} = & \dfrac{1}{3} \hat{F}^{-}(\bu_{i+1,j}) - \dfrac{5}{6} \hat{F}^{-}(\bu_{i,j}) + \dfrac{1}{6} \hat{F}^{-}(\bu_{i-1,j}),
	\end{split}
	\label{ENO_Fluxes_Plus}
\end{equation}
and similarly, $\hat{F}^{\pm}_{k,-}$ are defined by
\begin{equation}
	\begin{split}
		\hat{F}^{+}_{0,-} = & \dfrac{1}{3} \hat{F}^{+}(\bu_{i-3,j}) - \dfrac{7}{6} \hat{F}^{+}(\bu_{i-2,j}) + \dfrac{11}{6} \hat{F}^{+}(\bu_{i-1,j}), \\
		\hat{F}^{+}_{1,-} = & - \dfrac{1}{6} \hat{F}^{+}(\bu_{i-2,j}) + \dfrac{5}{6} \hat{F}^{+}(\bu_{i-1,j}) + \dfrac{1}{3} \hat{F}^{+}(\bu_{i,j}), \\
		\hat{F}^{+}_{2,-} = & \dfrac{1}{3} \hat{F}^{+}(\bu_{i-1,j}) + \dfrac{5}{6} \hat{F}^{+}(\bu_{i,j}) - \dfrac{1}{6} \hat{F}^{+}(\bu_{i+1,j}), \\
		\hat{F}^{-}_{0,-} = & \dfrac{1}{3} \hat{F}^{-}(\bu_{i+2,j}) - \dfrac{7}{6} \hat{F}^{-}(\bu_{i+1,j}) + \dfrac{11}{6} \hat{F}^{-}(\bu_{i,j}), \\
		\hat{F}^{-}_{1,-} = & - \dfrac{1}{6} \hat{F}^{-}(\bu_{i+1,j}) + \dfrac{5}{6} \hat{F}^{-}(\bu_{i,j}) + \dfrac{1}{3} \hat{F}^{-}(\bu_{i-1,j}), \\
		\hat{F}^{-}_{2,-} = & \dfrac{1}{3} \hat{F}^{-}(\bu_{i,j}) + \dfrac{5}{6} \hat{F}^{-}(\bu_{i-1,j}) - \dfrac{1}{6} \hat{F}^{-}(\bu_{i-2,j}).
	\end{split}
	\label{ENO_Fluxes_Minus}
\end{equation}
Finally, in \eqref{ENO_Fluxes_Plus} and \eqref{ENO_Fluxes_Minus} $\hat{F}^{\pm}(\bu)$ are given by the Rusanov-type flux splitting method defined in \eqref{Rusanov_Flux_Splitting}.

Since our aim is to reduce the computational cost and so for this we seek to use convolution, we observe that indeed these ENO fluxes \eqref{ENO_Fluxes_Plus} and \eqref{ENO_Fluxes_Minus} can be equivalently represented by discrete convolutions, \emph{i.e.,}
\begin{equation}
	\hat{F}^{\pm}_{k,+} = \widetilde{\K}^{\pm}_{k,+} \ast \F^{\pm}(\bu), \qquad
	\hat{F}^{\pm}_{k,-} = \widetilde{\K}^{\pm}_{k,-} \ast \F^{\pm}(\bu),
	\label{ENO_Fluxes_Convolution}
\end{equation}
where $\widetilde{\K}^{\pm}_{k,+}$ and $\widetilde{\K}^{\pm}_{k,-}$ are the induced vectors from \eqref{ENO_Fluxes_Plus} and \eqref{ENO_Fluxes_Minus}, and for completeness they are defined in Appendix \ref{WENO_Convolution_Vectors}.

Let us now shift our attention to the non-linear weights $\omega^{\pm}_{k,i\pm\frac{1}{2},j}$ that we used to construct the ZSWENO approximation in \eqref{WENO_Fluxes}. Following again the usual procedure \cite{Liu1994}, we define these weights as
\begin{equation}
	\begin{split}
		\omega^{\pm}_{k,i+\frac{1}{2},j} := \dfrac{\alpha^{\pm}_{k,i+\frac{1}{2},j}}{\sum_{k=0}^{2} \alpha^{\pm}_{k,i+\frac{1}{2},j}}, \qquad
		\alpha^{\pm}_{k,i+\frac{1}{2},j} := \dfrac{\omega^{o}_{k,+}}{\Big( \epsilon_{_{W}} + IS^{\pm}_{k,i+\frac{1}{2},j} \Big)^{p}}, \\
		\omega^{\pm}_{k,i-\frac{1}{2},j} := \dfrac{\alpha^{\pm}_{k,i-\frac{1}{2},j}}{\sum_{k=0}^{2} \alpha^{\pm}_{k,i-\frac{1}{2},j}}, \qquad
		\alpha^{\pm}_{k,i-\frac{1}{2},j} := \dfrac{\omega^{o}_{k,-}}{\Big( \epsilon_{_{W}} + IS^{\pm}_{k,i-\frac{1}{2},j} \Big)^{p}}, \\
	\end{split}
	\label{Non_Linear_Weights_WENO}
\end{equation}
where we take the usual $p = 2$, $\epsilon_{_{W}}=10^{-6}$ values \cite{Jiang1996}, and define the optimal weights $\omega^{o}_{k,+}$ and $\omega^{o}_{k,-}$ as: $\omega^{o}_{0,+}=0.1$, $\omega^{o}_{1,+}=0.6$, $\omega^{o}_{2,+}=0.3$ and $\omega^{o}_{0,-}=0.3$, $\omega^{o}_{1,-}=0.6$, $\omega^{o}_{2,-}=0.1$. Furthermore, in \eqref{Non_Linear_Weights_WENO} $IS^{\pm}_{k,i\pm\frac{1}{2},j}$ are the ZSWENO smoothness indicators \cite{Zhang2006}, and the explicit formulae first for $IS^{\pm}_{k,i+\frac{1}{2},j}$ are given by
\begin{equation}
	\begin{split}
		IS^{+}_{0,i+\frac{1}{2},j} = & \big( \hat{F}^{+}(\bu_{i-2,j}) - 4 \hat{F}^{+}(\bu_{i-1,j}) + 3 \hat{F}^{+}(\bu_{i,j}) \big)^{2}, \\
		IS^{+}_{1,i+\frac{1}{2},j} = & \big( \hat{F}^{+}(\bu_{i-1,j}) - \hat{F}^{+}(\bu_{i+1,j}) \big)^{2}, \\
		IS^{+}_{2,i+\frac{1}{2},j} = & \big( 3 \hat{F}^{+}(\bu_{i,j}) - 4 \hat{F}^{+}(\bu_{i+1,j}) + \hat{F}^{+}(\bu_{i+2,j}) \big)^{2}, \\
		IS^{-}_{0,i+\frac{1}{2},j} = & \big( \hat{F}^{-}(\bu_{i+3,j}) - 4 \hat{F}^{-}(\bu_{i+2,j}) + 3 \hat{F}^{-}(\bu_{i+1,j}) \big)^{2}, \\
		IS^{-}_{1,i+\frac{1}{2},j} = & \big( \hat{F}^{-}(\bu_{i+2,j}) - \hat{F}^{-}(\bu_{i,j}) \big)^{2}, \\
		IS^{-}_{2,i+\frac{1}{2},j} = & \big( 3\hat{F}^{-}(\bu_{i+1,j}) - 4 \hat{F}^{-}(\bu_{i,j}) + \hat{F}^{-}(\bu_{i-1,j}) \big)^{2},
	\end{split}
	\label{ZSWENO_Smoothness_Indicators_Plus}
\end{equation}
and then for $IS^{\pm}_{k,i-\frac{1}{2},j}$, we have
\begin{equation}
	\begin{split}
		IS^{+}_{0,i-\frac{1}{2},j} = & \big( \hat{F}^{+}(\bu_{i-3,j}) - 4 \hat{F}^{+}(\bu_{i-2,j}) + 3 \hat{F}^{+}(\bu_{i-1,j}) \big)^{2}, \\
		IS^{+}_{1,i-\frac{1}{2},j} = & \big( \hat{F}^{+}(\bu_{i-2,j}) - \hat{F}^{+}(\bu_{i,j}) \big)^{2}, \\
		IS^{+}_{2,i-\frac{1}{2},j} = & \big( 3 \hat{F}^{+}(\bu_{i-1,j}) - 4 \hat{F}^{+}(\bu_{i,j}) + \hat{F}^{+}(\bu_{i-1,j}) \big)^{2}, \\
		IS^{-}_{0,i-\frac{1}{2},j} = & \big( \hat{F}^{-}(\bu_{i+2,j}) - 4 \hat{F}^{-}(\bu_{i+1,j}) + 3 \hat{F}^{-}(\bu_{i,j}) \big)^{2}, \\
		IS^{-}_{1,i-\frac{1}{2},j} = & \big( \hat{F}^{-}(\bu_{i+1,j}) - \hat{F}^{-}(\bu_{i-1,j}) \big)^{2}, \\
		IS^{-}_{2,i-\frac{1}{2},j} = & \big( 3 \hat{F}^{-}(\bu_{i,j}) - 4 \hat{F}^{-}(\bu_{i-1,j}) + \hat{F}^{-}(\bu_{i-2,j}) \big)^{2}.
	\end{split}
	\label{ZSWENO_Smoothness_Indicators_Minus}
\end{equation}
Once again, in \eqref{ZSWENO_Smoothness_Indicators_Plus} and \eqref{ZSWENO_Smoothness_Indicators_Minus} the positive $\hat{F}^{+}(\bu)$ and negative $\hat{F}^{-}(\bu)$ parts are defined by using the Rusanov flux splitting \eqref{Rusanov_Flux_Splitting}. \dt{Furthermore, we can} observe that these smoothness indicators \eqref{ZSWENO_Smoothness_Indicators_Plus} and \eqref{ZSWENO_Smoothness_Indicators_Minus} can be equivalently expressed in terms of discrete convolutions, \emph{i.e.,}
\begin{equation}
	IS^{\pm}_{k,i+\frac{1}{2},j} = \big( \overline{\K}^{\pm}_{k,+} \ast \F^{\pm}(\bu) \big)^{2}, \qquad
	IS^{\pm}_{k,i-\frac{1}{2},j} = \big( \overline{\K}^{\pm}_{k,-} \ast \F^{\pm}(\bu) \big)^{2},
	\label{ZSWENO_Smoothness_Indicators_Convolutions}
\end{equation}
where the appropriately induced vectors $\overline{\K}^{\pm}_{k,+}$ and $\overline{\K}^{\pm}_{k,-}$ are defined in \ref{WENO_Convolution_Vectors}. This completes the description of the convolution-driven ZSWENO scheme for the inside differential operators $\nabla \cdot \F(\bu)^{I}$, \emph{i.e.,} the approximation of $\nabla \cdot \F(\bu)$ for all inside nodes.

However, since the adhesion operators $\nabla \cdot \F(\bu)$ were split in two parts in \eqref{Splitted_Adhesion}, it remains to describe the ZSWENO scheme for the layer operator $\nabla \cdot \F(\bu)^{L}$, accounting for all nodes that are considered to be in the layer part of the domain. As in Section~\ref{Diffusion_Operators_Numeics} for the boundary diffusion operators, here we also need to appropriately approximate the value of any point that is located outside of the tumour domain $\Omega(t_{0})$. To this end, we symmetrically reflect the values of the interior nodes to the ghost nodes in a node by node fashion due to the irregular tumour domain. This ultimately enables us to use the standard non-convolutional ZSWENO scheme \eqref{ENO_Fluxes_Plus}, \eqref{ENO_Fluxes_Minus}, \eqref{ZSWENO_Smoothness_Indicators_Plus} and \eqref{ZSWENO_Smoothness_Indicators_Minus} instead of the convolutional forms \eqref{ENO_Fluxes_Convolution} and \eqref{ZSWENO_Smoothness_Indicators_Convolutions} to approximate the layer operators $\nabla \cdot \F(\bu)^{L}$.

%--------------------------------------------------------------
%				  	 Propagation Speed
%--------------------------------------------------------------
\subsection{Approximation of the Propagation Speed $\alpha$}\label{Approximation_of_Propagation_Speed}
Due to the complexity of $\F(\cdot)$, calculating its Jacobian $d\hat{F}(\bu) / d \bu$ is extremely time consuming. Hence, to find the largest eigenvalue $\alpha$, we rather skip the exact evaluation of the Jacobian and choose to approximate the propagation speed in a Jacobian-free manner. For this, let us first define the eigenvalue problem by
\begin{equation}
	J v = \lambda v,
	\label{Eigenvalue_Problem}
\end{equation}
where $J = d\F(\bu) / d \bu$ is the Jacobian, $v$ is an eigenvector, and $\lambda$ is an eigenvalue of $J$. Then in order to find the largest eigenvalue $\alpha$, we solve \eqref{Eigenvalue_Problem} iteratively with convergence tolerance of $10^{-14}$ and with a random initial guess for $v$. Hence, similarly to for instance a Jacobian-free Newton-Krylov method \cite{Kelley2003,Knoll2004}, to find $\lambda v$ in \eqref{Eigenvalue_Problem} we evaluate the Jacobian-vector product $J v$ instead of $J$, which is proved to be a significantly less time-consuming task. To that end, the approximation of $J v$ is carried out via the first-order Taylor series expansion, and so it is given by
\begin{equation*}
	J v \approx \dfrac{\F(\bu + \epsilon_{p} v) - \F(\bu)}{\epsilon_{p}},
\end{equation*}
where $\epsilon_{p}$ is a small perturbation parameter. Since the precision is limited in the evaluation of the flux $\F(\cdot)$, a good choice to evaluate this small parameter $\epsilon_{p}$ is given by \cite{Knoll2004}
\begin{equation*}
	\epsilon_{p} = \dfrac{\sqrt{(1 + \nor{\bu})\epsilon_{mach}}}{\nor{v}},
\end{equation*}
where $\epsilon_{mach}>0$ is the machine precision.

%--------------------------------------------------------------
%	  						Results
%--------------------------------------------------------------
\section{Numerical Results}\label{section:Numerical_Results}
\begin{figure*}
	\centering
      \includegraphics[width=1\textwidth]{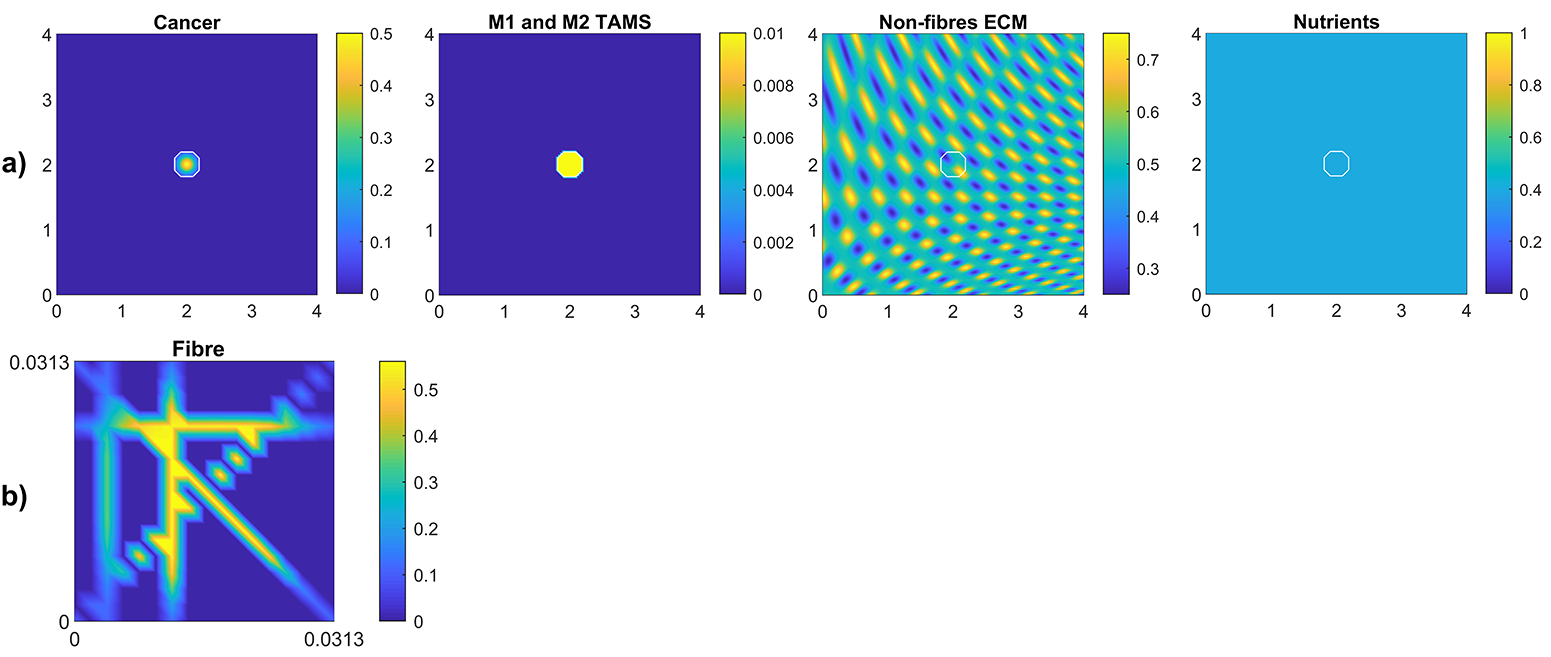}
    \caption{a) Initial conditions of the cancer cells, both M1 and M2 TAMs, non-fibre ECM and the nutrients, respectively. b) Initial condition of the micro-fibres on $\delta Y(x)$, which is repeated for all $x$.}
    \label{fig:Initial_Conditions}
\end{figure*}

%--------------------------------------------------------------
%			    Intro for the Numerical Results
%--------------------------------------------------------------
In this section, we present the numerical results for our multi-scale model. Hence, for the simulations let us consider a tissue domain $Y = [ 0 \times 4 ] \times [ 0 \times 4 ]$ and the following initial conditions:
\begin{equation}\label{initalConditionMacroSystem17Dec2020}
	\begin{split}
		c(x, 0) = & \dfrac{1}{2} \exp\bigg( \dfrac{-\nor{x}_{2}^{2}}{0.02} \bigg) \cdot \chi_{_{\Bila((2,2),0.25)}}, \\
        M_{1}(x, 0) = & 10^{-2} \cdot \chi_{_{\Bila((2,2),0.25)}}, \\
        M_{2}(x, 0) = & 10^{-2} \cdot \chi_{_{\Bila((2,2),0.25)}}, \\
        l(x, 0) = & \min \bigg( \dfrac{1}{2} + \dfrac{1}{4} \sin(7 \pi x_{1}x_{2})^{3} \cdot \sin \bigg(7 \pi \dfrac{x_{2}}{x_{1}} \bigg), 1 - c(x, 0) \bigg ), \\
        \sigma(x, 0) = & 0.4;
	\end{split}
\end{equation}
which are illustrated in Fig.~\ref{fig:Initial_Conditions}a). Here the white curves indicate the boundary of the tumour domain $\partial \Omega(0)$. Besides these macro-scale initial conditions, we also illustrate the initial state of one micro-scale fibre domain $\delta Y(x)$ in Fig.~\ref{fig:Initial_Conditions}b), which is repeated for all macro-scale points. Also, the ratio between the fibre and non-fibre ECM phases are assumed to be $20\%$ : $80\%$ for all simulations.

Finally, all presented simulations corresponds to time $50 \Delta t$. The baseline parameters values are provided in Appendix \ref{Parameter_Set}, and any alteration from these values will be stated accordingly.
%------------ Figure for the initial conditions
%------------

%--------------------------------------------------------------
%			   	     Spatial dependent results
%--------------------------------------------------------------
\subsection{Spatial Dependency of the Re-polarisation}
First, we investigate numerically the effects of changing the re-polarisation domain $\Omega_{p}(t, R_{p})$ used in \eqref{Macrophage_Re_Polarisation}, defined in Appendix \ref{Appendix_Re_Polarisation_Domain} and illustrated in Fig.~\ref{fig:Re_Polarisation_Domain}. Hence, here we study whether the success of a M2$\to$M1 re-polarisation strategy against the tumour is dependent on the spatial domain $\Omega_{p}(t, R_{p})$, specifically on $R_{p}$=the distance from the outer boundary $\partial \Omega_{o}(t)$. For this let us use the radii $R_{p} \in \{ 0, \Delta x, 2\Delta x, 3\Delta x, 4\Delta x \}$ for the re-polarisation domain $\Omega_{p}(t, R_{p})$, as well as let us start the process at time $0$, \emph{i.e.,} we take $t_{p} = 0$ in \eqref{Macrophage_Re_Polarisation}. To further study these spatial effects, we also consider multiple tissue conditions by changing the tissue environment controller $\beta \in \{ 0.75, 0.7875, 0.825 \}$ (see \ref{Appendix_MDE_Micro_Process} or Ref~\cite{Dumitru_et_al_2013}).

To compare the resulting tumours, we measure the overall tumour mass and spread at final time $50 \Delta t$, using the total density of the cancer cells as well as the area of the tumour. Ultimately, this enables us to quantify the changes, resulted by the modification of the re-polarisation domain $\Omega_{p}(t, R_{p})$. Specifically, the overall tumour mass is measured by integrating the cancer cell density $c(x,50 \Delta t)$ and the overall tumour spread as given by the area of the tumour domain $\Omega(50 \Delta t)$.

\begin{figure*}[!hp]
	\centering
      \includegraphics[width=0.75\textwidth]{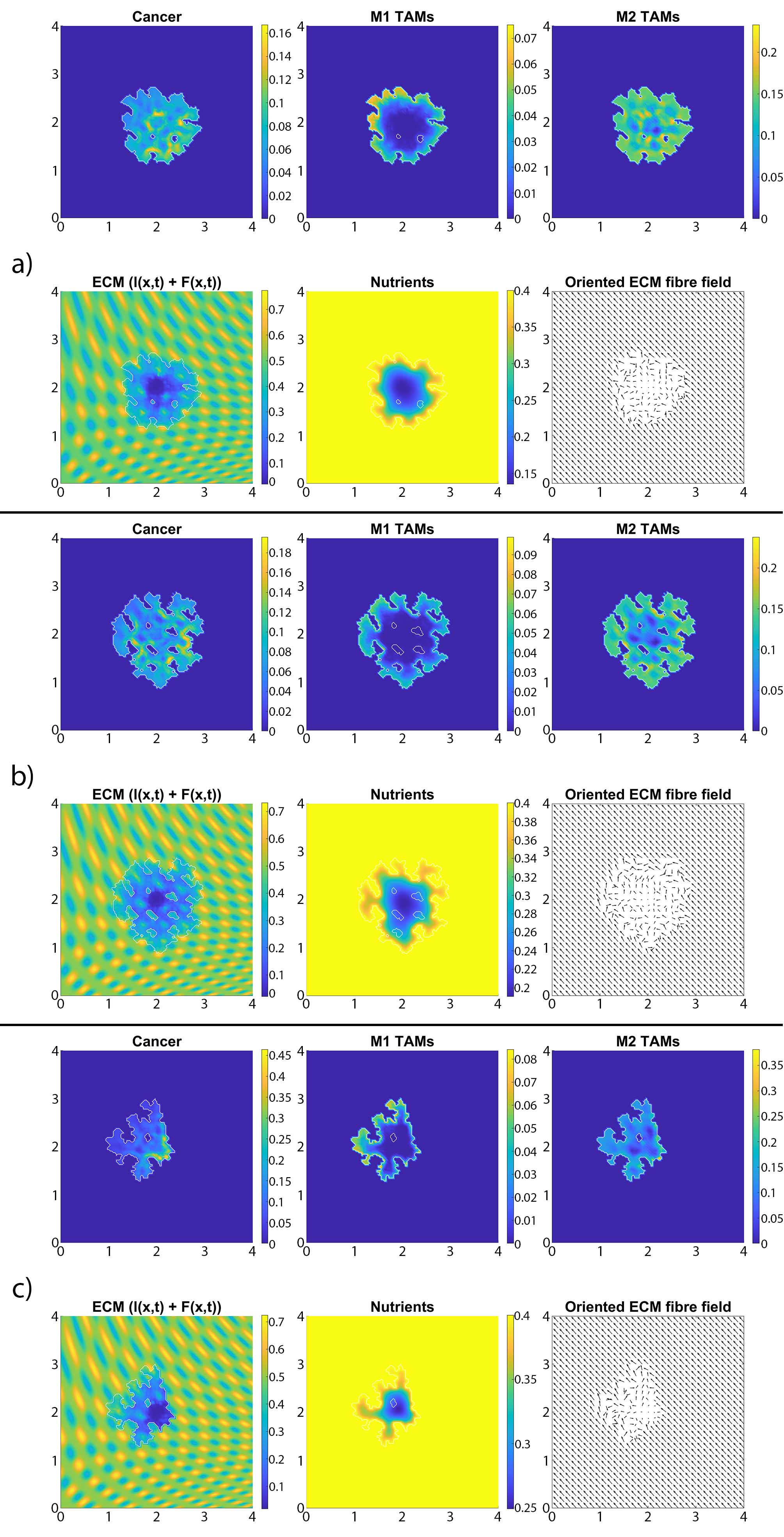}
    \caption{Baseline simulation with no macrophage re-polarisation (\emph{i.e.,} we set $R_{p} = 0$) at final time $50 \Delta t$. a) $\beta = 0.75$, b) $\beta =0.7875$, c) $\beta =0.825$.}
    \label{fig:No_Repolarisation}
\end{figure*}

\paragraph{Baseline cases:} In Fig.~\ref{fig:No_Repolarisation} \re{we present the numerical results for the baseline case characterised by the absence of M2$\to$M1} re-polarisation. Specifically, in Fig.~\ref{fig:No_Repolarisation} a)-c) we \re{investigate the no-repolarisation case in the context of different} tissue environment controllers: a) $\beta = 0.75$, b) $\beta =0.7875$, and c) $\beta =0.825$.
%------------ Figure for the baseline simulation with no re-polarisation
%------------
As we can see, at time $50 \Delta t$, the tumour is enlarged and has invaded some of its surroundings, and within the tumour domain $\Omega(50 \Delta t)$ we observed heterogeneity in all three cell populations (cancer cell, M1 and M2 TAMs). We also notice that near M2 TAMs accumulations sites, the density of cancer cells also seems to be higher. In contrast, higher M1 TAMs density usually corresponds to both low M2 TAMs as well as low cancer cell populations. \re{The movement and spatial distribution of tumour and immune cells is influenced directly and indirectly by the} ECM fibres. For illustrative purposes, the rearranged fibre structure is portrayed by a four-fold coarsened oriented ECM fibres field in Fig.~\ref{fig:No_Repolarisation}. The peritumoral degradation of the two-phase ECM (caused by the cancer cells and both TAMs) allows the tumour to expand and spread to the neighbouring tissues, resulting in some tumour fingering patterns that vary with the controller $\beta$: higher $\beta$ leads to more tumour fingering. This induces an irregular tumour domain, as well as the formation of "islands" inside the tumour, which \re{corresponds} to areas of initially low ECM density, \emph{i.e.,} the ECM level was too low in such areas to support tumour movement. Finally, in Fig.~\ref{fig:No_Repolarisation}, we also present the level of available nutrients. Hence, we can see that since the nutrients are only supplied through the outer boundary $\Omega_{o}(t)$, the initial normal level of nutrients is significantly depleted (by the three cell populations) inside the tumour. This can lead to hypoxia and then eventually create a necrotic core \re{(not modelled explicitly in this study)}.

\begin{figure*}
	\centering
      \includegraphics[width=0.75\textwidth]{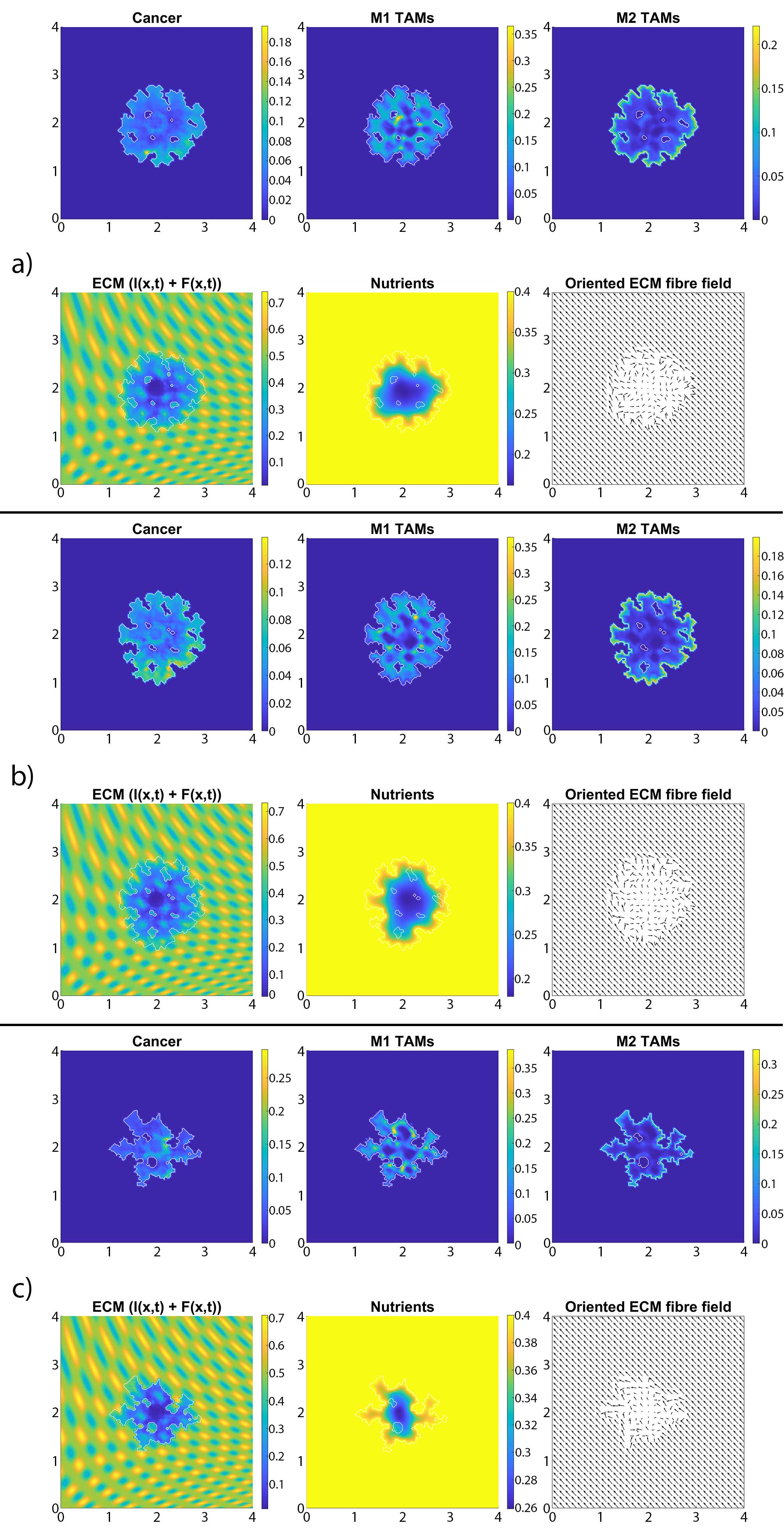}
    \caption{Simulation with re-polarisation domain $\Omega_{p}(t, \Delta x)$ and with starting time $t_{p} = 0$ at the final time $50 \Delta t$. a) $\beta = 0.75$, b) $\beta =0.7875$, c) $\beta =0.825$.}
    \label{fig:Repolarisation}
\end{figure*}

\paragraph{The impact of M2$\to$M1 re-polarisation:} In Fig.~\ref{fig:Repolarisation} a)-c), we introduce the re-polarisation of the M2 TAMs to M1 TAMs within the domain $\Omega_{p}(t, \Delta x)$ (\emph{i.e.,} $R_{p} = \Delta x$), and again we investigate the role of three controller values: a) $\beta = 0.75$, b) $\beta =0.7875$, c) $\beta =0.825$.
%------------ Figure for R_{p} = \delta x
%------------
%We first discuss the two macrophage populations and highlight the changes from the previous baseline case (in Fig.~\ref{fig:No_Repolarisation}) resulted by the introduction of the re-polarisation. 
While in Fig.~\ref{fig:No_Repolarisation} the M1 TAMs mainly accumulated near the tumour boundary, in Fig.~\ref{fig:Repolarisation} we see \re{that the re-polarisation of M2 TAMs leads to} an increase in the M1 TAMs population inside the tumour, leading to several accumulation sites located further away from the leading edge. In contrast, the presence of M2 TAMs inside of the tumour is decreased compared to Fig.~\ref{fig:No_Repolarisation}, and these immune cells now accumulate only in a small neighbourhood of the boundary \re{(because we do not re-polarise M2 TAMs into M1 TAMs in a $R_{p} =\Delta x$ neighbourhood of the boundary; see Fig~\ref{fig:Re_Polarisation_Domain}).} Since both macrophage populations interact with the cancer cells, we also see some differences in the cancer cell population \re{and these differences depend on the degradation level of ECM as controlled by the value of $\beta$}. %Consequently, the maximal cancer cell density is slightly reduces as well as these maximums occurs mainly at the tumour interface. Therefore, here we see again a harmony between the cancer cell and the M2 TAMs populations as well as we see that the increased M1 TAMs population reduced the overall amount of cancer cells and restricted their accumulation sites to a small $\Delta x$ neighbourhood of the tumour boundary. 
Using two measures (spread and mass), \re{we first focus on} the $\beta = 0.75$ case and compare Figs.~\ref{fig:No_Repolarisation} a) and \ref{fig:Repolarisation} a): \re{re-polarisation leads to} an $\approx 11\%$ increase in tumour spread as well as an $\approx 20\%$ decrease in tumour mass. Then for $\beta = 0.7875$ (Figs.~\ref{fig:No_Repolarisation} b) and \ref{fig:Repolarisation} b)) we observe an $\approx 5\%$ reduction in spread and an $\approx 34\%$ decrease in mass. Finally, for  $\beta = 0.825$ (Figs.~\ref{fig:No_Repolarisation} c) and \ref{fig:Repolarisation} c)), the tumour spread is increased by $\approx 15\%$ and the tumour mass is reduced by $\approx 31\%$. Hence, although re-polarising M2 TAMs into M1 TAMs does not show an overall improvement in terms of tumour spread, it significantly impeded the migration of the cancer cell mass from hypoxic regions to the proliferating rim as well as helped reducing the mass by killing the cancer cells. However, since $R_{p} = \Delta x$, the M2 TAMs are able to accumulate in a $\Delta x$ neighbourhood of the interface and so we can still observe a moderate density of cancer cells near the boundary.

\begin{figure*}
      \includegraphics[width=1\textwidth]{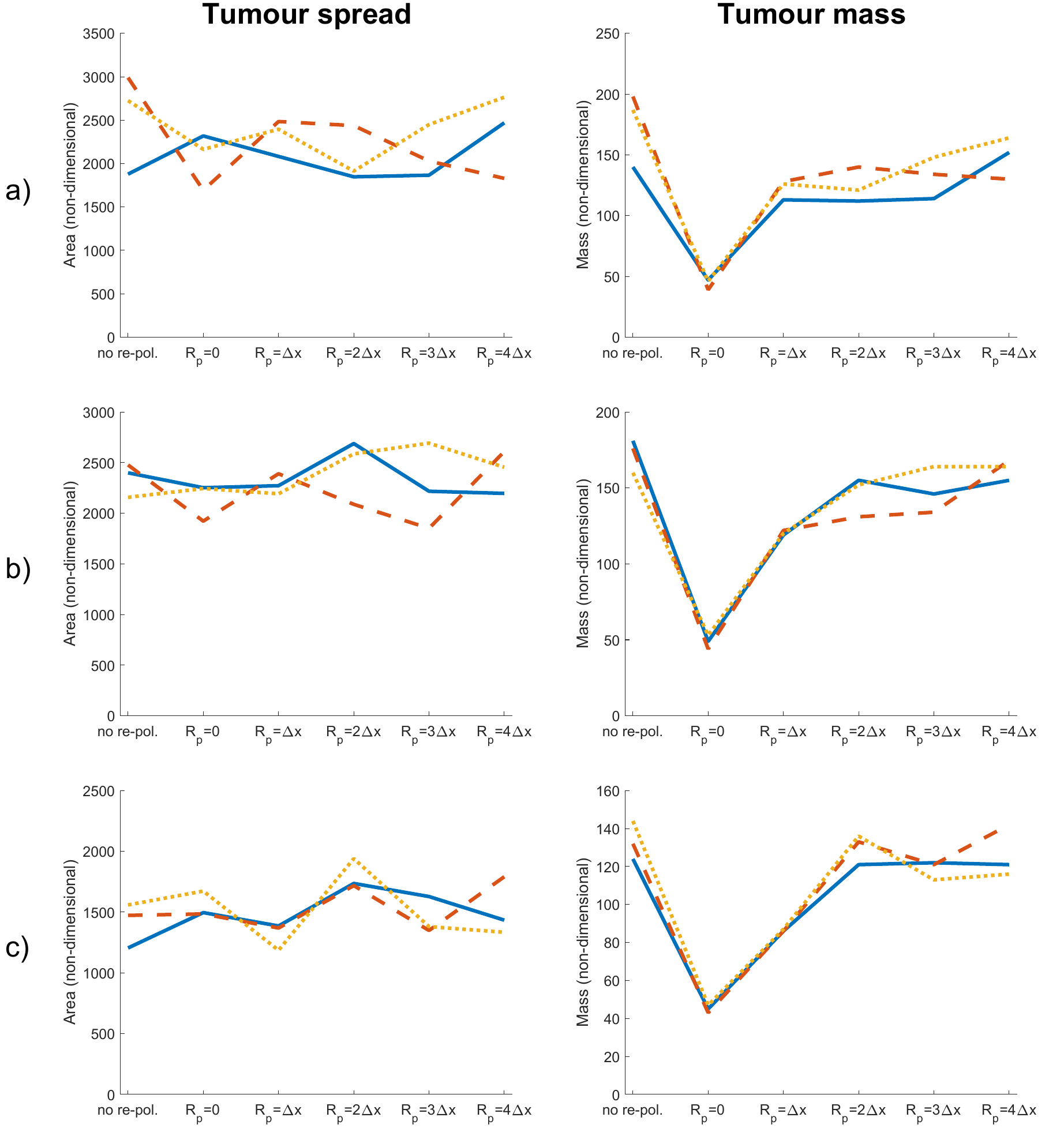}
    \caption{Results of varying the radius of the re-polarisation domain $R_{p} \in \{ 0, \Delta x, 2\Delta x, 3\Delta x, 4\Delta x \}$ with starting time $t_{p} = 0$. Left panel shows the change in the tumour spread and the right panel shows the change in tumour mass with respect to $R_{p}$. Here the solid-blue lines correspond to cancer cells with high MMP secretion rate ($\beta_{lc} = 3.0$, $\beta_{Fc} = 1.5$, $\alpha_{c} = 0.625$), the dashed-orange lines to a type with medium secretion rate ($\beta_{lc} = 2.0$, $\beta_{Fc} = 1.0$, $\alpha_{c} = 0.42$) and the dotted-yellow lines to low secretion rate type cancer ($\beta_{lc} = 1.0$, $\beta_{Fc} = 0.5$, $\alpha_{c} = 0.21$). The environment controller $\beta$ is set to a) $\beta = 0.75$, b) $\beta =0.7875$, c) $\beta =0.825$.}
    \label{fig:Spatial_Repolarisation_All}
\end{figure*}

\paragraph{Tumour spread/mass changes with respect to $R_{p}$:} In Fig.~\ref{fig:Spatial_Repolarisation_All} we vary $R_{p}$ and present the changes in the dimensionless tumour area (spread) and mass resulted by changing the radius $R_{p}$. 
%------------ Figure for the spread/mass of the tumour domain with varying distances
%------------
Specifically, in Fig.~\ref{fig:Spatial_Repolarisation_All} a)-c) we again consider the three previously used environment controllers $\beta = 0.75$, $\beta =0.7875$, $\beta =0.825$. Moreover, since various types of cancers secrete MMPs at different rates \cite{Gobin2019}, here we also consider high, medium and low MDE secretion rates. Finally, we vary the radii for the re-polarisation domain, $R_{p} \in \{ 0, \Delta x, 2\Delta x, 3\Delta x, 4\Delta x \}$, while we assume that the re-polarisation process starts at time $0$, \emph{i.e.,} we take $t_{p} = 0$. The left panels of Fig.~\ref{fig:Spatial_Repolarisation_All} show the changes in the dimensionless variable for tumour areas/spreads with respect to $R_{p}$, while the right panels show the changes in the dimensionless variables for tumour mass with respect to $R_{p}$. Here, the solid-blue lines correspond to high MMP secretion rates ($\beta_{lc} = 3.0$, $\beta_{Fc} = 1.5$, $\alpha_{c} = 0.625$), the dashed-orange lines correspond to medium secretion rates ($\beta_{lc} = 2.0$, $\beta_{Fc} = 1.0$, $\alpha_{c} = 0.42$), and the dotted-yellow lines correspond to low MMP secretion rates ($\beta_{lc} = 1.0$, $\beta_{Fc} = 0.5$, $\alpha_{c} = 0.21$). Comparing first the changes due to varying the tissue controller $\beta$, we can see a clear overall decrease in both tumour spread and mass as we increase $\beta$. This was seen also in Figs.\ref{fig:No_Repolarisation} and \ref{fig:Repolarisation}, where a more prominent tumour fingering pattern was present as we increased $\beta$ which resulted in a decrease in tumour spread. Furthermore, in Fig.~\ref{fig:Spatial_Repolarisation_All} a)-c) the overall behaviour of the tumour spread does not change significantly as we vary the radius of the re-polarisation domain $R_{p}$. Hence, even though the proteolytic molecular processes at the leading tumour edge change (via the MDE source \eqref{MDE_Source}) following the M2$\to$M1 re-polarisation, we cannot see a clear trend in tumour spread. This might be because of the complex interactions between tumour and infiltrating immune cell populations: although we see an increase in the M1 TAMs populations near the boundary (that also secrete more MDEs than the M2 cells \cite{Hedbrant2015}), the overall proteolytic molecular processes at the leading edge of the tumour might not change too much, which would mean similar tumour spread. 
Second, the rearrangement of the micro-fibre distribution also affects tumour spread, since the amount of fibres that are being re-located near the leading edge is dependent on the fluxes generated by the different cell populations. Therefore, our model suggest that merely re-polarising the M2 TAMs into the anti-tumoral M1 phenotype might not be enough to stop tumour spread.\\
 On the other hand, in Fig.~\ref{fig:Spatial_Repolarisation_All} a)-c) we see that the tumour mass is greatly reduced in all of the presented cases. These results also show that for the \re{largest} reduction we need to re-polarise the M2 TAMs inside the whole tumour domain (\emph{i.e.,} $R_{p}=0$, which means that the re-polarisation domain is $\Omega_{p}(t, 0) = \Omega(t)$). Hence, the presence of M2 TAMs in the proliferating rim may be enough to maintain the tumour mass and to induce tumour spread \re{by leading to a moderate presence of cancer cells in the proliferating rim}. 
 %This is particularly important since we considered a nutrient dependent cancer death term, meaning that the alive cancer cells that are mainly located in the hypoxic regions and the proliferating rim, have been brought to the fore. Thereby, the location of these cells become more significant and our model shows that an adequate density of the M2-like macrophages near the boundary (when we consider $R_{p} > 0$) can in fact promote a moderate amount cancer cells in the proliferating rim. This on the one hand, results in a higher tumour mass (compared to $R_{p} = 0$) and on the other hand, it seems to be sufficient for sustaining the tumour spread and growth.\\
 In conclusion, the results in Fig.~\ref{fig:Spatial_Repolarisation_All} suggest that re-polarising M2-like macrophages to the M1 phenotype has the ability to reduce tumour mass, but we may need a supplementary strategy primarily focusing on tumour spread or tumour stroma in order to restrain tumour development.

%--------------------------------------------------------------
%			   	     Temporal dependent results
%--------------------------------------------------------------
\subsection{Temporal Dependency of the Re-polarisation}
\begin{figure*}
	\centering
      \includegraphics[width=0.75\textwidth]{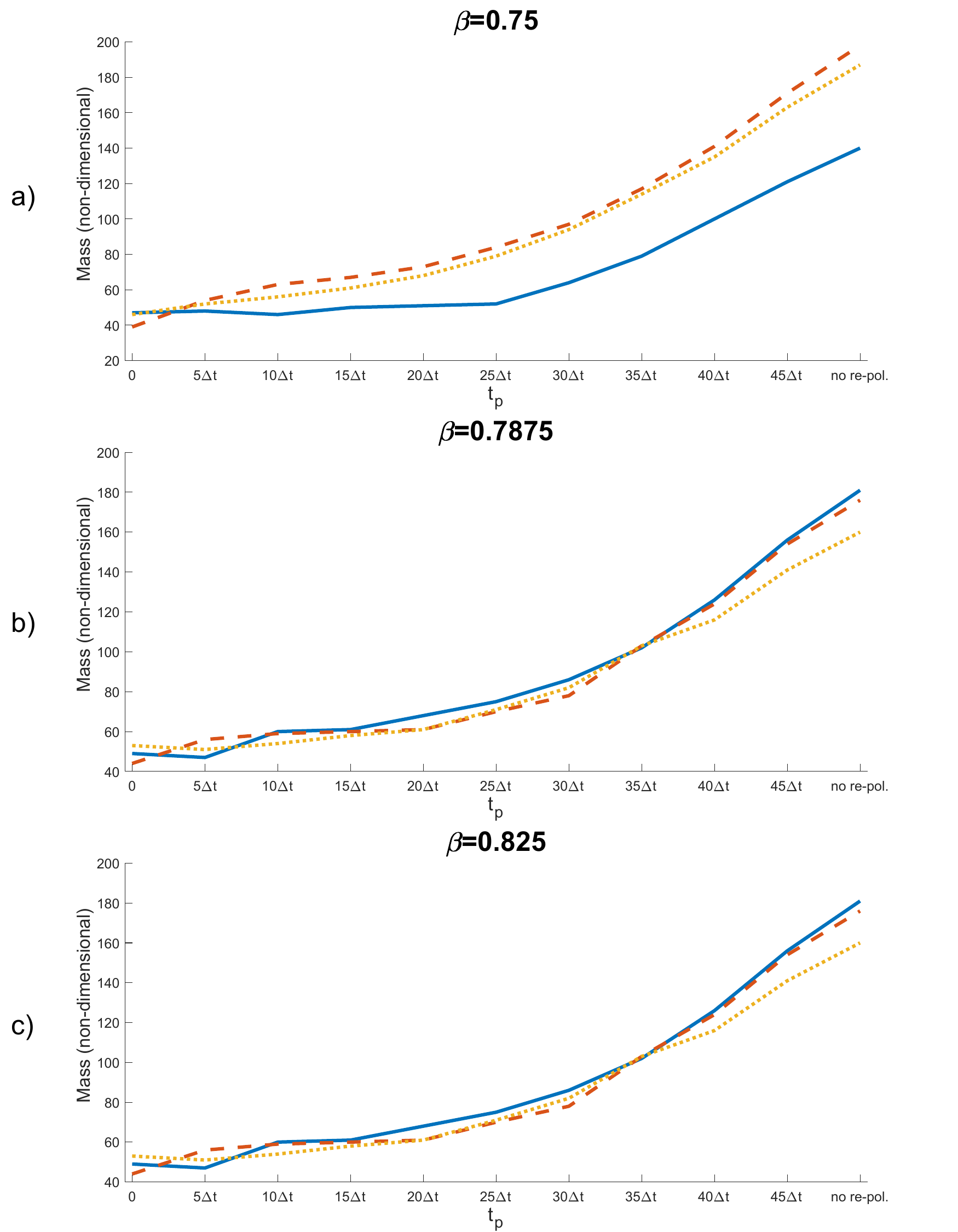}
    \caption{Results of varying the re-polarisation start time $t_{p} \in \{ 0, 5\Delta t, \dots, 45\Delta t, 50\Delta t \}$. Solid-blue lines correspond to a high MMP secretion rate cancer type ($\beta_{lc} = 3.0$, $\beta_{Fc} = 1.5$, $\alpha_{c} = 0.625$), the dashed-orange lines to a type with medium secretion rate ($\beta_{lc} = 2.0$, $\beta_{Fc} = 1.0$, $\alpha_{c} = 0.42$) and the dotted-yellow lines to low secretion rate type cancer ($\beta_{lc} = 1.0$, $\beta_{Fc} = 0.5$, $\alpha_{c} = 0.21$). The environment controller $\beta$ is set to a) $\beta = 0.75$, b) $\beta =0.7875$, c) $\beta =0.825$.}
    \label{fig:Temporal_Repolarisation_All}
\end{figure*}

Let us now concentrate on the temporal aspects of a possible re-polarisation strategy. Since in the previous section we have shown that the M2$\to$M1 re-polarisation impacts predominantly \re{tumour} mass, and since this effect was the strongest when $R_{p} = 0$, here we \re{investigate the effect of changing $t_{p}$ on tumour mass} when we re-polarise the M2 TAMs within the whole tumour domain. Thus, by varying $t_{p}$ used in \eqref{Macrophage_Re_Polarisation} we aim to investigate the effectiveness of such strategy when we introduce the re-polarisation of M2 TAMs at time $t_{p} > 0$. This is crucial since tumours are only detectable above a certain size, and so a potential treatment that uses re-polarisation cannot be started at time $t_{p} = 0$. To this end, in Fig.~\ref{fig:Temporal_Repolarisation_All} we present the change in the dimensionless tumour mass with respect to the re-polarisation start time $t_{p}$.
%------------ Figure for temporal changes
%------------
As before, Fig.~\ref{fig:Temporal_Repolarisation_All} a) corresponds to $\beta = 0.75$, Fig.~\ref{fig:Temporal_Repolarisation_All} b) to $\beta = 0.8785$ and Fig.~\ref{fig:Temporal_Repolarisation_All} c) to $\beta = 0.825$. \re{For each of these sub-cases,} the solid-blue line corresponds to high MMP secretion rate, the dashed-orange line to medium MMP secretion rate, and the dotted-yellow line to low MMP secretion rate. The results in Fig.~\ref{fig:Temporal_Repolarisation_All} \re{suggest that the lowest tumour mass} can be achieved by starting the re-polarisation as soon as possible, as one would expect. However, even at moderate times, for instance at $25 \Delta t$, \re{the tumour mass can be controlled at relatively low values. All these tumour values obtained following re-polarisation at different times $t_{p}$ are smaller than the tumour value obtained for no re-polarisation (see the last point on the horizontal axis in Fig.~\ref{fig:Temporal_Repolarisation_All}).} This indicates that a re-polarisation based strategy is not only viable when it is introduced at a very early stage, but it could also be effective even against a more advanced and larger tumour.
 %Furthermore, these results, presented in Figs.~\ref{fig:Spatial_Repolarisation_All} and \ref{fig:Temporal_Repolarisation_All}, indicate that a re-polarisation based strategy may be more general and could potentially work against different types of cancer.

%--------------------------------------------------------------
%	  					  Conclusions
%--------------------------------------------------------------
\section{Conclusions}\label{section:Conclusions}
\re{In this study} we have further developed and extended a multi-scale moving boundary framework for cancer invasion \cite{Shuttleworth_2019,Suveges_2020,Dumitru_et_al_2013} by considering also the dynamics of the anti-tumoral M1 TAMs as well as the nutrients. On one hand, we took into account the nutrients since every cell requires them to live and function properly and so their presence is indispensable. On the other hand, we focused on the classically activated M1-like macrophages since they are known to be capable of killing cancer cells. Moreover, since macrophages are one of the most abundant immune cells in the tumour micro-environment and their plasticity enables them to switch phenotype, they are prime candidates to assist in the fight against cancer. To this end, we investigated how the re-polarisation of the M2 TAMs into M1 TAMs can affect cancer development, by focusing especially on the macrophage populations near the leading edge of the tumour. Specifically, we studied the spatial aspect of the M2$\to$M1 re-polarisation through the definition of a re-polarisation domain $\Omega_{p}(t,R_{p})$, and the temporal aspect via the starting re-polarisation time $t_{p}$ used in \eqref{Macrophage_Re_Polarisation}.

To propose new hypotheses, we first introduced a macro-scale quasi-steady reaction-diffusion equation for the nutrients where we considered the spatial transport to be described by diffusion with constant-coefficient as well as a linear uptake rate. To account for the effect of nutrients on the different cell functions, we defined four effect-functions that we used \re{for} the rest of tumour dynamics. Furthermore, we introduced another macro-scale equation, describing the behaviour of the M1 TAMs where the motility of the M1 phenotype is represented both by random and directed movements. The rest of the equation involves an influx term, a nutrient-dependent proliferation and death laws as well as a nutrient-dependent polarisation and a re-polarisation terms that describe the switch between the two phenotypes. Similarly to the M2 TAMs, the M1 phenotype also secrete MDEs \re{that can degrade the ECM}. Therefore, this M1 phenotype directly contributes to the re-arrangement of the micro-fibres constituents, as well as serving as a source of MDEs for the proteolytic processes that occur on the invasive edge of the tumour.

We used this extended model to explore the possibilities of \re{macrophage }re-polarisation \re{to depend on spatial domain as well as on time}. First, in Figs.~\ref{fig:No_Repolarisation}, \ref{fig:Repolarisation} and \ref{fig:Spatial_Repolarisation_All} we presented the result of the spatial dependency by varying the re-polarisation domain $\Omega_{p}(t,R_{p})$. We concluded that even though the tumour spread does seems to be affected much \re{by the M2$\to$M1 re-polarisation} (which may be expected as biological studies \cite{Hedbrant2015,Ma2010,Ohri2009} have shown that M1 TAMs located in the stoma can promote cancer progression), the tumour mass can be reduced significantly. Therefore, even though we may need additional strategies directly targeting tumour spread, the tumour mass can potentially be reduced by the re-polarisation of M2 TAMs to M1 TAMs.

 Finally, \re{since tumours are only detectable above a certain size and therefore the M2$\to$M1 re-polarisation is usually applied at later stages in tumour development,} in Fig.~\ref{fig:Temporal_Repolarisation_All} we investigated the temporal dependency \re{of M2$\to$M1 re-polarisation}. There, we showed that 
 \re{while the smallest tumour mass can be obtained when the re-polarisation starts as soon as possible, in some cases it is possible to keep the tumour under control even when we re-polarise at later times. For example, in Fig.~\ref{fig:Temporal_Repolarisation_All}a), for high MMP secretion rates, tumour mass did not change when $t_{p}\leq 25\Delta t$. However, tumour mass slowly increased as we delayed the re-polarisation time for medium and low MMP secretion rates.}
 %for all considered cases we could achieve a considerable reduction in tumour mass even when we started the re-polarisation process at later stages. To this end, we concluded that although re-polarising M2 TAMs into M1 TAMs alone may not be sufficient for reducing tumour spread, it is adequate to reduce the mass of the cancer cells for multiple type of tumours. 

To conclude, we suggest that in addition to the re-polarisation of M2 TAMs, we also need additional strategies targeting tumour spread or tumour stroma in order to fully stop the tumour from advancing.

%--------------------------------------------------------------
%	  					   Funding
%--------------------------------------------------------------
\section*{Acknowledgment}
The authors would like to acknowledge that this research was supported by the EPSRC DTA EP/R513192/1 grant.

%--------------------------------------------------------------
%	  					   Appendix
%--------------------------------------------------------------
\appendix
%--------------------------------------------------------------
%						  Parameters
%--------------------------------------------------------------
\section{Parameter Values}
\label{Parameter_Set}
\begin{center}
\begin{longtable}{lllc}
\caption{Parameter set} \label{tab:long} \\

\hline \multicolumn{1}{l}{\textbf{Variable}} & \multicolumn{1}{l}{\textbf{Value}} & \multicolumn{1}{l}{\textbf{Description}} & \multicolumn{1}{c}{\textbf{Reference}} \\ \hline 
\endfirsthead

\multicolumn{4}{c}%
{{\bfseries \tablename\ \thetable{} -- continued from previous page}} \\
\hline \multicolumn{1}{l}{\textbf{Variable}} & \multicolumn{1}{l}{\textbf{Value}} & \multicolumn{1}{l}{\textbf{Description}} & \multicolumn{1}{c}{\textbf{Reference}} \\ \hline 
\endhead

\hline \multicolumn{4}{r}{{Continued on next page}} \\ \hline
\endfoot

\hline \hline
\endlastfoot

$D_{c}$ & $10^{-4}$ & Diffusion coeff. for the cancer cell \\
& & population $c$ & \cite{Domschke_et_al_2014} \\
$D_{cM_{1}}$ & $4.5$ & Coeff. of the M2 TAMs dependence \\
& & in the cancer diffusion & Estimated \\
$D_{cM_{2}}$ & $1.8$ & Coeff. of the M2 TAMs dependence \\
& & in the cancer diffusion & \cite{Hsu2012} \\
$D_{cF}$ & $8$ & Coeff. of the fibres dependence \\
& & in the cancer diffusion & \cite{Suveges_2020} \\
$D_{M_{1}}$ & $5 \times 10^{-5}$ & Diffusion coeff. for the M1 TAM \\
& & population & \cite{Hayenga_2015} \\
$D_{M_{1}F}$ & $16$ & Coeff. of the ECM stiffness \\
& & dependence in the M1 TAMs diffusion & \cite{Hayenga_2015} \\
$D_{M_{2}}$ & $5 \times 10^{-5}$ & Diffusion coeff. for the M2 TAM \\
& & population $M$ & \cite{Hayenga_2015} \\
$D_{M_{2}F}$ & $16$ & Coeff. of the ECM stiffness \\
& & dependence in the M2 TAMs diffusion & \cite{Hayenga_2015} \\
$D_{\sigma}$ & $1$ & Diffusion coeff. for the nutrients & Estimated \\
$D_{m}$ & $2.5 \times 10^{-3}$ & Diffusion coeff. for MDEs & \cite{Peng2016} \\
$\S_{max}$ & $0.5$ & Cell-cell adhesion coeff. & \cite{Shuttleworth_2019} \\
$\S_{min}$ & $0.01$ & Minimum level of cell-cell adhesion & \cite{Suveges_2020} \\
$\S_{cl}$ & $0.01$ & Cell-non-fibre adhesion coeff. & \cite{Shuttleworth_2019} \\
$\S_{cF}$ & $0.3$ & Cell-fibre adhesion coeff. & \cite{Domschke_et_al_2014} \\
$\S_{cM}$ & $0.125$ & Cell-macrophage adhesion \\
& & coeff. & \cite{Suveges_2020} \\
$\S_{M_{1}M}$ & $0.175$ & M1 TAMs self adhesion coeff. & \cite{Cui_2018} \\
$\S_{M_{2}M}$ & $0.05$ & M2 TAMs self adhesion coeff. & \cite{Cui_2018} \\
$\S_{Mc}$ & $0.125$ & Macrophage-cancer adhesion \\
& & coeff. & \cite{Suveges_2020} \\
$\S_{M\sigma}$ & $0.1$ & Strength of the Macrophage-nutrients \\
& & relationship & Estimated \\
$\mu_{c}$ & $0.25$ & Proliferation coeff. for cancer cell \\
& & population $c$ & \cite{Domschke_et_al_2014} \\
$\mu_{cM_{1}}$ & $4.0$ & Coeff. for the M1 TAMs dependence \\
& & in the cancer cell proliferation & Estimated \\
$\mu_{cM_{2}}$ & $1.4$ & Coeff. for the M2 TAMs dependence \\
& & in the cancer cell proliferation & \cite{Hu2015} \\
$M_{0}$ & $0.05$ & Influx of the M1 TAMs & Estimated \\
$\mu_{M}$ & $0.2$ & Inside tissue proliferation rate of M1 \\
& & and M2 TAMs & Estimated \\
$\mu_{MF}$ & $1.8$ & Coeff. for the ECM stiffness dependence \\
& & in the M1 and M2 TAMs prolifs. & \cite{Hayenga_2015} \\
$d_{c}$ & $0.1$ & Decay rate of the cancer cells & Estimated \\
$d_{cM_{1}}$ & $50.0$ & M1 TAMs killing rate of the cancer cells & Estimated \\
$d_{M}$ & $0.03$ & Decay rate of M1 and M2 TAMs & \cite{Strachan_2013} \\
$d_{\sigma}$ & $80.0$ & Nutrients uptake rate & Estimated \\
$\beta_{lc}$ & $3.0$ & Degradation coeff. for non-fibre ECM \\
& & due to the tumour & Estimated \\
$\beta_{lM_{1}}$ & $3.84$ & Degradation coeff. for non-fibre ECM \\
& & due to the M1 TAMs & Estimated \\
$\beta_{lM_{2}}$ & $0.96$ & Degradation coeff. for non-fibre ECM \\
& & due to the M2 TAMs & Estimated \\
$\beta_{Fc}$ & $1.5$ & Degradation coeff. for fibre ECM \\
& & due to the tumour & Estimated \\
$\beta_{FM_{1}}$ & $1.92$ & Degradation coeff. for fibre ECM \\
& & due to the M1 TAMs & Estimated \\
$\beta_{FM_{2}}$ & $0.48$ & Degradation coeff. for fibre ECM \\
& & due to the M2 TAMs & Estimated \\
$\gamma_{0}$ & $0$ & Baseline remodelling & \cite{Suveges_2020} \\
$\gamma_{M_{2}}$ & $0$ & Remodelling of the non-fibre ECM \\
& & due to the M2 TAMs & \cite{Suveges_2020} \\
$p_{12}$ & $14.0$ & Baseline polarisation rate & Estimated \\
$p_{21}$ & $6.0$ & Baseline re-polarisation rate & Estimated \\
$t_{p}$ & $0$ & Re-polarisation staring time & Estimated \\
$R_{p}$ & $0$ & Radius of the re-polarisation domain & Estimated \\
$\sigma_{nor}$ & $0.4$ & Normal level of nutrients & Estimated \\
$\sigma_{p}$ & $0.35$ & Maximal level of nutrients that \\
& & is still sufficient & Estimated \\
$\sigma_{n}$ & $0.2$ & Necrotic threshold for the \\
& & nutrients level & Estimated \\
$\Psi_{p,max}$ & $1.0$ & Maximal proliferation enhancement \\
& & rate & Estimated \\
$\Psi_{d,max}$ & $5.0$ & Maximal death enhancement rate & Estimated \\
$\Psi_{d,min}$ & $1.0$ & Minimal death enhancement rate & Estimated \\
$\Psi_{M,max}$ & $2.0$ & Maximal polarisation enhancement \\
& & rate & Estimated \\
$\Psi_{M,min}$ & $1.0$ & Minimal polarisation enhancement \\
& & rate & Estimated \\
$\alpha_{c}$ & $0.625$ & MDE secreting rate by the cancer cells & Estimated \\
$\alpha_{M_{1}}$ & $0.8$ & MDE secreting rate by the M1 TAMs & Estimated \\
$\alpha_{M_{2}}$ & $0.2$ & MDE secreting rate by the M2 TAMs & Estimated \\
$\beta$ & $0.75$ & Optimal level/tissue environment \\
& & controller of the ECM density & \cite{Dumitru_et_al_2013} \\
$R$ & $0.15$ & Sensing radius & \cite{Shuttleworth_2019} \\
$r$ & $0.0016$ & Width of micro-fibres & \cite{Shuttleworth_2019} \\
$f_{max}$ & $0.636$ & Maximum of micro-fibre density \\
& & at any point & \cite{Shuttleworth_2019} \\
$h_{L}$ & $0.03125$ & Macro-scale spatial step-size & \cite{Dumitru_et_al_2013} \\
$\epsilon$ & $0.0625$ & Size of the boundary micro-domain $\epsilon Y (\textbf{x})$ & \cite{Dumitru_et_al_2013} \\
$\delta$ & $0.03125$ & Size of the fibre micro-domain $\delta Y (\textbf{x})$ & \cite{Shuttleworth_2019} \\
\end{longtable}
\end{center}

%--------------------------------------------------------------
%						Outer Boundary
%--------------------------------------------------------------
\section{Definition of the Outer Boundary: $\partial\Omega_{o}(t)$}
\label{outerBoundarySet}
\begin{figure}
\centering
  \includegraphics[width=0.3\textwidth]{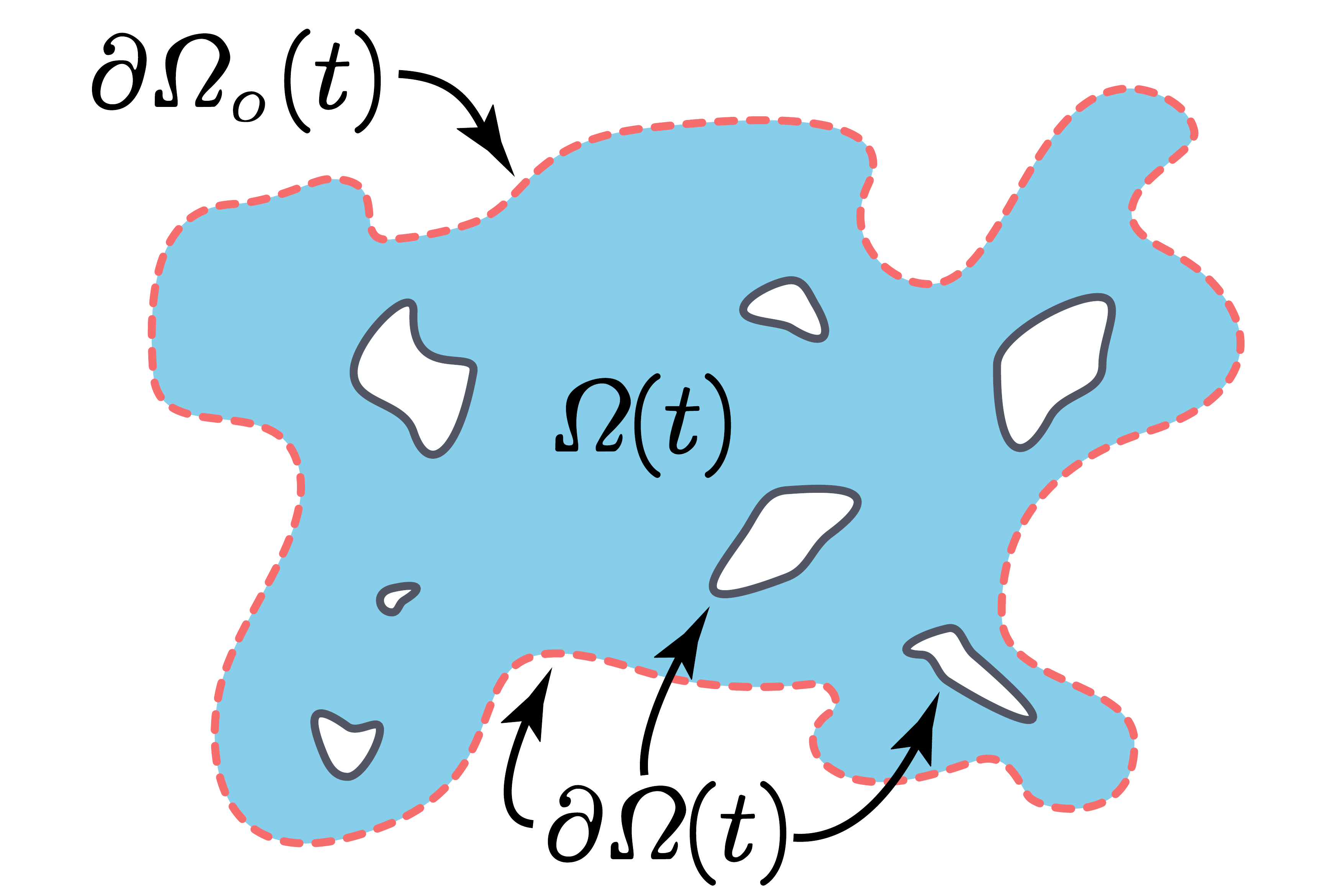}
\caption{Schematic of the outer boundary $\partial \Omega_{o}(t)$ that is highlighted with the \emph{dashed red} line.}
\label{fig:Outer_Boundary}
\end{figure}

Let $x\in \partial\Omega(t)$. Then, $x\in \partial\Omega_{o}(t)$ if and only if there exists $\phi_{_{x}}:[0,\infty)\to \R^{d}$ such that the following four properties hold true simultaneously:
\begin{align}
1)\; & \phi_{_{x}}(0)=x, \nonumber\\
2)\; & \phi_{_{x}}(s)\neq x, \quad, \forall s\in (0, \infty), \nonumber\\
3)\; & Im \phi_{_{x}}\setminus\{x\}\subset \complement  \Omega(t), \nonumber\\
4)\; & \lim\limits_{s \to\infty}dist (\phi(s), \partial \Omega(t))=\infty, \nonumber
\end{align}
where $\forall \, s \in (0,\infty)$, we have $dist(\phi(s), \partial \Omega(t)):=\inf\limits_{z\in\partial \Omega(t)}\nor{\phi(s)-z}_{_{2}}$ and represents the Euclidean distance from $\phi(s)$ to $\partial \Omega(t)$. For a visual representation, we refer the reader to Fig.~\ref{fig:Outer_Boundary}.
%------------ Figure for outre tumour boundary
%------------

%--------------------------------------------------------------
%					 Re-polarisation Domain
%--------------------------------------------------------------
\section{Definition of the Re-polarisation Domain: $\Omega_{p}(t, R_{p})$}
\label{Appendix_Re_Polarisation_Domain}
Here, we define the re-polarisation domain used in \eqref{Macrophage_Re_Polarisation}, which allows us to investigate the spatial dependency of a re-polarisation strategy. Hence, let us defined this new domain by
\begin{equation*}
	\Omega_{p}(t, R_{p}) := \{ x \in \Omega(t), x_{B} \in \partial \Omega_{o}(t) | \nor{x-x_{B}}_{2} \geq R_{p} \} \subset \Omega(t),
\end{equation*}
where $R_{p} > 0$ is a constant and $\partial \Omega_{o}(t)$ is the outer tumour boundary, defined in Appendix \ref{outerBoundarySet}. Again, this new domain is illustrated in Fig.~\ref{fig:Re_Polarisation_Domain}.

%--------------------------------------------------------------
%					Standard Mollifier
%--------------------------------------------------------------
\section{Standard Mollifier}
\label{Standard_Mollifier}
The smooth, compact support function $\psi : \R^{N} \rightarrow \R$ that is used throughout the paper which is indeed the standard mollifier, is defined as usual by
\begin{equation*}
	\psi(x) := 
	\begin{cases}
	\dfrac{\exp \Bigg( \dfrac{-1}{1 - \nor{z}_{2}^{2}}\Bigg)}{\int\limits_{\Bila(0,1)} \exp \Bigg( \dfrac{-1}{1 - \nor{z}_{2}^{2}}\Bigg) dz} & \text{if } x \in \Bila(0,1), \\[20pt]
	\hfil 0 & \text{if } x \notin \Bila(0,1).
	\end{cases}
\end{equation*}

%--------------------------------------------------------------
%					Rearrangement process
%--------------------------------------------------------------
\section{Further Details on the Micro-Fibre Rearrangement Process}
\label{Fibre_Reallocation_Process}
As we mentioned in Section~\ref{Fibre_Micro_Scale}, the micro-fibres rearrangement process is instigated by the emerging macro-scale fluxes that acts on the mass distribution of the micro-fibres $f(z,t)$ on $\delta Y(x)$. As a consequence, the micro-fibres are redistributed within $\delta Y(x)$ and its neighbouring micro-domains. As detailed in \cite{Shuttleworth_2019}, we calculate the new position $z^{*}$ of any micro-point $z$ via the reallocation vector $\nu_{\delta Y(x)}(z,t)$ which accounts both for the rearrangement vector $r(\delta Y(x),t)$ \eqref{Rearrangement_Vector} and upon the degree of alignment between $r(\delta Y(x),t)$ and the \emph{barycentral position vector} $x_{dir} := z-x$ as well as it takes into consideration the level of fibres at position $z$. Hence, the reallocation vector is defined by
\begin{equation*}
	\nu_{\delta Y(x)}(z,t) := \big[ x_{dir}(x) + r(\delta Y(x),t) \big] \cdot \dfrac{f(z,t) [f_{max} - f(z,t)]}{f^{*} + \nor{ r(\delta Y(x),t) - x_{dir}(x)}_{2}} \cdot \chi_{ \{f(z,t) > 0\} },
\end{equation*}
where $f_{max} > 0$ is the maximum possible level of fibres at any micro-point, $f^{*} := f(x,t) / f_{max}$ is the saturation level of micro-fibres and $\chi_{ \{f(z,t) > 0\} }$ represents the characteristic function of the micro-fibres support. Finally, the amount of fibres that is reallocated from $z$ to the target position $z^{*}$ is determined by a movement probability $p_{move}$ that monitors the available free space at the new position $z^{*}$ and it is given by
\begin{equation*}
	p_{move} := \max \Bigg( 0, \dfrac{f_{max} - f(z^{*},t)}{f_{max}} \Bigg).
\end{equation*}
Hence, the amount of fibres that moves to $z^{*}$ is given by $p_{move} \cdot f(z,t)$ and the amount of fibres that remains at position $z$ is consequently given by $(1 - p_{move}) \cdot f(z,t)$.

%--------------------------------------------------------------
%					Boundary Movement
%--------------------------------------------------------------
\section{Further details on the boundary MDE micro-dynamics}\label{Appendix_MDE_Micro_Process}
\dt{As derived in \cite{Dumitru_et_al_2013}, the MDE dynamics \eqref{MDE_Equation} determine the way the macro-scale tumour boundary moves, and in the following we will briefly detail the way the associated direction $\eta_{\epsilon Y}$ of movement  and displacement magnitude $\xi_{\epsilon Y}$ of the tumour boundary within the peritumoural region captured by any given boundary micro-domain $\epsilon Y\in \P(t_{0})$.   For this, we consider an appropriate dyadic decomposition $\{\D_{k}\}_{_{\I}}$ of the micro-domain $\epsilon Y$, and} denote by $y_{k}$ the barycentre of each $\D_{k}$. \dt{Then, we subsample this family of dyadic cubes and we select in a subfamily $\{\D_{k}\}_{_{\J^{*}}}$ only those dyadic cubes that are furthest away from the centre $x$ as well as are outside of the tumour domain $\Omega(t_{0})$ and carry an above the average mass of MDEs.}
Then, the direction and magnitude of the boundary movement on each $\epsilon Y\in \P(t_{0})$ are given by
\begin{equation}
	\begin{split}
		\eta_{\epsilon Y (x^{*}_{\epsilon Y})} & := x^{*}_{\epsilon Y} + \nu \sum_{l \in \mathcal{J}^{*}} \Bigg ( \int\limits_{\mathcal{D}_{l}} m(y, \tau) \; dy \Bigg ) \Big ( y_{l} - x^{*}_{\epsilon Y} \Big ), \qquad \nu \in  [0, \infty), \\
		\xi_{\epsilon Y (x^{*}_{\epsilon Y})} & := \sum_{l \in \mathcal{J}^{*}} \dfrac{\int\limits_{\mathcal{D}_{l}} m(y, \tau) \; dy}{\sum\limits_{l \in \mathcal{J}^{*}} \int\limits_{\mathcal{D}_{l}} m(y, \tau) \; dy} \Big | \overrightarrow{x^{*}_{\epsilon Y} y_{l}} \Big |,
	\end{split}
	\label{Direction_And_Magnitude}
\end{equation}
respectively. Ultimately, using \eqref{Direction_And_Magnitude}, we determine the new positions of each boundary point which leads to a new expanded tumour domain $\Omega(t_{0} + \Delta t)$. However, the movement of the boundary towards these newly computed locations is exercised provided that enough but not complete degradation of the peritunoural ECM took place, aspect that is quantified through the transitional probability of the boundary movement by
\begin{equation*}
	q(x^{*}_{\epsilon Y}) := \dfrac{\int\limits_{\epsilon Y(x^{*}_{\epsilon Y}) \setminus \Omega(t_{0})} m(y, \tau) dy}{\int\limits_{\epsilon Y(x^{*}_{\epsilon Y})} m(y, \tau) dy},.
\end{equation*}
The movement is exercised or not provided that $q(x^{*}_{\epsilon Y})$ exceed a certain tissue threshold $\omega(\beta, \epsilon Y)$ \cite{Dumitru_et_al_2013}. This tissue threshold allows the boundary to move to the new position only in favourable environment and rejects that not only when there is too much ECM remaining, but also when the level of ECM is not sufficient to support the tumour movement. Thus, denoting this threshold by $\omega (0,1) \times \epsilon Y \rightarrow [0,1]$, we define this by
\begin{equation*}
	\omega(\beta,\! \epsilon Y)\! :=\!
	\begin{cases}
		\sin \! \bigg[ \dfrac{\pi}{2} \bigg( 1 \!-\! \dfrac{v(x^{*}_{\epsilon Y}, t_{0} + \Delta t)}{\beta \cdot \sup\limits_{\xi \in \partial \Omega(t_{0})} \!\! v(\xi, t_{0} \!+\! \Delta t)} \bigg) \bigg] & \!\! \text{if } \dfrac{v(x^{*}_{\epsilon Y}, t_{0} \!+\! \Delta t)}{\sup\limits_{\xi \in \partial \Omega(t_{0})} \!\! v(\xi, t_{0} \!+\! \Delta t)} \leq \beta, \\[30pt]
		\sin \! \bigg[ \dfrac{\pi}{2(1 \!-\! \beta)} \bigg( \dfrac{v(x^{*}_{\epsilon Y}, t_{0} + \Delta t)}{\sup\limits_{\xi \in \partial \Omega(t_{0})} \!\! v(\xi, t_{0} \!+\! \Delta t)} - \beta \bigg) \bigg] & \!\! \text{if } \dfrac{v(x^{*}_{\epsilon Y}, t_{0} \!+\! \Delta t)}{\sup\limits_{\xi \in \partial \Omega(t_{0})} \!\! v(\xi, t_{0} \!+\! \Delta t)} > \beta,
	\end{cases}
\end{equation*}
where $\beta \in (0,1)$ controls the optimal level of ECM for cancer invasion, and $v(x,t) := l(x,t) + F(x,t)$ is the cumulative ECM. 

%--------------------------------------------------------------
%					   Adhesion Numerics
%--------------------------------------------------------------
\section{Approximation of the Adhesion Integrals}\label{Appendix_Adhesion}
\begin{figure*}[h!]
\centering
  \includegraphics[width=0.65\textwidth]{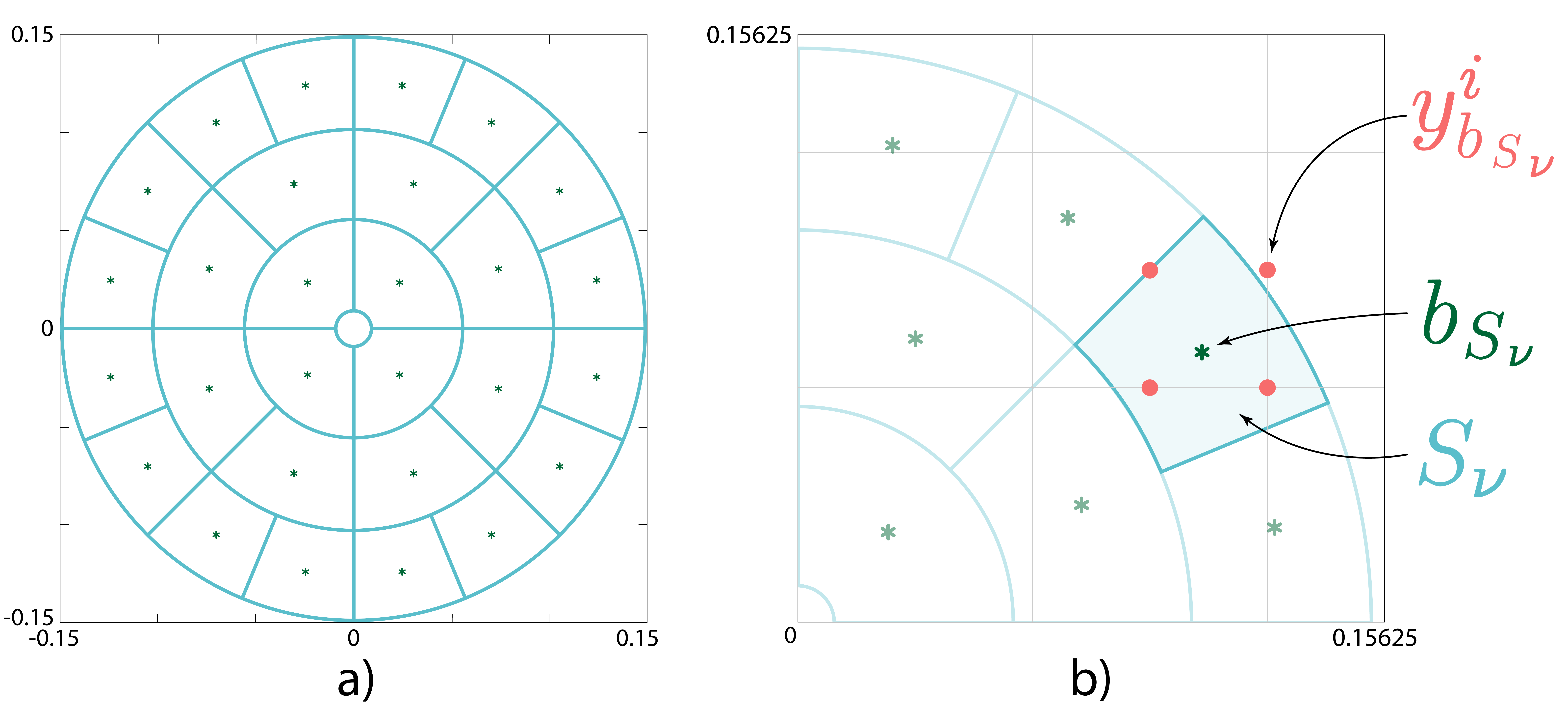}
\caption{Illustration of the sensing region $\Bila (0, R)$. a) Decomposition of the sensing region into annulus sectors $S_{\nu}$ with barycentres $b_{S_{nu}}$ (green stars). b) Illustration of one of the barycentre $b_{S_{\nu}}$ with the 4 on-grid neighbours $\{y^{k}_{b_{S_{\nu}}}\}_{k=1,\dots,4}$.}
\label{fig:Sensing_Region}
\end{figure*}

Each cell adhesion, cell-fibres ECM and cell-non-fibres ECM adhesions processes (modelled in \eqref{General_TAMs_Adhesion} and \eqref{Cancer_Adhasion}) take place within the \emph{sensing region} $x + \Bila(0,R)$. To partition this region, we adopt the technique developed in \cite{Shuttleworth_2019} and split the sensing region into $N_{s}$ annulus sectors, illustrated in Fig.~\ref{fig:Sensing_Region}a).
%------------ Figure for the sensing region
%------------
Hence, first the region $\Bila(0,R)$ is split into $s$ annuli and then each of these annuli are further split into $2^{m+(k-1)}$ uniformly distributed sectors (which is dyadically increases with $k = 1,\dots,s$, \emph{i.e.,} increases as we progress from the inner most annulus to the biggest last one) and consequently, the number of annulus sectors is given by
\begin{equation*}
	N_{s} := \sum_{k=1}^{s} 2^{m+(k-1)}.
\end{equation*}
Then, denoting these annulus sectors by $S_{\nu}$, with $\nu = 1,\dots,N_{s}$, we are able to approximate these adhesion integrals \eqref{General_TAMs_Adhesion} and \eqref{Cancer_Adhasion} as sum of the integrals of the step functions associated with each $S_{\nu}$. The values of each step functions are appropriately approximated as a linear combination of the mean values of cancer cell, M1 TAMs, M2 TAM population as well as both ECM components.

To approximate these mean values we need to integrate the corresponding densities and the fibre orientations on each $S_{\nu}$. For this, we follow \cite{Suveges_2020}, and denote the \emph{off-grid} barycentres of each annulus by $b_{S_{\nu}}$ and observe that each of them belongs to a rectangle defied by its four immediate on-grid macro-mesh neighbouring nodes $\{ y^{i}_{b_{S_{\nu}}} \}_{i = 1,\cdot,4}$. Hence, we have
\begin{equation*}
	b_{S_{\nu}} \in [y^{1}_{b_{S_{\nu}}}, y^{2}_{b_{S_{\nu}}}] \times [y^{3}_{b_{S_{\nu}}}, y^{4}_{b_{S_{\nu}}}],
\end{equation*}
which is illustrated in Fig.~\ref{fig:Sensing_Region}b). This enables us to approximate the values of cancer cell, M1 TAMs and M2 TAMs, both ECM components as well as the ECM fibres orientations on each $S_{\nu}$ at $b_{S_{\nu}}$ by bi-linear shape functions. Hence, each value on $b_{S_{\nu}}$ is approximated as a convex combination of the values at the neighbouring on-grid nodes $\{ y^{i}_{b_{S_{\nu}}} \}_{i = 1,\dots,4}$ with uniquely determined weights $\beta^{k}_{\nu}$, $k = 1,\dots,4.$

Then we use the convolution-driven approach \cite{Suveges_2020} and so using the weights $\beta^{k}_{\nu}$, we construct $N_{s}$ different matrices (one for each barycentre $S_{\nu}$). Considering an arbitrary annulus sector $S_{\nu}$ the associated matrix $\K^{S_{\nu}}_{\A}$ is given by
\begin{equation*}
	\K_{\A_{M}}^{S_{\nu}} =
	\begin{bmatrix} 
    	0 & \cdots & \cdots & \cdots & \cdots & 0 \\
    	\vdots & &\ddots & & \\
    	0 & \cdots & \beta_{\nu}^{1} & \beta_{\nu}^{3} & \cdots & 0 \\
    	0 & \cdots & \beta_{\nu}^{2} & \beta_{\nu}^{4} & \cdots & 0 \\
    	\vdots & & & \ddots & \\
    	0 & \cdots & \cdots & \cdots & \cdots & 0 \\
    \end{bmatrix},
\end{equation*}
which is a $P \times P$ matrix, with $P := \ceil{2 R / \Delta x} + 1$. The size of these matrices $\K_{\A_{M}}^{S_{\nu}}$ corresponds to the minimal squared region of macro-nodes $\{ y^{i}_{x} \}_{i=1\dots,P^{2}}$ that covers the entire sensing region $\Bila (0, R)$. Also, we observe that the locations of the four non-zero values in each $\K_{\A_{M}}^{S_{\nu}}$ are in harmony with the locations of the four neighbouring on-grid nodes associated with the barycentre $\S_{\nu}$ and that these values are given by the associated weights $\beta^{k}_{\nu}$.

Therefore, at any time $t_{n}$, the approximation of the adhesion integral terms $\A \in \{ \A_{c}(x,t,\bu,\theta_{f}), \A_{M}(x,t,\bu,\S_{MM}) \}$ (defined in \eqref{General_TAMs_Adhesion} and \eqref{Cancer_Adhasion}) are given by
\begin{equation}
	\A^{n} = \sum_{\nu=1}^{N_{s}} \A_{S_{\nu}}^{n}.
	\label{Attration_Adheison_Split}
\end{equation}
Finally, in \eqref{Attration_Adheison_Split}, $\A_{S_{\nu}}^{n}$ is defined for each adhesion integral as:
\begin{itemize}
	\item for $\A := \A_{c}(x,t,\bu,\theta_{f})$
	\begin{equation*}
		\begin{split}
			(\A_{c})_{S_{\nu}}^{n} \!=\! & \frac{K(b_{S_{\nu}})}{R} \! \bigg [ 
			\Big ( \widehat{n}_{\nu} \widetilde{\K}_{\A}^{S_{\nu}} \Big ) \! \ast \! \Big ( \big( \S_{cF} F^{n} \big) \circ \big( 1 \!-\! \rho(\bu^{n}) \big)^{+} \Big ) \\
			& + \Big ( n_{\nu} \widetilde{\K}_{\A}^{S_{\nu}} \Big ) \! \ast \! \Big ( \big( \S_{cc} c^{n} \!+\! \S_{cl} l^{n} \!+\! \S_{cM} (M^{n}_{1} \!+\! M^{n}_{2}) \big) \! \circ \! \big( 1 \!-\! \rho(\bu^{n}) \big)^{+} \Big ) \bigg ],
		\end{split}
	\end{equation*}
	\item for $\A := \A_{M}(x,t,\bu,\S_{MM})$
	\begin{equation*}
		(\A_{M})_{S_{\nu}}^{n} \!=\! \frac{K(b_{S_{\nu}})}{R} \! \bigg [ 
		\Big ( n_{\nu} \widetilde{\K}_{\A}^{S_{\nu}} \Big ) \ast \Big ( \big( \S_{Mc} c^{n} + \S_{MM} (M^{n}_{1} + M^{n}_{2}) \big) \circ \big( 1-\rho(\bu^{n}) \big)^{+} \Big ) \bigg ],
	\end{equation*}
\end{itemize}
where $\circ$ is the Hadamard product and $\widetilde{\K}_{\A}^{S_{\nu}}$ are the appropriately constructed matrices for the convolution operator $\ast$ using $\K_{\A}^{S_{\nu}}$, \emph{i.e.,} they are defined by
\begin{equation*}
	\widetilde{\K}_{\A}^{S_{\nu}} := J_{P} \cdot \K_{\A}^{S_{\nu}} \cdot J_{P},
\end{equation*}
where $J_{P}$ is a $P \times P$ exchange matrix given by
\begin{equation*}
	J_{P} =
	\begin{bmatrix} 
    	0 & 0 & \cdots & 0 & 1 \\
    	0 & 0 & & 1 & 0 \\
    	\vdots & & \ddots & & \vdots \\
    	0 & 1 & & 0 & 0 \\
    	1 & 0 & \cdots & 0 & 0 \\
    \end{bmatrix}.
\end{equation*}

%--------------------------------------------------------------
%				Boundary Micro-Scale Numerics
%--------------------------------------------------------------
\section{Numerical Approach of the Boundary Micro-Scale}\label{Appendix_MDE_Calculation}
As explained in Section \ref{MDE_Micro_Scale} the crucial source term $h(\cdot, \cdot)$ for the MDE micro-dynamic \eqref{MDE_Equation} arises through the contribution of the cancer cell, M1 and M2 TAMs populations via \eqref{MDE_Source}. To evaluate this source integral, we seek to use convolution \cite{Suveges_2020} with an appropriately chosen \emph{integral matrix} $\K_{I}$ (the form of which is detailed below). Similarly to the differential operators in Section~\ref{section:Numerical_Approach}, here we also split the source term into two parts, \emph{i.e.,}
\begin{equation}
	h(i,j) := 
	\begin{cases}
		h_{I}(i,j) & \text{if } \I(i,j) = 1, \\
		\hfil 0 & \text{otherwise},
	\end{cases}
	\label{Numerics_Source_Split}
\end{equation}
where $\I$ is defined in \eqref{Indicator_All}. Therefore, the source inside the tumour at any node $(x_{i},y_{j})$ and any time-step $n$ is given by
\begin{equation*}
	h_{I}^{n} = 4 \gamma_{h}^{2} \dfrac{(N_{f} - 1)^{2}}{\I^{n} \ast \mathbbm{1}} \circ \Big [ \big( \alpha_{c} c^{n} + \alpha_{M_{1}} M_{1}^{n} + \alpha_{M_{2}} M_{2}^{n} \big) \ast \K_{I} \Big ],
\end{equation*}
where $\circ$ denotes the Hadamard product, $\ast$ denotes the discrete convolution, $\gamma_{h}$ is defined in \eqref{MDE_Source}, $\alpha_{c}$, $\alpha_{M_{1}}$ and $\alpha_{M_{2}}$ are the MDE secreting rates of the cancer cell, M1 and M2 TAMs respectively. Furthermore, $\I^{n}$ is the discretised indicator function defined in \eqref{Indicator_All}, $N_{f} := 2 \gamma_{h} / \Delta x + 1$ is the size of the integral matrix $\K_{I}$ and $\mathbbm{1}$ is an $N_{f} \times N_{f}$ matrix of ones. In this work, we use the trapezoidal rule to construct the integral matrix $\K_{I}$  \cite{Suveges_2020}, and so we define this as
\begin{equation*}
	\K_{I} = 
	\begin{bmatrix} 
    	1 & 2 & \dots & 2 & 1 \\
    	2 & 4 & \dots & 4 & 2 \\
    	\vdots & &\ddots & & \\
    	2 & 4 & \dots & 4 & 2 \\
    	1 & 2 & \dots & 2 & 1 \\
    \end{bmatrix}.
\end{equation*}
Hence, using \eqref{Numerics_Source_Split} we can calculate the emerging source term however, it is still defined on the macro-scale. Since we use this source term in the MDE micro-scale dynamics \eqref{MDE_Equation}, we need to approximate its values at any micro-node $y$ appropriately. For this, we use bi-linear shape functions on a square micro-mesh \cite{Dumitru_et_al_2013}, which enables us to solve the MDE micro-scale dynamics \eqref{MDE_Equation} by the method of lines. Here, due to the simplicity of the micro-dynamics, we discretise the space using central differences, and for time-marching, we use the backward Euler method.

%--------------------------------------------------------------
%			ENO Approximation Convolution Vectors
%--------------------------------------------------------------
\section{Vectors Induced by the ZSWENO Scheme}\label{WENO_Convolution_Vectors}
Here, for completeness we define the vectors used for the discrete convolutions in \eqref{ENO_Fluxes_Convolution} and \eqref{ZSWENO_Smoothness_Indicators_Convolutions}. Let us start with the vectors that were induced by the ENO fluxes \eqref{ENO_Fluxes_Plus} and \eqref{ENO_Fluxes_Minus}. Hence, first the vectors $\widetilde{\K}^{\pm}_{k,+}$ with $k = 0,\dots,2$ are given by
\begin{equation*}
	\widetilde{\K}^{+}_{0,+} \!\!=\!\!\!
	\begin{bmatrix}
    	0 \\
    	0 \\
    	0 \\
    	\frac{11}{6} \\[0.5em]
    	-\frac{7}{6} \\[0.5em]
    	\frac{1}{3} \\[0.5em]
    	0 \\
    \end{bmatrix}\!\!,
    \;
    \widetilde{\K}^{+}_{1,+} \!\!=\!\!\!
	\begin{bmatrix} 
    	0 \\
    	0 \\
    	\frac{1}{3} \\[0.5em]
    	\frac{5}{6} \\[0.5em]
    	-\frac{1}{6} \\[0.5em]
    	0 \\
    	0 \\
    \end{bmatrix}\!\!,
    \;
    \widetilde{\K}^{+}_{2,+} \!\!=\!\!\!
	\begin{bmatrix} 
    	0 \\
    	-\frac{1}{6} \\[0.5em]
    	\frac{5}{6} \\[0.5em]
    	\frac{1}{3} \\[0.5em]
    	0 \\
    	0 \\
    	0 \\
    \end{bmatrix}\!\!,
    \;
    \widetilde{\K}^{-}_{0,+} \!\!=\!\!\!
	\begin{bmatrix} 
    	\frac{1}{3} \\[0.5em]
    	-\frac{7}{6} \\[0.5em]
    	\frac{11}{6} \\[0.5em]
    	0 \\
    	0 \\
    	0 \\
    	0 \\
    \end{bmatrix}\!\!,
    \;
    \widetilde{\K}^{-}_{1,+} \!\!=\!\!\!
	\begin{bmatrix} 
    	0 \\
    	-\frac{1}{6} \\[0.5em]
    	\frac{5}{6} \\[0.5em]
    	\frac{1}{3} \\[0.5em]
    	0 \\
    	0 \\
    	0 \\
    \end{bmatrix}\!\!,
    \;
    \widetilde{\K}^{-}_{2,+} \!\!=\!\!\!
	\begin{bmatrix} 
    	0 \\
    	0 \\
    	\frac{1}{3} \\[0.5em]
    	\frac{5}{6} \\[0.5em]
    	-\frac{1}{6} \\[0.5em]
    	0 \\
    	0 \\
    \end{bmatrix}\!\!.
\end{equation*}
Similarly, the vectors $\widetilde{\K}^{\pm}_{k,-}$ with $k = 0,\dots,2$ are given by
\begin{equation*}
	\widetilde{\K}^{+}_{0,-} \!\!=\!\!\!
	\begin{bmatrix} 
		0 \\
    	0 \\
    	0 \\
    	0 \\
    	\frac{11}{6} \\[0.5em]
    	-\frac{7}{6} \\[0.5em]
    	\frac{1}{3} \\[0.5em]
    \end{bmatrix}\!\!,
    \;
    \widetilde{\K}^{+}_{1,-} \!\!=\!\!\!
	\begin{bmatrix} 
    	0 \\
    	0 \\
    	0 \\
    	\frac{1}{3} \\[0.5em]
    	\frac{5}{6} \\[0.5em]
    	-\frac{1}{6} \\[0.5em]
    	0 \\
    \end{bmatrix}\!\!,
    \;
    \widetilde{\K}^{+}_{2,-} \!\!=\!\!\!
	\begin{bmatrix} 
    	0 \\
    	0 \\
    	-\frac{1}{6} \\[0.5em]
    	\frac{5}{6} \\[0.5em]
    	\frac{1}{3} \\[0.5em]
    	0 \\
    	0 \\
    \end{bmatrix}\!\!,
    \;
    \widetilde{\K}^{-}_{0,-} \!\!=\!\!\!
	\begin{bmatrix} 
    	0 \\
    	\frac{1}{3} \\[0.5em]
    	-\frac{7}{6} \\[0.5em]
    	\frac{11}{6} \\[0.5em]
    	0 \\
    	0 \\
    	0 \\
    \end{bmatrix}\!\!,
    \;
    \widetilde{\K}^{-}_{1,-} \!\!=\!\!\!
	\begin{bmatrix} 
    	0 \\
    	0 \\
    	-\frac{1}{6} \\[0.5em]
    	\frac{5}{6} \\[0.5em]
    	\frac{1}{3} \\[0.5em]
    	0 \\
    	0 \\
    \end{bmatrix}\!\!,
    \;
    \widetilde{\K}^{-}_{2,-} \!\!=\!\!\!
	\begin{bmatrix} 
    	0 \\
    	0 \\
    	0 \\
    	\frac{1}{3} \\[0.5em]
    	\frac{5}{6} \\[0.5em]
    	-\frac{1}{6} \\[0.5em]
    	0 \\
    \end{bmatrix}\!\!.
\end{equation*}

Then, the other vectors induced by the smoothness indicators \eqref{ZSWENO_Smoothness_Indicators_Plus} and \eqref{ZSWENO_Smoothness_Indicators_Minus} are $\overline{\K}^{\pm}_{k,+}$ and $\overline{\K}^{\pm}_{k,-}$ which were used in the convolutional form \eqref{ZSWENO_Smoothness_Indicators_Convolutions}. First, $\overline{\K}^{\pm}_{k,+}$ are defined by
\begin{equation*}
	\overline{\K}^{+}_{0,+} \!\!=\!\!\!
	\begin{bmatrix}
    	0 \\
    	0 \\
    	0 \\
    	3 \\
    	-4 \\
    	1 \\
    	0 \\
    \end{bmatrix}\!\!,
    \;
    \overline{\K}^{+}_{1,+} \!\!=\!\!\!
	\begin{bmatrix} 
    	0 \\
    	0 \\
    	-1 \\
    	0 \\
    	1 \\
    	0 \\
    	0 \\
    \end{bmatrix}\!\!,
    \;
    \overline{\K}^{+}_{2,+} \!\!=\!\!\!
	\begin{bmatrix} 
    	0 \\
    	1 \\
    	-4 \\
    	3 \\
    	0 \\
    	0 \\
    	0 \\
    \end{bmatrix}\!\!,
    \;
    \overline{\K}^{-}_{0,+} \!\!=\!\!\!
	\begin{bmatrix} 
    	1 \\
    	-4 \\
    	3 \\
    	0 \\
    	0 \\
    	0 \\
    	0 \\
    \end{bmatrix}\!\!,
    \;
    \overline{\K}^{-}_{1,+} \!\!=\!\!\!
	\begin{bmatrix} 
    	0 \\
    	1 \\
    	0 \\
    	-1 \\
    	0 \\
    	0 \\
    	0 \\
    \end{bmatrix}\!\!,
    \;
    \overline{\K}^{-}_{2,+} \!\!=\!\!\!
	\begin{bmatrix} 
    	0 \\
    	0 \\
    	3 \\
    	-4 \\
    	1 \\
    	0 \\
    	0 \\
    \end{bmatrix}\!\!, 
\end{equation*}
and similarly $\overline{\K}^{\pm}_{k,-}$ are given by
\begin{equation*}
	\overline{\K}^{+}_{0,-} \!\!=\!\!\!
	\begin{bmatrix}
    	0 \\
    	0 \\
    	0 \\
    	0 \\
    	3 \\
    	-4 \\
    	1 \\
    \end{bmatrix}\!\!,
    \;
    \overline{\K}^{+}_{1,-} \!\!=\!\!\!
	\begin{bmatrix} 
    	0 \\
    	0 \\
    	0 \\
    	-1 \\
    	0 \\
    	1 \\
    	0 \\
    \end{bmatrix}\!\!,
    \;
    \overline{\K}^{+}_{2,-} \!\!=\!\!\!
	\begin{bmatrix} 
    	0 \\
    	0 \\
    	1 \\
    	-4 \\
    	3 \\
    	0 \\
    	0 \\
    \end{bmatrix}\!\!,
    \;
    \overline{\K}^{-}_{0,-} \!\!=\!\!\!
	\begin{bmatrix} 
    	0 \\
    	1 \\
    	-4 \\
    	3 \\
    	0 \\
    	0 \\
    	0 \\
    \end{bmatrix}\!\!,
    \;
    \overline{\K}^{-}_{1,-} \!\!=\!\!\!
	\begin{bmatrix} 
    	0 \\
    	0 \\
    	1 \\
    	0 \\
    	-1 \\
    	0 \\
    	0 \\
    \end{bmatrix}\!\!,
    \;
    \overline{\K}^{-}_{2,-} \!\!=\!\!\!
	\begin{bmatrix} 
    	0 \\
    	0 \\
    	0 \\
    	3 \\
    	-4 \\
    	1 \\
    	0 \\
    \end{bmatrix}\!\!, 
\end{equation*}

%--------------------------------------------------------------
%		    Hadamard product and Frobenius inner product
%--------------------------------------------------------------
\section{Hadamard product}
\label{Appendix_Hadamard}
The Hadamard product is the matrix entry-wise product of two square matrices $A$ and $B$ of the same dimensions \emph{i.e.,}
\begin{equation*}
	(A \circ B)_{i,j} = A_{i,j} \cdot B_{i,j},
\end{equation*}
which we use to derive our numerical scheme.\\

%%%%%%%%%%%%%%%%%%%%%%%%%

\begin{acknowledgements}
The authors would like to acknowledge the EPSRC DTA EP/R513192/1 grant that supported this research. 
\end{acknowledgements}

%--------------------------------------------------------------
%	  						Bibliography
%--------------------------------------------------------------

\bibliographystyle{spmpsci}       % mathematics and physical sciences

\bibliography{ThesisReferences}   % name your BibTeX data base

\end{document}